\documentclass[11pt,a4paper,twoside,final]{scrartcl}
\usepackage{a4wide}
\usepackage{amsfonts}
\usepackage{amsmath}
\usepackage{amssymb}
\usepackage{amsthm}
\usepackage[numbers]{natbib}
\usepackage{booktabs}
\usepackage{algorithm,algpseudocode}
\usepackage{pgfplots}
\pgfplotsset{compat=newest}
\usepackage{graphicx}
\usepackage{color}
\usepackage{colortbl}
\usepackage{multirow}
\usepackage{boxedminipage}
\usepackage{enumerate}
\usepackage[caption=false,font=footnotesize]{subfig}
\usepackage{rotating}

\usepackage{hyperref}

\algnewcommand{\algorithmicgoto}{\textbf{go\;to}}%
\algnewcommand{\Goto}[1]{\algorithmicgoto~\ref{#1}}%

\newcommand{\N}{\ensuremath{\mathbb{N}}}

\newcommand{\T}{\ensuremath{\mathbb{T}}}

\newcommand{\Z}{\ensuremath{\mathbb{Z}}}

\newcommand{\R}{\ensuremath{\mathbb{R}}}
\newcommand{\A}{\ensuremath{\operatorname{A}}}
\newcommand{\C}{\ensuremath{\mathbb{C}}}
\newcommand{\nchoosek}[2]{\left(\begin{array}{c}#1\\#2\end{array}\right)}
\newcommand{\ii}{\textnormal{i}}
\newcommand{\e}{\textnormal{e}}

\newcommand{\ceil}[1]{\left\lceil#1\right\rceil}
\newcommand{\floor}[1]{\left\lfloor#1\right\rfloor}

\newcommand{\boldk}{{\ensuremath{\boldsymbol{k}}}}

\newcommand{\boldh}{{\ensuremath{\boldsymbol{h}}}}

\newcommand{\boldx}{{\ensuremath{\boldsymbol{x}}}}

\newcommand{\boldz}{{\ensuremath{\boldsymbol{z}}}}

\newcommand{\boldnu}{{\ensuremath{\boldsymbol{\nu}}}}

\newcommand{\boldgamma}{{\ensuremath{\boldsymbol{\gamma}}}}

\newcommand{\boldzero}{{\ensuremath{\boldsymbol{0}}}}

\newtheorem{theorem}{Theorem}[section]
\newtheorem{lemma}[theorem]{Lemma}
\newtheorem{remark}[theorem]{Remark}
\newtheorem{generalisation}[theorem]{Generalisation}
\newtheorem{definition}[theorem]{Definition}
\newtheorem{example}[theorem]{Example}
\newtheorem{corollary}[theorem]{Corollary}
\newtheorem{proposition}[theorem]{Proposition}

\newenvironment{Theorem}{\goodbreak \begin{theorem}\normalfont \slshape}{\end{theorem}}
\newenvironment{Lemma}{\goodbreak \begin{lemma}\normalfont \slshape}{\end{lemma}}
\newenvironment{Remark}{\goodbreak \begin{remark}\normalfont \rmfamily}{\bend\end{remark}}

\newenvironment{Corollary}{\goodbreak \begin{corollary}\normalfont \slshape}{\end{corollary}}

\makeatletter
\def\imod#1{\allowbreak\mkern10mu({\operator@font mod}\,\,#1)}
\makeatother

\numberwithin{equation}{section}
\numberwithin{table}{section}
\numberwithin{figure}{section}

\newcommand{\bend}{\hspace*{0ex} \hfill \hbox{\vrule height
    1.5ex\vbox{\hrule width 1.4ex \vskip 1.4ex\hrule  width 1.4ex}\vrule
    height 1.5ex}}

\long\def\symbolfootnote[#1]#2{\begingroup%
\def\thefootnote{\fnsymbol{footnote}}\footnote[#1]{#2}\endgroup}
\newcommand{\sspan}{\textnormal{span}}

\newcommand{\OO}[1]{\mathcal{O}\left(#1\right)}

\renewcommand{\mathbf}[1]{\ensuremath{\boldsymbol{#1}}}
\renewcommand{\textbf}[1]{{\ensuremath{\boldsymbol{#1}}}}
\renewcommand{\thefootnote}{\fnsymbol{footnote}}
\pgfplotscreateplotcyclelist{MR1LOF}{%
dashed, mark=diamond, mark options={solid, scale=1.5}\\%
dashed, mark=o, mark options={solid, scale=1.5}\\%
densely dotted, mark=triangle, mark options={solid, scale=1.5}\\%
densely dotted, mark=square, mark options={solid, scale=1.25}\\%
dashdotted,mark=x, mark options={solid, scale=1.5}\\%
dashdotted, mark=+, mark options={solid, scale=1.5}\\%
loosely dashed, mark=square, mark options={solid,scale=0.75}\\%
loosely dashdotted, mark=asterisk, mark options={solid,scale=0.5}, opacity=0.5\\%
dasdotdotted, every mark/.append style={solid},mark=star\\%
densely dashdotted,every mark/.append style={solid, fill=gray},mark=diamond*\\%
}%

\title{A fast probabilistic component--by--component construction of exactly integrating rank-1 lattices and applications}

\date{\today}
\author{
Lutz K\"ammerer\footnotemark[1]}

\hypersetup{pdfauthor=Lutz K\"ammerer,
          pdftitle={A fast probabilistic component--by--component construction of exactly integrating rank-1 lattices and applications},
          pdfsubject={},
          pdfcreator = {pdflatex},
          plainpages=false,
          pdfstartview=FitH,       
          pdfview=FitH,            
          pdfpagemode=UseOutlines, 
          bookmarksnumbered=true, 
          bookmarksopen=false,     
          bookmarksopenlevel=0,   
          colorlinks=true,       
          linkcolor=black,
          citecolor=black,
          urlcolor=black}

\newif\ifshowextendedpaperversion
\showextendedpaperversionfalse

\begin{document}

\maketitle

\begin{abstract}
\small
During the last decades, rank-1 lattices as sampling schemes for numerical integration and approximation became highly popular due to their universal applicability in high dimensional problems.
This popularity is essentially based on the pioneering observation that suitable rank-1 lattices can be constructed component--by--component (CBC).
However, precisely this necessary construction is the decisive point at which the practical applicability may fail.

For that reason, several more and more efficient CBC constructions were developed during the last decades.
On the one hand, there exist constructions that are based on minimizing some error functional, whereby this functional requires some dimension incremental structure.
On the other hand, there is the possibility to construct rank-1 lattices whose corresponding cubature rule, a rank-1 lattice rule, exactly integrates all elements within a space of multivariate trigonometric polynomials.

In this paper, we focus on the second approach, i.e., the exactness of rank-1 lattice rules.
The main contribution is the analysis of a probabilistic version of an already known algorithm that realizes a CBC construction of such rank-1 lattices.
It turns out that the computational effort of the known deterministic algorithm can be considerably improved in average by means of a simple randomization.
Moreover, we give a detailed analysis of the computational costs with respect to a certain failure probability, which then leads to the development of a probabilistic CBC algorithm.
In particular, the presented approach will be highly beneficial for the construction of so-called reconstructing rank-1 lattices, that are practically relevant for function approximation.
Subsequent to the rigorous analysis of the presented CBC algorithms, we present an algorithm that determines reconstructing rank-1 lattices of reasonable lattice sizes with high probability.
We provide estimates on the resulting lattice sizes and bounds on the occurring failure probability.
Furthermore, we discuss the computational complexity of the presented algorithm.

Various numerical tests illustrate the efficiency of the presented algorithms. 
Among others, we demonstrate how to exploit the efficiency of our algorithm even for the construction of
exactly integrating rank-1 lattices, provided that a certain property of the treated space of trigonometric polynomials is known.


\medskip

{\small%
}
\end{abstract}
\footnotetext[1]{
  Chemnitz University of Technology, Faculty of Mathematics, 09107 Chemnitz, Germany\\
  kaemmerer@mathematik.tu-chemnitz.de
}

\medskip

\ifshowextendedpaperversion
\tableofcontents
\newpage
\fi

\section{Introduction}
\begin{table}
\small
\setlength{\tabcolsep}{15pt}
\begin{tabular}{p{1cm}p{3.5cm}p{4cm}p{2.75cm}}
&\multicolumn{2}{c}{computational complexity for}\\
&\multicolumn{2}{c}{constructing rank-1 lattices that allow for}\\[.5em]
& exact integration & unique reconstruction & comments\\[.5em]
\midrule
\rule{0em}{2em}Alg.~\ref{alg:heuristic}
& $d\,\log^{1+\varepsilon}(d/\delta)\,|I|\log|I|$\newline
{\footnotesize (cf.~Remark~\ref{rem:heuristic_complex_int})}
& $d\,\log^{1+\varepsilon}(d/\delta)\,|I|\,\log^2|I|$\newline {\footnotesize (cf.~Remark~\ref{rem:heuristic_complex_reco})}&
probabilistic%
\newline $M<4M_{\operatorname{lb}}$\newline success w.h.p.\\
\rule{0em}{1.5em}\cite{KuoMiNoNu19} &$d\,|I|\,\log|I|$
&$d\,|I|^2\,\log|I|$&
non-probabilistic%
\newline $M<2M_{\operatorname{lb}}$\\[1.em]
\bottomrule
\end{tabular}
\caption{Computational complexities of most efficient component--by--component constructions of rank-1 lattices that fulfill properties \eqref{eq:exact_int_property} and \eqref{eq:reco_property}, respectively,  in comparison. The parameter $\delta\in(0,1)$ is an upper bound on the failure probability of Algorithm~\ref{alg:heuristic}, $\varepsilon>0$ arbitrarily small, and the subset relation $I\subset[|I|,|I|]^d$ is assumed.}\label{tab:intro_complexities}
\end{table}

The advantageous general approximation properties, cf.\,e.g.~\cite{KuSlWo08,KuWaWo09,KaPoVo13,ByKaUlVo16}, of rank-1 lattice sampling strategies in combination with the possibility of component--by--component (CBC) constructions of suitable rank-1 lattices have led to intensive studies on rank-1 lattices during the last decades.
In particular, the beneficial structure of rank-1 lattices allows for using one dimensional fast Fourier transforms in order to efficiently reconstruct signals
from sampling points, cf.~\cite{LiHi02,MuSo11}.

However, suitable rank-1 lattices need to be determined in order to guarantee for excellent approximation properties of this signal reconstruction.
For a given rank-1 lattice size $M$, i.e., the number of sampling points one spends for the reconstruction of the signal, the usage of CBC strategies for determining suitably generating vectors $\boldz\in\{0,\ldots,M-1\}^d$ heavily reduces the size of the search space from $M^d$ to $dM$, whereby this restricted search space does not cause
significant disadvantages with respect to the achievable approximation errors in specific situation.

The runtime of relevant CBC algorithms proved to be polynomially in the spatial dimension $d$ and the lattice size $M$, \cite{SlRe02, CoNu07, Kae2013}.
A recently developed, non-probabilistic algorithm \cite{KuoMiNoNu19} for the construction of exactly integrating rank\mbox{-}1 lattices has
a computational complexity that is linear in the dimension $d$, the lattice size $M$, and the cardinality $|I|$ of the frequency set $I$, $I\subset\Z^d$, $|I|<\infty$, which is highly efficient
in this setting. In this context, exactly integrating means that the quasi Monte Carlo method based on the samples along the rank-1 lattice is a cubature rule that exactly integrates all signals within the linear space $\Pi_I:=\sspan\{\e^{2\pi\ii\boldk\cdot\circ}\colon\boldk\in I\}$, cf.\ Section~\ref{sec:pre:CBC}.

A slightly modified approach can be used to construct rank-1 lattices that allow for an exact reconstruction of all signals from linear spaces $\Pi_I$ of multivariate trigonometric polynomials, cf. Section~\ref{sec:pre:reco} for details.
We call such rank-1 lattices \emph{reconstructing rank-1 lattices for signals in $\Pi_I$} or -- in short -- \emph{reconstructing rank-1 lattices for the frequency set $I$}.
The crucial advantages of such reconstructing rank-1 lattices are their applicability
for approximating functions that are not necessarily in $\Pi_I$, the corresponding highly efficient algorithms using a permutation and a single one-dimensional fast Fourier transform, and the associated known error guarantees, cf.~\cite{Tem86, KaPoVo13, ByKaUlVo16, KuoMiNoNu19}.

In this paper, we develop a probabilistic approach that can be used for the construction of exactly integrating rank-1 lattices and for the construction of reconstructing rank-1 lattices as well.
The crucial innovation of the presented algorithm is that the computational costs do not depend on the lattice size $M$ anymore and that they are almost linear in the cardinality of the frequency set $I$ under consideration. Roughly speaking, the
computational costs depend almost linear on the number $|I|$ of multivariate monomials which one can exactly integrate simultaneously using the sampling points of the resulting rank-1 lattice or, in the case of signal reconstruction, on the number $|I|$ of multivariate monomials that builds up the space of multivariate trigonometric polynomials under consideration from which all signals can be reconstructed uniquely using the sampling points of the resulting rank-1 lattice.
In particular for reconstructing rank-1 lattices, these two numbers $M$ and $|I|$ may necessarily differ widely, i.e., up to $|I|\lesssim \sqrt{M}$, cf.~\cite{KaKuPo10,ByKaUlVo16}.
Precisely in these situations, the presented probabilistic algorithm will reduce the computational costs for the construction of reconstructing rank-1 lattices significantly compared to known algorithms.

Table~\ref{tab:intro_complexities} presents best known upper bounds on the computational complexities of state-of-the-art non-probabilistic algorithms, cf.~\cite{KuoMiNoNu19}, and Algorithm~\ref{alg:heuristic}, that, or rather its analysis, is one important contribution of this paper.
It turns out that the crucial advantage of the presented probabilistic approach is in the construction of reconstructing rank-1 lattices. On the one hand, the new approach suffers from a certain, small failure probability and the upper bound on the resulting lattice sizes $M$ is doubled compared to known non-probabilistic approaches. On the other hand, the computational complexity is highly reduced since it depends almost linear on the cardinality $|I|$ of the considered frequency set $I$ whereas the computational complexity of the non-probabilistic approach depends quadratically on this cardinality $|I|$.

Besides the specification of all used algorithms, the main contribution of this work is the analysis of the computational complexities and, in particular, the success probabilities in full detail, cf. Sections~\ref{sec:main} to~\ref{sec:prob_M}.
More specifically, we analyze the CBC approach for determining a suitable generating vector $\boldz\in\Z^d$ for fix and large enough lattice size $M$ in Section~\ref{sec:main}. In Section~\ref{sec:check_exactness}
we discuss the difference between the presented algorithms when determining exactly integrating rank-1 lattices and reconstructing rank-1 lattices. Furthermore, we present Algorithm~\ref{alg:heuristic} in Section~\ref{sec:prob_M} and determine its resulting rank-1 lattice size in conjunction with associated success rate estimates.
Various numerical tests demonstrate the practical relevance of the newly developed approaches in Section~\ref{sec:numerics}.

\section{Prerequisites}

First, we collect some notations and basic facts about rank-1 lattices.
For $\boldz\in\Z^d$ and $M\in\N$, we call the set of sampling nodes
$$
\Lambda(\boldz,M):=\left\{\boldx_j:=\frac{j}{M}\boldz\bmod{1}\colon j=0,\ldots,M-1\right\}\subset\T^d\simeq[0,1)^d
$$
a rank-1 lattice of size $M$ with generating vector $\boldz$.
The main goal of this paper is the construction of so-called exactly integrating lattice rules with respect to a given frequency set $I\subset\Z^d$, $|I|<\infty$.
In short words, this means that 
each $d$-variate trigonometric polynomial $p\in\Pi_I:=\operatorname{span}\{\e^{2\pi\ii\boldk\cdot\circ}\colon\boldk\in I\}$ of the form
\begin{equation}
p(\boldx):=\sum_{\boldk\in I}\hat{p}_\boldk\e^{2\pi\ii\boldk\cdot\boldx},
\label{eq:def_trig_poly}
\end{equation}
$\left(\hat{p}_\boldk\right)_{\boldk\in I}\in \C^{|I|}$, 
can be exactly integrated using the cubature rule
\begin{equation}
Q_{\boldz,M}(p):=M^{-1}\sum_{j=0}^{M-1}p\left(\frac{j}{M}\boldz\right).
\label{eq:qmc_rule}
\end{equation}
Since the set of sampling nodes used by $Q_{\boldz,M}$ is a rank-1 lattice, 
$Q_{\boldz,M}$ is called (rank-1) lattice rule.
In detail, the lattice rule $Q_{\boldz,M}$ exactly integrates all polynomials $p$ with frequency support in $I$ whenever
$Q_{\boldz,M}(p)=\hat{p}_\boldzero$
holds for all $p\in\Pi_I$, which
is equivalent to
$$
Q_{\boldz,M}(\e^{2\pi\ii\boldk\cdot\circ})=\begin{cases}
1, & \boldk=\boldzero,\\
0, & \boldk\in I\setminus\{\boldzero\}\,.
\end{cases}
$$
As a consequence of \eqref{eq:qmc_rule}, we obtain
$$Q_{\boldz,M}(\e^{2\pi\ii\boldk\cdot\circ})=
M^{-1}\sum_{j=0}^{M-1}\left(\e^{2\pi\ii\frac{\boldk\cdot\boldz}{M}}\right)^j
=\begin{cases}
1, & \frac{\boldk\cdot\boldz}{M}\in\Z,\\
0, & \text{otherwise}\,,
\end{cases}
$$
which leads directly to the observation that $Q_{\boldz,M}$ exactly integrates
all $p\in\Pi_I$ iff the generating vector $\boldz\in\Z^d$ and the lattice size $M\in\N$ fulfill
\begin{equation}
\boldk\cdot\boldz\not\equiv 0\imod{M}\text{ for all }\boldk\in I\setminus\{\boldzero\}.
\label{eq:exact_int_property}
\end{equation}
In cases where $\boldz\in\Z^d$ and $M\in\N$ are chosen such that \eqref{eq:exact_int_property}
holds for a given frequency set $I$, we say that $\Lambda(\boldz,M)$ fulfills the exact integration property for $I$.

One highly interesting question is how to choose a suitable generating vector $\boldz$ and an adequate lattice size $M$ for a given frequency set $I$. In particular, for too small $M$, we cannot find any $\boldz\in\Z^d$ such that
\eqref{eq:exact_int_property} holds. For instance, for any $I$ with $|I|\ge 2$ and fixed $M=1$ there exists no generating vector $\boldz\in\Z^d$ such that \eqref{eq:exact_int_property} is fulfilled. Moreover, the lattice size $M$ is (an upper bound on) the number of sampling values one uses in order to apply the lattice rule $Q_{\boldz,M}$.
Accordingly, one is interested in lattice sizes $M$ as small as possible.

We define 
$\max(I):=\max\{\|\boldk\|_\infty\colon\boldk\in I\}$
which is the largest component of all frequency vectors within $I$ in absolute value.
Clearly, using $M=(\max(I)+1)^d$ and the generating vector $\boldz=(z_1,\ldots,z_d)^\top$, $z_t:=(\max(I)+1)^{t-1}$, $t=1,\ldots,d$,
yields a rank-1 lattice $\Lambda(\boldz,M)$ that fulfills the exact integration property \eqref{eq:exact_int_property} for~$I$ in general.
However, this lattice size $M$ is far too big -- at least in cases where $I$ is not a cartesian grid.

Nevertheless, the mere existence of the last mentioned exactly integrating rank-1 lattice for~$I$
implies that the finite number of possible rank-1 lattice sizes not greater than $M=(\max(I)+1)^d$
contains a global minimum such that 
\eqref{eq:exact_int_property} holds true for some generating vector~$\boldz$. However, in general a search for the smallest possible lattice size~$M$ and a corresponding generating vector~$\boldz$ is practically impossible even for moderate dimensions~$d$ and frequency sets~$I$ of small cardinality~$|I|$ due to the huge search space.
Even for known suitable lattice size~$M$, the (brute force) search for an appropriate generating vector~$\boldz$ is hard due to the huge search space consisting of~$M^d$ elements.

\subsection{Component--by--component construction of exactly integrating rank-1 lattices}\label{sec:pre:CBC}

Nevertheless, there exists a general component--by--component (CBC) strategy for determining suitable components $z_j$, $j=1,\ldots,d$, of the generating vector $\boldz$ one after another.
Clearly, this strategy may not allow for finding a generating vector $\boldz$ for the global minimum of all possible lattice sizes $M$, but there exist simple and reasonable estimates on primes $M$ such that this CBC strategy will successfully determine a corresponding generating vector such that 
\eqref{eq:exact_int_property} holds true.

We present a short sketch on this strategy, cf.~\cite{KuoMiNoNu19} for details.
Let $M$ be an integer such that
$$
k_1\not\equiv 0\imod{M}\text{ for all } \boldk\in I\setminus\{\boldh\in I\colon h_1=0\} 
$$
and set $z_1:=1$. In addition, for any $s\in{2,\ldots,d}$, let the integer $M$ fulfill
\begin{equation}
\begin{split}
&\Bigg|\{y\in\{0,\ldots,M-1\}\colon\\
&\quad -k_sy\equiv\sum_{j=1}^{s-1}k_jz_j\imod{M}\text{ for any }\boldk\in I\setminus\{\boldh\in I\colon h_1=\ldots=h_s=0\}\}\Bigg|<M \label{eq:cbc_M_estimate_detail}
\end{split}
\end{equation}
for any admissible choice of $z_2,\ldots,z_{s-1}\in\{0,\ldots,M-1\}$, i.e., for any choice of $z_2,\ldots,z_{s-1}\in\{0,\ldots,M-1\}$ such that
\begin{equation}
\begin{split}
\sum_{j=1}^{s-1}&k_jz_j\not\equiv 0\imod{M}\\
&\text{ for all }(k_1,\ldots,k_{s-1})^\top\in\left\{(h_1,\ldots,h_{s-1})^\top\colon\boldh\in I, \sum_{j=1}^{s-1}|h_j|\neq 0\right\}
\end{split}
\label{eq:admissible_zs}
\end{equation}
is valid.
Obviously, even for this strategy it is practically impossible to determine the smallest possible $M$ that fulfills \eqref{eq:cbc_M_estimate_detail} for any admissible choice 
of the first $s-1$ components of the generating vector since the number of admissible choices may be extremely large.

Nevertheless, the CBC strategy succeeds with suitably large primes $M$, which one estimates using the just mentioned approach. Usually, recent research results determine bounds such that for each prime $M$ larger than this bound this CBC strategy will succeed.
This paper is not concerned about estimates on these lower bounds.
However, we just define an abstract lower bound $M_{\operatorname{lb}}\in\N$
for our theoretical considerations in the following way.
For each $M\in\N$ and each $s\in\{2,\ldots,d\}$ we define the sets
\begin{equation}
\begin{split}
&Y^{(s)}(z_1,\ldots,z_{s-1},M,I):=\Bigg\{y\in\{0,\ldots,M-1\}\colon\\
&\quad -k_sy\equiv\sum_{j=1}^{s-1}k_jz_j\imod{M}\text{ for any }\boldk\in I
\setminus\{\boldh\in I\colon h_1=\ldots=h_s=0\}\Bigg\}
\end{split}
\end{equation}
and
\begin{equation}
\begin{split}
&Z^{(s-1)}(M,I):=\bigg\{(z_1,\ldots,z_{s-1})^\top\in\{0,\ldots,M-1\}^{s-1}\colon\\
&\qquad\qquad\qquad\qquad\qquad\text{\eqref{eq:admissible_zs} holds for } (z_1,\ldots,z_{s-1})^\top, M, \text{ and } I\bigg\}\,.
\end{split}
\end{equation}
Then, we define
\begin{equation*}
N'(M,I)=\max\limits_{s=2,\ldots, d}\,\max\limits_{(z_1,\ldots,z_{s-1})^\top\in Z^{(s-1)}(M,I)}\left|Y^{(s)}(z_1,\ldots,z_{s-1},M,I)\right|
\end{equation*}
and we take into account $s=1$
\begin{equation}
N(M,I):=
\begin{cases}
M & k_1\equiv 0\imod{M}\text{ for any } \boldk\in \{\boldh\in I\colon h_1 \neq 0\} \\
N'(M,I)
&\text{otherwise.}
\end{cases}\label{eq:def_NMI}
\end{equation}

Certainly, the bound $N(M,I)\le M$ holds true and values $N(M,I)$ less than $M$ guarantees
success for CBC approaches that checks all the possiblities for choosing the components
of the generating vector in each CBC step.

\begin{sloppypar}
Now, the usual approach is to restrict $M$ to prime numbers larger than $\max\{\|\boldk\|_\infty\colon \boldk\in I\}=:\max(I)$, which has the crucial advantage that the cardinality of the sets of failing candidates
$Y^{(s)}(z_1,\ldots,z_{s-1},M,I)$ can be simply estimated by
$|Y^{(s)}(z_1,\ldots,z_{s-1},M,I)|\le |I|$, cf.~\cite{KuoMiNoNu19} for more detailed estimates.
Taking into account that $M\ge 1+\max\{|k_1|\colon \boldk\in I\}$ we observe $N(M,I)=N'(M,I)$ and simply estimate
$N(M,I)\le |I|$ for all large enough primes $M$.
\end{sloppypar}

Due to the fact that $N(M,I)$ is uniformly bounded for all large enough primes $M$, we can find a minimal prime $P_{\min}(I)\in\N$
such that
$$
P_{\min}(I)> \max\left(\max_{\substack{M \text{ prime}\\M\ge P_{\min}(I)}}N(M,I),\; \max(I)\right)\,.
$$
Clearly, $P_{\min}(I)$ is at most the smallest prime number larger
than $\max\left(|I|,\,\max(I)\right)$.
Constructing this abstract prime number $P_{\min}(I)$, we can define
$$M_{\operatorname{lb}}=M_{\operatorname{lb}}(I):=P_{\min}(I)-1\,.$$
Taking \eqref{eq:def_NMI} into account, one observes that
for each prime number $P>M_{\operatorname{lb}}(I)$ and each $s=2,\ldots,d$, we can estimate the cardinalities 
$|Y^{(s)}(z_1,\ldots,z_{s-1},P,I)|\le M_{\operatorname{lb}}(I)$.

We read this from another point of view.
In all cases where the lattice size $M$ is chosen as prime with $M\ge P_{\min}(I)$,
we can find at least $M-N(I,M) \ge M-M_{\operatorname{lb}}(I)$ different appropriate choices for
the component $z_s$ of the generating vector $\boldz$ in the current $s$th CBC step.
The last observation holds true for all $s=2,\ldots,d$, due to the definition of
$M_{\operatorname{lb}}(I)$.
At this point, we would like to point out that the number of appropriate choices of the components
of the vectors can significantly exceed the number $M-M_{\operatorname{lb}}(I)$ due to the fact that
$M_{\operatorname{lb}}(I)$ is a very universal upper bound on the cardinalities $|Y^{(s)}(z_1,\ldots,z_{s-1},M,I)|$
that determines the amount of integers within $\{0,\ldots,M-1\}$ that are no appropriate choice of the component $z_s$
in the case, where $M$, $s$, and also the first $s-1$ components of the generating vector $\boldz$ are already fixed.

\subsection{Reconstructing rank-1 lattices}\label{sec:pre:reco}

Up to now, we discussed the exact integration property of rank-1 lattices.
In the following, we focus on the reconstruction property of rank-1 lattices.
To this end, we consider the trigonometric polynomial $p$ as in \eqref{eq:def_trig_poly}.
Since $\left(\e^{2\pi\ii\boldk\cdot\circ}\right)_{\boldk\in\Z^d}$ is an orthonormal system in $L_2(\T^d)$, we
can reconstruct the active Fourier coefficients $\hat{p}_\boldk$, $\boldk\in I$, by solving
the integrals
\begin{equation}
\hat{p}_\boldk=\int_{\T^d}p(\boldx)\e^{-2\pi\ii\boldk\cdot\boldx}\mathrm{d}\boldx, \quad \boldk\in I\,.
\label{eq:int_fc_p}
\end{equation}
Similar as above, we can treat these integrals numerically using the lattice rule
$$
Q_{\boldz,M}\left(p(\circ)\e^{-2\pi\ii\boldk\cdot\circ}\right):=M^{-1}\sum_{j=0}^{M-1}p\left(\frac{j}{M}\boldz\right)
\e^{-2\pi\ii\frac{j}{M}\boldk\cdot\boldz}.
$$
We call the rank-1 lattice $\Lambda(\boldz,M)$ a reconstructing rank-1 lattice for the frequency set $I$ iff 
$Q_{\boldz,M}\left(p(\circ)\e^{-2\pi\ii\boldk\cdot\circ}\right)=\hat{p}_\boldk$ holds for all $\boldk\in I$ and all $p\in\Pi_I$. In other words, the lattice rule $Q_{\boldz,M}$ allows for the exact integration of all $|I|$ integrals in
\eqref{eq:int_fc_p}.
Consequently, we require
\begin{equation}
\begin{split}
\hat{p}_\boldk&=Q_{\boldz,M}\left(p(\circ)\e^{-2\pi\ii\boldk\cdot\circ}\right)=M^{-1}\sum_{j=0}^{M-1}\sum_{\boldh\in I}\hat{p}_\boldh
\e^{2\pi\ii\frac{j}{M}(\boldh-\boldk)\cdot\boldz}\\
&=
\sum_{\boldh\in I}\hat{p}_\boldh M^{-1}\sum_{j=0}^{M-1}
\left(\e^{2\pi\ii\frac{(\boldh-\boldk)\cdot\boldz}{M}}\right)^j
=\sum_{\substack{\boldh\in I\\\frac{(\boldh-\boldk)\cdot\boldz}{M}\in \Z}}\hat{p}_\boldh
\end{split}\label{eq:aliasing_formula}
\end{equation}
for all $\boldk\in I$. Thus equation \eqref{eq:aliasing_formula} yields that
a reconstructing rank-1 lattice $\Lambda(\boldz,M)$ necessarily (and sufficiently) fulfills
\begin{equation}
(\boldh-\boldk)\cdot\boldz\not\equiv 0\imod{M}\text{ for all }\boldh,\boldk\in I, \boldh\neq\boldk
\label{eq:reco_diff_property}
\end{equation}
or, equivalently,
\begin{equation}
\boldh\cdot\boldz\not\equiv \boldk\cdot\boldz\imod{M}\text{ for all }\boldh,\boldk\in I, \boldh\neq\boldk\,.
\label{eq:reco_property}
\end{equation}
We will call \eqref{eq:reco_property} reconstruction property of $\Lambda(\boldz,M)$ for the frequency set $I$.
At this point, we would like to emphasize the close connection between exactly integrating rank-1 lattices and reconstructing rank-1 lattices.
To this end, we define the difference set
$D(I):=\{\boldh-\boldk\colon\boldk,\boldh\in I\}$
of the frequency set $I$, which leads to the observation that \eqref{eq:reco_diff_property}
is equivalent to the exact integration property of $\Lambda(\boldz,M)$ for the frequency set $D(I)$, cf.~\eqref{eq:exact_int_property}.
In other words, a rank-1 lattice that satisfies the reconstruction property for $I$ equivalently fulfills the
exact integration property for the difference set $D(I)$ of $I$.

As a consequence, one can apply a CBC strategy to $D(I)$ as explained in Section~\ref{sec:pre:CBC}
in order to construct reconstructing rank-1 lattices for the frequency set $I$.
At this point we define
$$
N_I:=\max(D(I))=\max_{t=1,\ldots,d}\left(\max\{k_t\colon \boldk\in I\}-\min\{k_t\colon\boldk\in I\}\right)\le 2\max(I)
$$
the \emph{expansion of a frequency set $I$} for later usage.

\section{Main algorithms}\label{sec:main}

\begin{algorithm}[th]
\small
\caption{Probabilistic deterministic CBC construction of exactly integrating rank-1 lattices}\label{alg:construct_rand_reco_r1l_I_basic}
  \begin{tabular}{p{2cm}p{5cm}p{7.3cm}}
    Input: 	& $I\subset[-N,N]^d$ 	& frequency set\\
    		& $M$	& lattice size $M\ge 2N+1$\\
    		& $T\le M$	& switching point of algorithm
  \end{tabular}
		
  \begin{algorithmic}[1]
	\State $z_1=1$;
	\For{$\ell=2:d$}
		\State choose a set $\mathcal{T}$ of $T$ different elements from $\{0,\ldots,M-1\}$, uniformly at random\label{alg:construct_rand_reco_r1l_I_basic:choose_T}
		\State shuffle the elements of $\mathcal{T}$ and save $(z_\ell^{(r)})_{r=1}^T$, i.e., $\{z_\ell^{(r)}, r=1,\ldots,T\}=\mathcal{T}$\label{alg:construct_rand_reco_r1l_I_basic:shuffle}
		\For{$r=1:T$}
			\If{{\texttt{check\_exactness}}($I_\ell$, $(z_1,\ldots,z_{\ell-1},z_\ell^{(r)})^\top$,$M$)$==$\texttt{true}}\label{alg:construct_rand_reco_r1l_I_basic:check_reco_property}
				\State $z_\ell=z_\ell^{(r)}$\quad and \quad \Goto{alg:construct_rand_reco_r1l_I_basic:gotopoint}\label{alg:construct_rand_reco_r1l_I_basic:break1}
			\ElsIf{$r==T$}\label{alg:construct_rand_reco_r1l_I_basic:r_larger_T}
				\State construct random permutation array $\left(z_\ell^{(r)}\right)_{r=T+1}^M$ of %
				 $\{0,\ldots,M-1\}\setminus \mathcal{T}$\label{alg:construct_rand_reco_r1l_I_basic:perm_M}
				\For{$r=T+1:M$}
					\If{{\texttt{check\_exactness}}($I_\ell$, $(z_1,\ldots,z_{\ell-1},z_\ell^{(r)})^\top$,$M$)$==$\texttt{true}}%
						\State $z_\ell=z_\ell^{(r)}$\quad and \quad \Goto{alg:construct_rand_reco_r1l_I_basic:gotopoint}
					\ElsIf{$r==M$}
						\State raise \texttt{ERROR}
					\EndIf
				\EndFor\label{alg:construct_rand_reco_r1l_I_basic:endinnerfor}
			\EndIf
		\EndFor\label{alg:construct_rand_reco_r1l_I_basic:gotopoint}
	\EndFor
  \end{algorithmic}
  \begin{tabular}{p{2cm}p{5cm}p{7.3cm}}
    Output: & $\boldz\in\{0,\ldots,M-1\}^d\subset\Z^d$ & generating vector of rank\mbox{-}1 lattice such that\\
	    & $\Lambda(\boldz,M)$ &  is an exactly integrating rank\mbox{-}1 lattice for $I$
\\    \cmidrule{1-3}
	Complexity:&
    \multicolumn{2}{l}{$\OO{d\,T\log T+ d\,\log(d/\delta)\, h(I,M)}$,}\\
    &\multicolumn{2}{l}{\qquad with probability $1-\delta$ for $M\ge c M_{\operatorname{lb}}$ and $T\ge\log_c((d-1)/\delta)$}\\
    &\multicolumn{2}{l}{$\OO{d\,T\log T+ d\, h(I,M)}$,}\\
    &\multicolumn{2}{l}{\qquad in expectation for $M\ge c M_{\operatorname{lb}}$ and $T\ge\log_c M$}\\
    &\multicolumn{2}{l}{$\OO{d\,M\,h(I,M)}$,}\\
    &\multicolumn{2}{l}{\qquad in worst case, assuming that $h(I,M)\gtrsim \log{M}$}
    \\\cmidrule{1-3}
    \multicolumn{3}{r}{\footnotesize$h(I,M)$ is an upper bound on the computational complexity of a single call of}\\ \multicolumn{3}{r}{\footnotesize\texttt{check\_exactness}($I_\ell$, $(z_1,\ldots,z_{\ell-1},y)^\top$,$M$) for fixed $\ell$ and arbitrary, fixed $y\in\{0,\ldots,M-1\}$ }
     \end{tabular}
\end{algorithm}

First we start with a fixed and suitably large lattice size $M$ and develop an appropriate strategy
that determines a generating vector such that $\Lambda(\boldz,M)$ is an exactly integrating rank\mbox{-}1 lattice for a given frequency
set $I$. Certainly, the lattice size $M$ depends on the frequency set $I$, cf.\ Section~\ref{sec:pre:reco}.
Since we will benefit from randomly checking component candidates for the generating vector $\boldz\in\{0,\ldots, M-1\}^d$,
we need a suitable strategy for fixing a random sequence of these candidates in order to construct a 
probabilistic deterministic algorithm. However, generating the whole sequence of $M$ candidates in advance
results in runtimes that are linear in $M$, which we need to avoid in order to obtain
runtimes of our algorithm that may be sublinear in $M$.
For that reason, we suggest a two-step strategy for determining a random permutation of the elements in $\{0,\ldots, M-1\}$.

\begin{Lemma}\label{lem:two_step_perm}
Let $T\in\{1,\ldots,M\}$ and we construct a permuted vector $v$ of all $\{0,\ldots,M-1\}$ numbers using
the following approach:

First we take a set $\mathcal{T}$ of $T$ distinct elements in $\{1,\ldots,M-1\}$ uniformly at random. Then, we shuffle the set $\mathcal{T}$ uniformly at random in order to obtain the first $T$ elements $v_1,\ldots,v_T$ of $v\in\N_0^M$. Subsequently, we determine $v_{T+1},\ldots,v_M$ as a random permutation of 
$\{0,\ldots,M-1\}\setminus\mathcal{T}$.
Then $v$ is a permutation of $\{0,\ldots,M-1\}$, and each permutation of $\{0,\ldots,M-1\}$ is obtained with probability $1/M!$.
\end{Lemma}
\begin{proof}
Clearly, each permutation may be reached  due to the construction.
A fixed set $\mathcal{T}$ is obtained with probability $\frac{T!(M-T)!}{M!}$. Each permutation of the set $T$ has a probability of $1/T!$ and each permutation of $\{0,\ldots,M-1\}\setminus\mathcal{T}$ has a probability of $1/(M-T)!$. Clearly shuffling $T$ or $M-T$ elements is independent from the concrete set, which yields a probability of the fixed vector $v$ that is $1/M!$.
\end{proof}

\begin{sloppypar}
As mentioned above, we investigate a probabilistic deterministic version of the brute force CBC strategy, cf.~\cite{KuoMiNoNu19}, using randomization.
We present the specific approach as Algorithm~\ref{alg:construct_rand_reco_r1l_I_basic}.
Randomly selecting permutations as suggested in Lemma~\ref{lem:two_step_perm} yields heavily improved arithmetic complexities with high probability and in average as well
when we assume that the lattice size $M$ is large enough and the parameter $T$ is chosen appropriately.
\end{sloppypar}

\begin{Theorem}\label{thm:deterministic}
Let $1<c\in\R$, $I\subset\Z^d$, $\delta>0$, and $M$ being a prime number satisfying $M\ge cM_{\operatorname{lb}}$.
Then,
\begin{itemize}
\item Algorithm~\ref{alg:construct_rand_reco_r1l_I_basic} determines an exactly integrating rank-1 lattice for the frequency set $I$ for arbitrarily chosen parameter $T\in\{1,\ldots,M\}$,
\item choosing $T\in\N$, $T\ge\frac{\log(d-1)-\log{\delta}}{\log{c}}$, Algorithm~\ref{alg:construct_rand_reco_r1l_I_basic} determines an exactly integrating rank-1 lattice for the frequency set $I$ checking a number of at most $\left\lceil\frac{\log(d-1)-\log{\delta}}{\log{c}}\right\rceil$ component candidates $z_\ell^{(r)}$ in each CBC step with probability at least $1-\delta$, which immediately yields an arithmetic complexity of Algorithm~\ref{alg:construct_rand_reco_r1l_I_basic} in
$$
\OO{d\,T\log T+ d\,\log(d/\delta)\, h(I,M)}
$$
with probability at least $1-\delta$,
\item for arbitrary $T\le M$, we achieve an expected arithmetic complexity of Algorithm~\ref{alg:construct_rand_reco_r1l_I_basic} in
$\OO{d(T\log T+c^{-T}M+h(I,M)}$,
\item assuming $h(I,M)\gtrsim \log{M}$, the worst case arithmetic complexity of Algorithm~\ref{alg:construct_rand_reco_r1l_I_basic} can be bounded by terms in $\OO{d M h(I,M)}$.
\end{itemize}
\end{Theorem}

\begin{proof}
Since the suggested algorithm realizes the usual brute force CBC search, cf.\ \cite[Rem.~24]{KuoMiNoNu19}, and the lattice size $M$ guarantees the existence of at least one exactly integrating rank-1 lattice for the frequency set $I$, which can be found using a CBC construction, Algorithm~\ref{alg:construct_rand_reco_r1l_I_basic} will determine one exactly integrating rank-1 lattice.

In order to show the other assertions we consider each CBC step separately, i.e., we fix $\ell\in\{2,\ldots,d\}$. Clearly, there may be elements $y\in\{0,\ldots,M-1\}$ such that 
\[
{\texttt{check\_exactness}}(I_\ell, (z_1,\ldots,z_{\ell-1},y)^\top,M)==\texttt{false},
\]
which we call failing component candidates. The number of these failing candidates is bounded from above by $M_{\operatorname{lb}}$, cf.\ Section~\ref{sec:pre:CBC}.
We define the sets of failing candidates in iteration $r$ by
\begin{align*}
\mathcal{Z}_\ell^r&:=\{y\in\{0,\ldots,M-1\}\setminus\{z_\ell^{(1)},\ldots,z_\ell^{(r-1)}\}\colon
\\
&\hspace*{5em}
{\texttt{check\_exactness}}(I_\ell, (z_1,\ldots,z_{\ell-1},y)^\top,M)==\texttt{false}\}\,.
\end{align*}
We know that in each CBC step, we have at most a number of $M_{\textnormal{lb}}$ of the potential $z_\ell$'s in the set
$\mathcal{Z}_\ell^1$ and we need to consider iteration $r$ iff $z_\ell^{(1)},\ldots,z_\ell^{(r-1)}$ are all failing candidates and each of them does not appear anymore.
Otherwise, $z_\ell^{(1)},\ldots,z_\ell^{(r-1)}$ are not all failing candidates, which means that
at least one of them is a non-failing, thus working, candidate and we can continue with the next CBC step.

Accordingly, we need to consider the set $\mathcal{Z}_\ell^r$ iff 
the set $\left\{z_\ell^{(1)},\ldots,z_\ell^{(r-1)}\right\}$ contains exactly $r-1$ of failing $z_\ell$'s which means that there are at most $M_{\operatorname{lb}}-r+1$ failing candidates left in the set 
$\mathcal{Z}_\ell^r$, in formula we have
$$
|\mathcal{Z}_\ell^r|\le \max(0,M_{\textnormal{lb}}-r+1)\le M/c\,.
$$
Thus, when reaching iteration $r$ we obtain
\begin{align*}
\mathbb{P}&\left(\left\{z_\ell^{(r)} \textnormal{ is a failing candidate}\right\}\bigg|\left\{z_\ell^{(1)},\ldots, z_\ell^{(r-1)} \textnormal{ are all failing candidates}\right\}\right)\\
&\quad=\mathbb{P}\left(\left\{z_\ell^{(r)}\in\mathcal{Z}_\ell^{r}\right\}\bigg|\bigcap_{j=1}^{r-1}\left\{z_\ell^{(j)}\in\mathcal{Z}_\ell^j\right\}\right)\le\frac{|\mathcal{Z}_\ell^r|}{M}\le \frac{M_{\textnormal{lb}}-r+1}{M}\le1/c
\end{align*}
and in particular $|\mathcal{Z}_\ell^{M_{\textnormal{lb}}+1}|=0$
which yields for arbitrary $r\le M_{\textnormal{lb}}$
\begin{equation}
\mathbb{P}\left(\bigcap_{j=1}^r\left\{z_\ell^{(j)}\in\mathcal{Z}_\ell^j\right\}\right)\le\prod_{j=1}^r \frac{M_{\textnormal{lb}}-j+1}{M}\le c^{-r}
\label{eq:prob_rth_iteration_fails}
\end{equation}
and 
\begin{equation}
\mathbb{P}\left(\bigcap_{j=1}^{M_{\textnormal{lb}}+1}\left\{z_\ell^{(j)}\in\mathcal{Z}_\ell^j\right\}\right)=0\,.
\label{eq:prob_zero}
\end{equation}

For fixed $T\in\N$, we use the union bound and estimate the probability that Algorithm~\ref{alg:construct_rand_reco_r1l_I_basic} checks in at least one of the $d-1$ CBC steps more than $T$
candidates $z_\ell^{(j)}$
\begin{equation}\label{eq:union_bound_basic}
\mathbb{P}\left(\bigcup_{\ell=2}^d\bigcap_{j=1}^T\left\{z_\ell^j\in\mathcal{Z}_\ell^j\right\}\right)\le(d-1)\prod_{j=1}^T \frac{M_{\textnormal{lb}}-j+1}{M}\le (d-1)c^{-T}
\end{equation}
Using any $T\ge\frac{\log{(d-1)}-\log{\delta}}{\log{c}}$ in \eqref{eq:union_bound_basic}, we observe $\delta$ as upper bound. Thus the assertion on the success probability holds.

Due to the last considerations, we observe that each CBC step will succeed with probability $1-\frac{\delta}{d-1}$ at latest in its $r_\delta$'s, $r_\delta=\ceil{\log_c\{\frac{d-1}{\delta}\}}$, iteration.
Due to the fact that $r_\delta\le T$, the CBC step of Algorithm~\ref{alg:construct_rand_reco_r1l_I_basic}
will successfully finish in line~\ref{alg:construct_rand_reco_r1l_I_basic:break1} with probability at least $1-\frac{\delta}{d-1}$. In this case, the arithmetic complexity of the corresponding CBC step is dictated by lines~\ref{alg:construct_rand_reco_r1l_I_basic:choose_T}, \ref{alg:construct_rand_reco_r1l_I_basic:shuffle}, and~\ref{alg:construct_rand_reco_r1l_I_basic:check_reco_property}.
For the first, we observe an arithmetic complexity that is bounded by $C_1 T\log T$, which is caused by randomly choosing a set $\mathcal{T}$ in line~\ref{alg:construct_rand_reco_r1l_I_basic:choose_T} when applying the strategy by Floyd, cf.~\cite{BF1987}. The  arithmetic complexity of line~\ref{alg:construct_rand_reco_r1l_I_basic:shuffle} is bounded by $C_2T$ when using the shuffling strategy by Fisher-Yates, cf.\ \cite[Algorithm~P]{Kn98a}.
When considering line~\ref{alg:construct_rand_reco_r1l_I_basic:check_reco_property}, we need to keep in mind that this line will be called up to $r_\delta$ times, which yields an arithmetic complexity in $\OO{\log(d/\delta)\,h_\ell(I,M)}$, where $h_\ell(I,M)$ is an upper bound
on the arithmetic complexity of a single call of the algorithm
${\texttt{check\_exactness}}(I_\ell, (z_1,\ldots,z_{\ell-1},y)^\top,M)$
that realizes the test for the exact integration property for fixed $\ell$ and arbitrary, fixed $y\in\{0,\ldots,M-1\}$.

Altogether, we obtain an arithmetic complexity in 
$
\OO{T\log T+h_\ell(I,M)\,\log\frac{d}{\delta}}
$
with probability $1-\frac{\delta}{d-1}$ for a single CBC step.
Using the complementary events and the union bound for all $d-1$ CBC steps yields that
the whole algorithm has an arithmetic complexity in 
$
\OO{d(T\log T+\log(d/\delta)\,h(I,M))}
$
with probability $1-\delta$, where $h(I,M)$ is a universal upper bound on $\max_{\ell\in\{2,\ldots,d\}}h_\ell(I,M)$.

In order to estimate the expected computational complexity, again we first consider a single CBC step.
Obviously, lines~\ref{alg:construct_rand_reco_r1l_I_basic:choose_T} and~\ref{alg:construct_rand_reco_r1l_I_basic:shuffle} will cause an arithmetic complexity in $\OO{T\log T}$ independent from the number of candidates that are checked during the CBC step.
If the $\ell$th CBC step will check more than $T$ different candidates, i.e., $r>T$, then line~\ref{alg:construct_rand_reco_r1l_I_basic:perm_M} can be realized using Fisher-Yates in an
arithmetic complexity in $\OO{M}$. The probability that this line is reached is bounded from above by $c^{-T}$.

It remains to estimate the expected number of calls of the function 
$${\texttt{check\_exactness}}(I_\ell, (z_1,\ldots,z_{\ell-1},z_\ell^{(r)})^\top,M)\,.$$
To this end, we estimate
\begin{align*}
\mathbb{P}({\texttt{check\_exactness}}&\textnormal{ is called $r$ times})\\
&=\mathbb{P}\left(\left\{z_\ell^{(1)},\ldots, z_\ell^{(r-1)} \textnormal{ are all failing candidates}\right\}\right)\\
&=
\mathbb{P}\left(\bigcap_{j=1}^{r-1}\left\{z_\ell^{(j)}\in\mathcal{Z}_\ell^j\right\}\right)\le c^{-r+1}
\end{align*}
Accordingly, we call the function {\texttt{check\_exactness}} $r$ times with a probability bounded from above by $c^{-r+1}$, which yields an expected number $R$ of calls
$$
R:=\sum_{r=1}^{M_{\operatorname{lb}}+1} r\, \mathbb{P}\left(\bigcap_{j=1}^{r-1}\left\{z_\ell^{(j)}\in\mathcal{Z}_\ell^j\right\}\right)<c\sum_{r=1}^\infty r\,c^{-r}=\left(\frac{c}{c-1}\right)^2\,.
$$
Together, we obtain an expected arithmetic complexity in 
$$\OO{T\log T+c^{-T}M+h(I,M)}$$
for a single CBC step. Due to the linearity of the expectation, we achieve
an additional factor $d-1$ for the expected complexity of Algorithm~\ref{alg:construct_rand_reco_r1l_I_basic}.

Analyzing the worst case complexity differs from analyzing the expected complexity only in the number of calls of 
$${\texttt{check\_exactness}}(I_\ell, (z_1,\ldots,z_{\ell-1},z_\ell^{(r)})^\top,M)\,,$$
which is $M_{\operatorname{lb}}+1\le M$ in the worst case. Thus, we achieve
$\OO{T\log T + M h(I,M)} = \OO{M h(I,M)}$ for a single CBC step when assuming $h(I,M)\gtrsim \log{M}$. Summing up the complexities of $d-1$ CBC steps yields the last assertion.
\end{proof}

\begin{algorithm}[t]
\small
\caption{Probabilistic CBC construction of exactly integrating rank-1 lattices}\label{alg:construct_rand_reco_r1l_I}
  \begin{tabular}{p{2cm}p{5cm}p{7.3cm}}
    Input: 	& $I\subset[-N,N]^d$ 	& frequency set\\
    		& $M$	& lattice size\\
    		& $T$	& maximal number of tested components
  \end{tabular}
		
  \begin{algorithmic}[1]
	\State $z_1=1$;
	\For{$\ell=2:d$}
		\State choose a set $\mathcal{T}$ of $T$ different elements from $\{0,\ldots,M-1\}$, uniformly at random\label{alg:construct_rand_reco_r1l_I:choose_Z_l}
		\State shuffle the elements of $\mathcal{T}$ and save $(z_\ell^{(r)})_{r=1}^T$, i.e., $\{z_\ell^{(r)}, r=1,\ldots,T\}=\mathcal{T}$\label{alg:construct_rand_reco_r1l_I:shuffle}
		\For{$r=1:T$}
			\If{{\texttt{check\_exactness}}($I_\ell$, $(z_1,\ldots,z_{\ell-1},z_\ell^{(r)})^\top$,$M$)$==$\texttt{true}}\label{alg:construct_rand_reco_r1l_I:check_reco_property}
				\State $z_\ell=z_\ell^{(r)}$
				\State \texttt{break}
			\ElsIf{$r==T$}
				\State raise \texttt{ERROR}
			\EndIf
		\EndFor
	\EndFor
  \end{algorithmic}
  \begin{tabular}{p{2cm}p{5cm}p{7.3cm}}
    Output: & $\boldz$ & generating vector of rank\mbox{-}1 lattice such that\\
	    & $\Lambda(\boldz,M)$ &  is a reconstructing rank\mbox{-}1 lattice for $I$
\\    \cmidrule{1-3}
	Complexity:&
    \multicolumn{2}{l}{$\OO{T\,d\,\left(\log{T}+h(I,M)\right)}$}
    \\\cmidrule{1-3}
    \multicolumn{3}{r}{\footnotesize$h(I,M)$ is an upper bound on the computational complexity of a single call of}\\ \multicolumn{3}{r}{\footnotesize\texttt{check\_exactness}($I_\ell$, $(z_1,\ldots,z_{\ell-1},y)^\top$,$M$) for fixed $\ell$ and arbitrary, fixed $y\in\{0,\ldots,M-1\}$ }
  \end{tabular}
\end{algorithm}

The last theorem motivates the modification of Algorithm~\ref{alg:construct_rand_reco_r1l_I_basic} in order to significantly reduce the worst case arithmetic complexity when accepting a prespecified failure probability.
The strategy is just to check only a reasonable number $T$ of component candidates in each CBC step, i.e., 
leaving out the computationally expensive parts of Lines~\ref{alg:construct_rand_reco_r1l_I_basic:perm_M} to~\ref{alg:construct_rand_reco_r1l_I_basic:endinnerfor} of Algorithm~\ref{alg:construct_rand_reco_r1l_I_basic}.
Certainly, this yields a non-deterministic and probabilistic approach.
For the sake of completeness we present the simpler strategy as Algorithm~\ref{alg:construct_rand_reco_r1l_I}.
At this point, we stress on the fact that Algorithm~\ref{alg:construct_rand_reco_r1l_I} is a probabilisitc algorithm, but nevertheless terminates in finite time.
Moreover, the algorithm returns admissible output in cases of success only, i.e., 
non-successful runs will be reported.
The insights gained in the proof of Theorem~\ref{thm:deterministic}
directly yield results on Algorithm~\ref{alg:construct_rand_reco_r1l_I}.

\begin{Corollary}\label{thm:alg1}
Let $1<c\in\R$, $I\subset\Z^d$, $\delta>0$, and $M$ being a prime number satisfying $M\ge cM_{\operatorname{lb}}$.
Choosing $T\in\N$, $T\ge\frac{\log(d-1)-\log{\delta}}{\log{c}}$, Algorithm~\ref{alg:construct_rand_reco_r1l_I} determines an exactly integrating rank-1 lattice for the frequency set $I$
with probability at least $1-\delta$.

The runtime of Algorithm~\ref{alg:construct_rand_reco_r1l_I} is in
$\OO{d\,T(\log T+ h(I,M))}$, where
$h(I,M)$ is an upper bound on the computational complexity of a single call of the function \texttt{check\_exactness}($I_\ell$, $(z_1,\ldots,z_{\ell-1},y)^\top$,$M$) for arbitrary $\ell\in\{2,\ldots,d\}$ and arbitrary $y\in\{0,\ldots,M-1\}$.
\end{Corollary}

\begin{proof}

The same argumentations as in the first part of the proof of Theorem~\ref{thm:deterministic} and taking $T\ge\frac{\log{(d-1)}-\log{\delta}}{\log{c}}$ into account yields $\delta$ as an upper bound for the failure probability. Thus the first assertion holds.

Let us comment on the computational complexity of Algorithm~\ref{alg:construct_rand_reco_r1l_I}. Due to the precomputation of the candidate sets $\{z_{\ell}^{(r)}, r=1,\ldots,T\}$ using the algorithm from~\cite{BF1987}, and some sorted structure for saving the set $\{z_{\ell}^{(r)}, r=1,\ldots,\tilde{t}\}$, $\tilde{t}<T$, we achieve a computational complexity for Line~\ref{alg:construct_rand_reco_r1l_I:choose_Z_l} in $\OO{T\log T}$. Furthermore, Line~\ref{alg:construct_rand_reco_r1l_I:shuffle} can be realized using Fisher-Yates~\cite[Algorithm~P]{Kn98a} in $\OO{T}$. Due to the definition of $h(I,M)$ we obtain computational costs caused by the maximum of $T$ passes of line~\ref{alg:construct_rand_reco_r1l_I:check_reco_property} that is bounded in $\OO{T\, h(I,M)}$ for each $\ell=2,\ldots,d$. Altogether, we achieve that
Algorithm~\ref{alg:construct_rand_reco_r1l_I} terminates and has a computational complexity in $\OO{d\,T (\log T+h(I,M))}$ in the worst case.
\end{proof}

\begin{Remark}\label{rem:discuss_T}
Let us comment on the expected computational complexity of Algorithm~\ref{alg:construct_rand_reco_r1l_I}, which obviously is less than the expected complexity of Algorithm~\ref{alg:construct_rand_reco_r1l_I_basic}, cf. Theorem~\ref{thm:deterministic}. In detail, the term that depend linearly on $M$ vanishes, i.e., we achieve a computational complexity which is in $\OO{d(T\log T+h(I,M)}$ for large enough $M>M_{\operatorname{lb}}$,
 in expectation. However, please note that Algorithm~\ref{alg:construct_rand_reco_r1l_I} may fail with certain probability for $T<M_{\operatorname{lb}}$ even in those cases where the success of the brute force approach in Algorithm~\ref{alg:construct_rand_reco_r1l_I_basic} is guaranteed.

Fixing $c>1$ and choosing $T:=\left\lceil\log_c{\frac{d-1}{\delta}}\right\rceil$ in Corollary~\ref{thm:alg1} causes a computational worst case complexity of Algorithm~\ref{alg:construct_rand_reco_r1l_I} which is in
$$
\OO{d\,\log(d/\delta)\,\left(\log\log(d/\delta)+h(I,M)\right)}
$$
and the Algorithm will succeed with probability at least $1-\delta$.
Since the probability estimates in the proof of Theorem~\ref{thm:deterministic} are very basic and pessimistic, we may achieve much better behavior in practice.

Moreover, due to \eqref{eq:prob_zero}, we observe that the failure probability is specified by $\delta=0$ for
$T>M_{\operatorname{lb}}$, e.g., $T=M$.
However, the computational complexity of Algorithm~\ref{alg:construct_rand_reco_r1l_I} is in $\Omega\left( d\,T \right)$ due to the shuffling in Line~\ref{alg:construct_rand_reco_r1l_I:shuffle}.
Fixing $T$ too large, this term may dominate the computational complexity of the whole algorithm even though a small number of randomly chosen $z_\ell^{(r)}$ is tested later on.
\end{Remark}

\begin{Remark}
We stress on the fact that the \emph{expected computational complexity} can be further decreased when using a slight modification of Algorithm~\ref{alg:construct_rand_reco_r1l_I}. This modification concerns Lines~\ref{alg:construct_rand_reco_r1l_I:choose_Z_l} and~\ref{alg:construct_rand_reco_r1l_I:shuffle}, which is the random preselection of the numbers $z_\ell^{(r)}$. Simply choosing the candidates $z_\ell^{(r)}\in\{0,\ldots,M-1\}$ uniformly distributed at random within the inner loop
is an alternative strategy, which yields that the term that depends on $T$ in the expected computational complexity vanishes.
In doing so, one could start the modified Algorithm~\ref{alg:construct_rand_reco_r1l_I} using a suitable parameter $T$ and analyze this approach in more detail similar to the proof of Corollary~\ref{thm:alg1}. This will lead to an improved expected computational complexity of this algorithm which does not suffer from logarithmic terms in the dimension $d$, i.e., the runtime will be strictly linear in the dimension $d$ in expectation.
Certainly, the modified algorithm would also suffer from a possible failure probability which can be estimated in a similar way as specified in the proof of Corollary~\ref{thm:alg1} for Algorithm~\ref{alg:construct_rand_reco_r1l_I}.
However, we decided to present Algorithm~\ref{alg:construct_rand_reco_r1l_I} as it is in order to
guarantee that we may check exactly $T$ different candidates and that the algorithm is still deterministic.
In particular, we are willing to accept the additional computational costs caused by Lines~\ref{alg:construct_rand_reco_r1l_I:choose_Z_l} and~\ref{alg:construct_rand_reco_r1l_I:shuffle}.
In Section~\ref{sec:check_exactness}, it turns out that these additional costs will be negligible since usually
the computational complexity $h(I,M)$ of \texttt{check\_exactness} is the strongly dominating part of the algorithm when choosing suitable parameters $T$.

On the other hand one can also improve the \emph{worst case computational complexity} by taking the recent observation of \cite{KuoMiNoNu19} into account. The strategy could be to estimate the computational effort that is needed to apply one CBC construction step from \cite{KuoMiNoNu19} and switching to that approach at that point, where Algorithm~\ref{alg:construct_rand_reco_r1l_I} has already spent a similar computational effort in a single CBC step. Thus, the worst case computational complexity would be in the order of two times of those that the CBC construction from \cite{KuoMiNoNu19} needs, i.e., $\OO{dM}$. The expected computational complexity is still bounded by the term that is given in Corollary~\ref{thm:alg1} when appropriately choosing the lattice size $M$ and the (possibly slightly modified) parameter $T$.
\end{Remark}

\begin{Remark}\label{rem:compcomp_alg12}
Both, Algorithms~\ref{alg:construct_rand_reco_r1l_I_basic} and~\ref{alg:construct_rand_reco_r1l_I}, are predestined
for parallelization since the function \texttt{check\_exactness} is independently called up to (at least) $T$ times in each CBC step.
However, since the expected number of \texttt{check\_exactness} calls is extremely small for large enough $M$, 
only a reasonable small number of parallel function calls is worthwhile.
A parallelization of \texttt{check\_exactness} itself seems to be more appropriate in order
to accelerate single thread performance.
\end{Remark}

\section{Discussion on \texttt{check\_exactness}}
\label{sec:check_exactness}

All algorithms presented in the last section apply the---up to now unexplained---function \texttt{check\_exactness} and use an
unspecified upper bound $h(I,M)$ on its computational costs in the theoretical investigations.
In this section, we will analyze this function in more detail. First,
we will consider the case where the full frequency set $I$ is given and we search for an
exactly integrating rank-1 lattice for this set of frequencies.
A second approach is the construction of reconstructing rank-1 lattices, which is in fact the
construction of an exactly integrating rank-1 lattice for a difference set $I=\mathcal{D}(I')$.
It turns out that for the latter approach the presented strategies from Section~\ref{sec:main}
will be most beneficial since for checking the exactness of the rank-1 lattice the knowledge of the
difference set is not necessary. The exactness can be checked dealing only with the original set $I'$.

In Section~\ref{sec:main} we particularly simplified the input of the function \texttt{check\_exactness}
in order to provide a general view on the main strategies presented in Algorithms~\ref{alg:construct_rand_reco_r1l_I_basic} and~\ref{alg:construct_rand_reco_r1l_I}.
Going into more detail, the exactness check for the $\ell$th component can even be done without 
knowing the full frequency set.
Certainly, the lattice size $M$ and the designated component candidate $y$ need to be fixed.
Then, it is enough to know some temporary result from the successful exactness check of the $(\ell-1)$th CBC step,
which directly led to the $\ell$th CBC step,
together with the corresponding $\ell$th components of the frequencies within $I$
in order to decide whether or not $y$ is an admissible CBC choice, cf.\ also e.g.\ \cite[Remark~24]{KuoMiNoNu19}.
Since the determination of exactly integrating rank-1 lattices and reconstructing rank-1 lattices differ in the
used approaches which put the function \texttt{check\_exactness} into practice,
and these differences are crucial for observing the computational complexities given in Table~\ref{tab:intro_complexities}, we believe that it is necessary or at least very helpful to discuss in detail the two different approaches we utilize.

\begin{algorithm}[t]
\small
\caption{\texttt{check\_exactness} of the rank-1 lattice $\Lambda((z_1,\ldots,z_{\ell-1},y)^\top,M)$ for the frequency set $I_\ell$ under the assumption that $\Lambda((z_1,\ldots,z_{\ell-1})^\top,M)$ is an exactly integrating rank-1 lattice for the frequency set $I_{\ell-1}$.}\label{alg:check_exactness1}
  \begin{tabular}{p{1.75cm}p{5.25cm}p{7.3cm}}
      Input: 	& $(k_{j,\ell})_{j=1}^{|I|}$ 	& \small $\ell$th components of the frequency vectors in the set $I$\\
    		& $\boldsymbol{\nu}=\nu(I,(z_1,\ldots,z_{\ell-1},\boldzero)^\top,M)$ &  \small residues modulo $M$ of inner products of the elements in $I$ with generating vector $(z_1,\ldots,z_{\ell-1},\boldzero)^\top$\\
    		& $M$	& \small lattice size\\
    		& $y \in [0,M)\cap\Z$ & \small component candidate of generating vector
  \end{tabular}
  \begin{algorithmic}[1]
  	\State $\operatorname{exact}=\texttt{false}$
  	\State determine the sets $J:=\{j\in\{1,\ldots,|I|\}\colon k_{j,\ell}\neq 0\}$
  	  	\label{alg:check_exactness1_construct_J}
  	\State compute $\boldsymbol{\nu}:=\left(\boldsymbol{\nu}+\left(y\, k_{j,\ell}\right)_{j=1}^{|I|}\right)\bmod M$\label{alg:check_exactness1_compute_nu}
\If{$0\not\in \{\nu_j\colon j\in J\}$}\label{alg:check_exactness1_findzeros}
	  	\State $\operatorname{exact}=\texttt{true}$
  	\EndIf
  \end{algorithmic}
  \begin{tabular}{p{1.75cm}p{5.25cm}p{7.3cm}}
    Output: & $\operatorname{exact}$ & \small boolean, \texttt{true} in case that $\Lambda((z_1,\ldots,z_{\ell-1},y)^\top,M)$ exactly integrates each $p\in\Pi_{I_\ell}$\\
	    & $\nu(I,(z_1,\ldots, z_{\ell-1},y,\boldzero)^\top,M)=\boldsymbol{\nu}$ &  \small current scalar products in case that $\operatorname{exact}=\texttt{true}$
\\    \cmidrule{1-3}
  \end{tabular}
  \begin{tabular}{p{1.75cm}p{5.25cm}p{7.3cm}}
	Complexity:&
$\OO{|I|}$&
  \end{tabular}
\end{algorithm}

We start with the algorithm that can be used in CBC constructions in order to determine the exact integration property of rank-1 lattices for given frequency sets $I$.
For practical implementations, we assume that the elements of the frequency set $I$ have some unique numbering, i.e., $I=\{\boldk_j\in\Z^d\colon j=1,\ldots,|I|\}$, which can be a ``natural order'' of the elements of $I$ given by
\begin{equation}
\begin{split}
\boldk_j<\boldk_{j+1} \Leftrightarrow\; \exists &\,t_0\in\{1,\ldots d\} \text{ s.t. }\\
& k_{j,t}=k_{j+1,t}\, \forall t\in\{1,\ldots,t_0-1\} \text{ and } k_{j,t_0}<k_{j+1,t_0},\label{eq:natural_sort}
\end{split}
\end{equation}
$j=1,\ldots,|I|-1$.
For a given rank-1 lattice $\Lambda(\boldz,M)$, $\boldz\in\Z^d$, $M\in \N$, we define the
vector of inner products $\nu(I,\boldz,M)$ of all elements of $I$ with the generating vector $\boldz$ modulo $M$ as
$$
\nu(I,\boldz,M):=\left(\boldk_j\cdot\boldz\bmod M\right)_{j=1}^{|I|}.
$$
A simple exercise using the compatibility of the modulo operation with addition proves
$$\nu(I,(z_1,\ldots,z_\ell,\boldzero)^\top,M)=\left(
\nu(I,(z_1,\ldots,z_{\ell-1},\boldzero)^\top,M)+\left(k_{j,\ell}z_\ell\right)_{j=1}^{|I|}\right)\bmod M\,,$$
from which follows that computing the inner products required for checking a single component candidate $y$ in a CBC construction has a computational complexity in $\OO{|I|}$ for each CBC step, cf. line~\ref{alg:check_exactness1_compute_nu} of Algorithm~\ref{alg:check_exactness1}.
In each $\ell$th CBC step, $\ell=2,\ldots,d$, we assume that the rank-1 lattice $\Lambda((z_1,\ldots,z_{\ell-1})^\top,M)$ has the exact integration property for the frequency set 
$$I_{\ell-1}:=\{(k_1,\ldots,k_{\ell-1})^\top\in\Z^{\ell-1}\colon\boldk\in I\},$$
from which follows that the rank-1 lattice $\Lambda((z_1,\ldots,z_{\ell-1},y)^\top,M)$, $y\in\Z$, has the 
exact integration property for the set $I_{\ell}$, iff $\left(\nu(I,(z_1,\ldots,z_{\ell-1},y,\boldzero)^\top,M)\right)_j\neq 0$
for all those $j$ where $k_{j,l}\neq 0$ holds.
We denote the set of those indices by $J$, cf. line~\ref{alg:check_exactness1_construct_J} of Algorithm~\ref{alg:check_exactness1}.
Clearly, $J$ can be determined in a runtime which is linear in $|I|$, i.e., $\OO{|I|}$.
The remaining effort of Algorithm~\ref{alg:check_exactness1} is caused by line~\ref{alg:check_exactness1_findzeros}
and is just searching the first zero in a vector of length $|I|$, which has a complexity of $\OO{|I|}$.
Accordingly, the total computational complexity of Algorithm~\ref{alg:check_exactness1} is linear in $|I|$.

It remains to show, that Algorithm~\ref{alg:check_exactness1} does what it is supposed to.
At this point we stress on the fact that the exact integration property of the rank-1 lattice $\Lambda((z_1,\ldots,z_{\ell-1})^\top,M)$ for the projected frequency set $I_{\ell-1}$,
is an essential assumption in order to obtain the correctness of Algorithm~\ref{alg:check_exactness1}.

\begin{Lemma}
For given frequency set $I\subset\Z^d$, $|I|<\infty$, lattice size $M\in\N$ and $\ell\in\N$, $\ell\ge 2$, we assume $\Lambda((z_1,\ldots,z_{\ell-1})^\top,M)\subset\T^{\ell-1}$ is an exactly integrating rank-1 lattice for $I_{\ell-1}$, i.e.,
$$
\boldh\cdot(z_1,\ldots,z_{\ell-1})^\top\not\equiv 0\mod{M}\quad\text{for all }\boldh\in I_{\ell-1}\setminus\{\boldzero\}\,.
$$
Then, for given $y\in\Z$ Algorithm~\ref{alg:check_exactness1} determines whether or not
the $\ell$-dimensional rank-1 lattice
$\Lambda((z_1,\ldots,z_{\ell-1},y)^\top,M)\subset\T^\ell$
is an exactly integrating one for $I_\ell$.
\end{Lemma}

\begin{proof}
We work directly on the set $I$ instead of $I_\ell$ in order to avoid the computational effort caused by the detection of duplicates.
As a consequence, we may check a single element of $I_\ell$ several times, but that does not matter. 
Moreover, the exact integration property of $\Lambda((z_1,\ldots,z_{\ell-1},y,0,\ldots,0)^\top,M)\subset\T^\ell\times\{\boldzero\}\subset\T^d$ for all frequencies $\boldk\in I\setminus\{\boldh\in I\colon h_1=\ldots,=h_\ell=0\}$ implies the exact integration property of
$\Lambda((z_1,\ldots,z_{\ell-1},y)^\top,M)\subset\T^\ell$ for $I_\ell:=\{(h_{1},\ldots,h_{\ell})^\top\in\Z^\ell\colon \boldh\in I\}$.

For $\boldk_j\in I$, we distinguish three different cases:
\begin{enumerate}
\item \underline{$(k_{j,1},\ldots,k_{j,\ell})^\top=\boldzero$:}
We observe $\sum_{s=1}^{\ell-1}k_{j,s}z_s+yk_{j,\ell}=0$ as required for all $y\in\Z$ without doing anything in Algorithm~\ref{alg:check_exactness1}.
\item \underline{$(k_{j,1},\ldots,k_{j,\ell-1})^\top\neq\boldzero$ and $k_{j,\ell}=0$:}\newline
We get $\sum_{s=1}^{\ell-1}k_{j,s}z_s+yk_{j,\ell}=\sum_{s=1}^{\ell-1}k_{j,s}z_s\neq 0$ as required for all $y\in\Z$, 
due to the fact that $\Lambda((z_1,\ldots,z_{\ell-1})^\top,M)\subset\T^{\ell-1}$ is an exactly integrating rank-1 lattice for $I_{\ell-1}$. Similar as in the last case, we do not need to check anything in Algorithm~\ref{alg:check_exactness1}.
\item \underline{$k_{j,\ell}\neq0$:}
In this case we have $j\in J:=\{i\in\{1,\ldots,|I|\}\colon k_{i,\ell}\neq 0\}$. Algorithm~\ref{alg:check_exactness1}
checks whether or not at least one of the inner products $\sum_{s=1}^{\ell-1}k_{j,s}z_s+yk_{j,\ell}\neq 0$, $j\in J$,
is zero modulo $M$.
If at least one of them is zero, the rank-1 lattice $\Lambda((z_1,\ldots,z_{\ell-1},y)^\top,M)$ can not exactly integrate at least one of the monomials $\e^{2\pi\ii(k_{j,1},\ldots,k_{j,\ell})^\top\cdot\circ}$, $j\in J$, and Algorithm~\ref{alg:check_exactness1} returns \texttt{false}.
Otherwise the rank-1 lattice $\Lambda((z_1,\ldots,z_{\ell-1},y)^\top,M)$ exactly integrates all monomials $\e^{2\pi\ii(k_{j,1},\ldots,k_{j,\ell})^\top\cdot\circ}$, $j\in J$, and Algorithm~\ref{alg:check_exactness1} returns \texttt{true}.
\end{enumerate}
Accordingly, Algorithm~\ref{alg:check_exactness1} realizes the exactness check of a CBC step under the assumption that
$\Lambda((z_1,\ldots,z_{\ell-1})^\top,M)$ is already an exactly integrating rank-1 lattice for the frequency set $I_{\ell-1}$.
\end{proof}

As mentioned above, there is a more clever way of checking the exact integration property for a difference set $D(I')$, which directly uses the reconstruction property of $I'$, cf.~\eqref{eq:reco_property}.
Algorithm~\ref{alg:check_exactness} states the strategy in detail.
The crucial advantage is that the reconstruction property can be checked using only the inner products of the generating vector with all frequencies in $I'$ instead of all frequencies in $D(I')$, which leads to substantially reduced computational costs -- at least in cases where $|I'|$ and $|D(I')|$ differ significantly.

\begin{algorithm}[t]
\small
\caption{Check exactness of the rank-1 lattice $\Lambda((z_1,\ldots,z_\ell)^\top,M)$ for the difference set of $D(I_\ell)$ under the assumption that $\Lambda((z_1,\ldots,z_{\ell-1})^\top,M)$ is an exactly integrating rank-1 lattice for the difference set $D(I_{\ell-1})$}\label{alg:check_exactness}
  \begin{tabular}{p{1.75cm}p{5.25cm}p{7.3cm}}
    Input: 	& $(k_{j,\ell})_{j=1}^{|I|}$ 	& \small $\ell$th components of the frequency vectors in the set $I$\\
    		& $\nu=\nu(I,(z_1,\ldots,z_{\ell-1},\boldzero)^\top,M)$ &
    		\small residues modulo $M$ of inner products of the elements in $I$ with generating vector $(z_1,\ldots,z_{\ell-1},\boldzero)^\top$\\
    		& $M$	& \small lattice size\\
    		& $y \in [0,M)\cap\Z$ & \small component candidate of generating vector
  \end{tabular}
  \begin{algorithmic}[1]
  	\State $\operatorname{exact}=\texttt{false}$
  	\State construct the set $\tilde{I}:=\{(\nu_j,k_{j,\ell})\colon j=1,\ldots,|I|\}\subset\Z^2$\label{alg:check_exactness_construct_tilde_I}
  	\If{$|\{\boldh\cdot(1,y)^\top\mod M\colon \boldh\in\tilde{I}\}|==|\tilde{I}|$}\label{alg:check_exactness_scalar_products_tilde_I}
	  	\State $\operatorname{exact}=\texttt{true}$
  		\State $\nu(I,(z_1,\ldots,z_{\ell-1},y,\boldzero)^\top,M):=\left(\nu+\left(y\,k_{j,\ell}\right)_{j=1}^{|I|}\right)\bmod M$\label{alg:check_exactness_scalar_products}
  	\EndIf
  \end{algorithmic}
  \begin{tabular}{p{1.75cm}p{5.25cm}p{7.3cm}}
    Output: & $\operatorname{exact}$ & \small boolean, \texttt{true} in case that $\Lambda((z_1,\ldots,z_{\ell-1},y)^\top,M)$ exactly reconstructs each $p\in\Pi_{I_\ell}$\\
	    & $\nu(I,(z_1,\ldots, z_{\ell-1},y,\boldzero)^\top,M)=\boldsymbol{\nu}$ &  \small current scalar products in case that $\operatorname{exact}=\texttt{true}$
\\    \cmidrule{1-3}
  \end{tabular}
  \begin{tabular}{p{1.75cm}p{5.25cm}p{7.3cm}}
	Complexity:& $\OO{|I|\log|I|}$ &
  \end{tabular}
\end{algorithm}

\begin{Lemma}
For given frequency set $I\subset\Z^d$, $|I|<\infty$, lattice size $M\in\N$ and $\ell\in\N$, $\ell\ge 2$,
we assume $\Lambda((z_1,\ldots,z_{\ell-1})^\top,M)\subset\T^{\ell-1}$ is a reconstructing rank-1 lattice for $I_{\ell-1}$, i.e.,
$$
\boldh\cdot(z_1,\ldots,z_{\ell-1})^\top\not\equiv \boldk\cdot(z_1,\ldots,z_{\ell-1})^\top\mod{M}\quad\text{for all }\boldh,\boldk\in I_{\ell-1}, \boldh\neq\boldk\,.
$$
Then, for given $y\in\Z$ Algorithm~\ref{alg:check_exactness} determines whether or not
$\Lambda((z_1,\ldots,z_{\ell-1},y)^\top,M)\subset\T^\ell$
is a reconstructing rank-1 lattice for $I_\ell$.
\end{Lemma}

\begin{proof}
\sloppy
Let $\ell>1$.
The frequency set $I_\ell\subset\Z^\ell$ can be determined in the following way
$$I_\ell:=\{(k_{j,1},\ldots,k_{j,\ell})^\top\in\Z^\ell\colon\boldk\in I\}\,.$$
Moreover, we consider the set $\tilde{I}:=\{(\nu_j,k_{j,\ell})\colon j=1,\ldots,|I|\}$
and observe that the mapping 
$$
(k_{j,1},\ldots,k_{j,\ell})^\top
\mapsto(\nu_j,k_{j,l})
$$
is a bijection between the sets $I_\ell$ and $\tilde{I}$ due to the fact that 
$$
(k_{j,1},\ldots,k_{j,\ell-1})^\top
\mapsto \nu_j
$$
is a bijective mapping from $I_{\ell-1}$ to the set 
$$\{\nu_j:=(k_{j,1},\ldots,k_{j,\ell-1})^\top\cdot(z_1,\ldots,z_{\ell-1})^\top\bmod M\colon(k_{j,1},\ldots,k_{j,\ell-1})^\top\in I_{\ell-1} \},$$
which holds because of the reconstruction property of  $\Lambda((z_1,\ldots,z_{\ell-1})^\top,M)$ for $I_{\ell-1}$.
Due to the definition of $\nu_j$ and the modulo compatibility of the addition, we observe
$$
(k_{j,1},\ldots,k_{j,\ell})^\top\cdot(z_1,\ldots,z_{\ell-1},y)^\top\equiv \nu_j+yk_{j,\ell}\mod{M},
$$
i.e., the set $\{\boldh\cdot(1,y)^\top\mod M\colon \boldh\in\tilde{I}\}$ contains exactly the values $(k_{1},\ldots,k_{\ell})^\top\cdot(z_1,\ldots,z_{\ell-1},y)^\top$, $(k_{1},\ldots,k_{\ell})^\top\in I_\ell$, modulo $M$.
In the case that these values are all unique, we observe
$$|\{\boldh\cdot(1,y)^\top\mod M\colon \boldh\in\tilde{I}\}|=|\tilde{I}|=|I_\ell|$$
distinct values of the inner products modulo $M$ and thus the reconstruction property of $\Lambda((z_1,\ldots,z_{\ell-1},y)^\top,M)$ for $I_\ell$ has proved to be true. On the other hand, in cases where we observe
$$|\{\boldh\cdot(1,y)^\top\mod M\colon \boldh\in\tilde{I}\}|<|\tilde{I}|=|I_\ell|,$$
at least two distinct elements in $I_\ell$ have the same inner product with $(z_1,\ldots,z_{\ell-1},y)$ modulo $M$ and thus the reconstruction property, cf.~\eqref{eq:reco_property}, for $I_\ell$ does not hold true.
\end{proof}

As mentioned above, Algorithm~\ref{alg:check_exactness} substantially reduces the computational effort of the
exactness check of a difference set $D(I')$ to $\OO{|I'|\log|I'|}$ instead of $\OO{|D(I')|}$ when using Algorithm~\ref{alg:check_exactness1}.
The cost-intensive part of Algorithm~\ref{alg:check_exactness} is the elimination of the duplicates in order to determine 
the set $\tilde{I}$ in line~\ref{alg:check_exactness_construct_tilde_I} for which we need to sort the set $\tilde{I}$. A similar worst-case runtime is caused by line~\ref {alg:check_exactness_scalar_products_tilde_I} in order to determine the number of elements within the set $\{\boldh\cdot(1,y)^\top\mod M\colon \boldh\in\tilde{I}\}$ of current inner products.
Both sets are constructed in $\OO{|I'|\log|I'|}$ arithmetic operations.
The remaining costs of Algorithm~\ref{alg:check_exactness} are dominated by line~\ref{alg:check_exactness_scalar_products}, which
is in $\OO{|I'|}$.

At this point,  we would like to mention that for frequency sets $I'$ ordered as mentioned in \eqref{eq:natural_sort} the sets
$\tilde{I}$ can be constructed in linear time $\OO{|I'|}$. However, precomputing this order of the frequency set $I'$ causes
an arithmetic complexity in $\OO{|I'|\log|I'|+d|I'|}$ at least once, but it may save computing time in practice.

\begin{Remark}\label{rem:diff_insetad_of_direct}
At first glance, Algorithms~\ref{alg:check_exactness1} and~\ref{alg:check_exactness} can be used for determining exactly integrating and reconstructing rank-1 lattices, respectively.
Due to the fact that the reconstruction property for the frequency set $I'$ coincides with the exact integration property for the corresponding difference set $D(I')$,
one can combine both algorithms in order to construct exactly integrating rank-1 lattices. In particular for frequency sets $I$ that fulfill the subset relation
$I\subset D(I^{(1)})\cup I^{(2)}$ for some $I^{(1)}, I^{(2)}\subset\Z^d$ with $\sqrt{|I|}\le|I^{(1)}|+|I^{(2)}|\ll |I|$, it can be advantageous to check the exact integration property
of the rank-1 lattice $\Lambda((z_1,,\ldots,z_{\ell-1},y),M)$ for $I_\ell$ by just checking the reconstruction property for $I^{(1)}_\ell$ using Algorithm~\ref{alg:check_exactness}
and checking the exact integration property for $I^{(2)}_\ell$ using Algorithm~\ref{alg:check_exactness1}. However, a suitable partitioning
needs to be known in order to apply this approach and -- roughly speaking -- $|D(I^{(1)})\cup I^{(2)}|\lesssim|I|$ should hold in order to observe comparable values of
$M_{\operatorname{lb}}$ for both frequency sets $D(I^{(1)})\cup I^{(2)}$ and $I$.
\end{Remark}

\begin{Remark}
We close this section discussing some implementation issues.
When choosing $M$ reasonably small, we may utilize Algorithm~\ref{alg:check_exactness1} multiple times for fixed frequency set $I$, fixed lattice size $M$, fixed $\ell$, and several $y$. Then, it is enough to perform line~\ref{alg:check_exactness1_construct_J} of Algorithm~\ref{alg:check_exactness1} only once, since the set $J$ is independent of $y$. A similar redundancy can be avoided when applying line~\ref{alg:check_exactness_construct_tilde_I} of Algorithm~\ref{alg:check_exactness} only once for fixed $I$, $M$, $\ell$, and varying $y$.

Furthermore, it is possible to parallelize Algorithms~\ref{alg:check_exactness1} and~\ref{alg:check_exactness}
as mentioned in Remark~\ref{rem:compcomp_alg12}.
In Algorithm~\ref{alg:check_exactness1}, one could simply chop the input vectors $(k_{j,\ell})_{j=1}^{|I|}$
and $\boldnu$ in a suitable way and work on the separate parts in parallel. Then, the output of Algorithm~\ref{alg:check_exactness1} is just the logical conjunction of all partial results.

In Algorithm~\ref{alg:check_exactness}, we need to eliminate duplicates (line~\ref{alg:check_exactness_construct_tilde_I}) or search for a duplicate (line~\ref{alg:check_exactness_scalar_products_tilde_I})
within integer vectors which does not allow for simply chopping the input in several parts. However, the efficient elimination of duplicates is based on sorting integer sequences, which can be done using parallel sorting algorithms, e.g., a parallel radix sort, cf.~\cite{ObKaFaSh19,AxWiFeSa20} and references therein.
Then, the duplicates can be found similar to the search for zeros in line~\ref{alg:check_exactness1_findzeros} of Algorithm~\ref{alg:check_exactness1}.
Thus, this step can also be parallelized.
\end{Remark}

\section{Probabilistic approach for determining suitable lattice sizes}\label{sec:prob_M}

The theoretical results in Section~\ref{sec:main} mainly depend on the choice of suitable lattice sizes $M$ for which the knowledge of $M_{\operatorname{lb}}$ is required. However, there exist simple upper bounds on $M_{\operatorname{lb}}$, which can be used in order to choose $M$ such that $M\ge cM_{\operatorname{lb}}$ is guaranteed. Of course, two crucial issues may arise using these upper bounds. On the one hand, simple to determine bounds may suffer from low accuracy, i.e., one might choose the lattice sizes $M$ much to large. On the other hand, computing these bounds should not cause excessive computational costs. In particular in the case where one is interested in reconstructing rank-1 lattices for $I'$, known well-suited bounds depend mainly on the difference set $D(I')$ and not on the frequency set $I'$ itself, cf.~\cite{Kae2013, KuoMiNoNu19}. Accordingly, one has to determine the difference set $D(I')$ from $I'$ which cause computational costs in $\mathcal{O}(d|I'|^2)$ in general. In particular for unstructured frequency sets $I'$, this leads to significantly higher computational costs than those caused by the whole CBC construction for reconstructing rank-1 lattices in Algorithm~\ref{alg:construct_rand_reco_r1l_I} using Algorithm~\ref{alg:check_exactness} for checking the reconstruction property for $I'$ (or, equivalently, the exact integration property for $D(I')$).

Moreover, even the smallest known upper bounds on $M_{\operatorname{lb}}$ are based on considerations that assume worst case scenarios, which might not be expected in general.
For that reason, we suggest a simple algorithm that just starts with a lattice size $M$ based on a very rough estimate on 
$M_{\operatorname{lb}}$ and tries to determine a suitable generating vector using Algorithm~\ref{alg:construct_rand_reco_r1l_I}.
Then we continue iteratively. As long as Algorithm~\ref{alg:construct_rand_reco_r1l_I} succeeds, we continue with a new lattice size $M$ that is roughly half of the lattice size $M$ of the last successful lattice search.
The systematic reduction of the lattice size $M$
is an obvious approach that can decisively reduce the sampling values one spends
for the exact integration or the reconstruction problem,
in particular in cases where suitable upper bounds on $M_{\operatorname{lb}}$ are hard to determine.
A similar approach is also suggested
in \cite[Remark~24]{KuoMiNoNu19} without going into details.
The strategy we propose is indicated in Algorithm~\ref{alg:heuristic}.
Since our approach will be affected from a certain small failure probability,
we analyze Algorithm~\ref{alg:heuristic} in more detail and discuss the computational complexities 
later on. Moreover, we compare the approach to the non-probabilistic state of the art algorithm from \cite{KuoMiNoNu19}.

First, we prepare the discussion on the computational complexity of Algorithm~\ref{alg:heuristic}
by considering the $\operatorname{\texttt{nextprime}}$ function defined as
$$
\operatorname{\texttt{nextprime}}(x):=\min\{p \text{ prime}\colon p>x\}.
$$

\begin{Lemma}\label{lem:primes}
Let $p_1>2$ be a prime number. Then the smallest prime $p_2$ larger than $p_1/2$ fulfills
$$p_2\le\frac{5p_1}{7}.$$
\end{Lemma}
\begin{proof}
For the first six odd prime numbers $p_1$, one easily proves the assertion using Table~\ref{tab:primes}
and observes the maximal ratio for $p_1=7$ and $p_2=5$.
So let $p_1> 17$, i.e., $p_1\ge 19$. 

Due to \cite{Na52}, for each $x\ge 9$, there exists a prime number in the interval
$(x,4x/3)$. For $p_1 \ge 19$, we have $x=(p_1-1)/2 \ge 9$ and the last mentioned fact
yields that there exists a prime $p_2$ in the interval
$\big((p_1-1)/2, 4(p_1-1)/6\big)$ and since $p_1$ is odd, we have $p_2>p_1/2$.
We estimate
$$
p_2:=\operatorname{\texttt{nextprime}}(p_1/2)=\min\{p \text{ prime}\colon p>p_1/2\}
<4(p_1-1)/6< 2p_1/3 < 5p_1/7\,.
$$
\end{proof}
\begin{table}[tb]
\centering
\begin{tabular}{lrrrrrr}
\toprule
$p_1$ & 3 & 5 & 7 & 11 & 13 & 17\\
$p_2$ & 2 & 3 & 5 &  7 &  7 & 11 \\
\bottomrule
\end{tabular}
\caption{Table of first primes illustrating the statement of Lemma~\ref{lem:primes}}
\label{tab:primes}
\end{table}

Again, we stress on the fact that $M\ge cM_{\operatorname{lb}}$ is one essential assumption in Theorem~\ref{thm:deterministic} and Corollary~\ref{thm:alg1}. A simple upper bound on $M_{\operatorname{lb}}$ is given by
\begin{equation}
M_{\operatorname{lb}}\le \max\big(|I\setminus\{\boldzero\}|+1,\max(I)\big)\le \max\big(|I|+1,\max(I)\big)=:\tilde{M}_{\text{int}}(I)=:\tilde{M}_{\text{int}} \label{def:tildeM_int}
\end{equation}
when constructing exactly integrating rank-1 lattices, cf.\ \cite[Lemma 4]{KuoMiNoNu19}, or 
\begin{align}
M_{\operatorname{lb}}&\le  \max((|D(I')|+1)/2,N_{I'})\nonumber\\
&\le \underbrace{\max(|I'|^2/2,N_{I'})}_{=:\tilde{M}_{\text{reco}}(I')=:\tilde{M}_{\text{reco}}} \le\max\left(|I'|^2/2,2\,\max(I')\right),\label{def:tildeM_reco}
\end{align}
when constructing reconstructing rank-1 lattices for $I'$ with $D(I')>1$, 
cf.\ \cite[Corollary~1]{Kae2013} together with \cite[Lemma~2.2]{Kae17}, and \cite[Theorem~2.1]{PoVo14} for the last estimate.

In the following, we focus on the construction of exactly integrating rank-1 lattices for the frequency set $I$ using  Algorithm~\ref{alg:heuristic}.
Slight modifications for constructing reconstructing rank-1 lattices are later discussed in Remark~\ref{rem:heuristic_complex_reco}.

As a consequence of Lemma~\ref{lem:primes} and \eqref{def:tildeM_int}, we observe that Algorithm~\ref{alg:heuristic} tests a number of lattice sizes that is at most logarithmic in 
$\tilde{M}_{\text{int}}$.
Moreover, 
the lattice sizes $\tilde{M}$ decrease and, especially for smaller lattice sizes~$\tilde{M}$, we expect significant failure probabilities when applying Algorithm~\ref{alg:construct_rand_reco_r1l_I}, which may be reduced by starting several, i.e., at most $K$, independent instances of Algorithm~\ref{alg:check_exactness_construct_tilde_I}.
Accordingly, the complexity of Algorithm~\ref{alg:heuristic} is bounded by terms in $\OO{K\log \tilde{M}_{\text{int}}}$ times the complexity of Algorithm~\ref{alg:check_exactness_construct_tilde_I} plus the computational expense of at most $\OO{\log\tilde{M}_{\text{int}}}$ many calls of the $\operatorname{\texttt{nextprime}}$ function. The latter can be realized by sequentially testing growing integers starting from the input of $\operatorname{\texttt{nextprime}}$ for primality. Using an advanced AKS primality test, the computational complexity of each test is  in $\OO{\left(\log\tilde{M}_{\text{int}}\right)^c}$, $c<7$, cf.~\cite{AKS04, LePo19}.
Due to~\cite{BaHaPi01}, the number of primality tests necessary to find the next prime is in $\OO{\tilde{M}_{\text{int}}^{0.525}}$.
Accordingly, the total number of arithmetic operations Algorithm~\ref{alg:heuristic} spends for all $\operatorname{\texttt{nextprime}}$ calls is in $\OO{\tilde{M}_{\text{int}}^{0.525+\varepsilon}}$ for arbitrarily small $\varepsilon>0$. 

\begin{algorithm}[t]
\small
\caption{Heuristic CBC construction of exactly integrating rank-1 lattices}\label{alg:construct_rand_reco_r1l_I_2}
\label{alg:heuristic}
  \begin{tabular}{p{1.75cm}p{5.0cm}p{7.55cm}}
    Input: 	& $I\subset[-\max(I),\max(I)]^d$ 	& frequency set\\
    		& $K$					& number of iterations per lattice size candidate $\tilde{M}$\\
    		& $T$					& maximal number of tests in each CBC step for the components of the generating vector $\boldz$
  \end{tabular}

  \begin{algorithmic}[1]
  		\State $\tilde{M}:=\operatorname{\texttt{nextprime}}(2\,\max(|I|+1, \max(I)))$\label{alg:heuristic_nextprime1}
  		\State $k=0$
  		\While{$k<K$}
  			\State $k=k+1$
	  		\While{Algorithm~\ref{alg:construct_rand_reco_r1l_I} with input $I$, $\tilde{M}$, $T$ succeeds with output $\tilde{\boldz}$}
					\State $M=\tilde{M}$, $\boldz=\tilde{\boldz}$
					\If{$\tilde{M}\neq2$}
		  		  		\State $\tilde{M}:=\operatorname{\texttt{nextprime}}(\tilde{M}/2)$\label{alg:heuristic_nextprime2}
		  		  		\State $k=0$
			  		\EndIf
			\EndWhile
  		\EndWhile
  \end{algorithmic}
  \begin{tabular}{p{1.75cm}p{5.0cm}p{7.55cm}}
    Output: & $\boldz$ & generating vector and\\
			& $M$ & lattice size of rank-1 lattice such that\\
	    & $\Lambda(\boldz,M)$ &  is an exactly integrating rank-1 lattice for $I$
\\    \cmidrule{1-3}
    	Complexity:&
    $\OO{K\,T\,d\,\left(\log{T}+h(I,M)\right)\log m_I + m_I^{0.525+\varepsilon}}$;&
    \multicolumn{1}{r}{$m_I:=\max(|I|,\max(I))$}
        \\\cmidrule{1-3}
        \multicolumn{3}{r}{\footnotesize$h(I,M)$ is an upper bound on the computational complexity of a single call of}\\ \multicolumn{3}{r}{\footnotesize\texttt{check\_exactness}($I_\ell$, $(z_1,\ldots,z_{\ell-1},y)^\top$,$M$) for fixed $\ell$ and arbitrary, fixed $y\in\{0,\ldots,M-1\}$ }
  \end{tabular}
\end{algorithm}

According to the last considerations, Algorithm~\ref{alg:heuristic} is a probabilistic algorithm that determines
reasonable lattice sizes as well as corresponding generating vectors of exactly integrating rank-1 lattices in cases of success.
Please note that the algorithm, although
probabilistic, is guaranteed to terminate after a finite number of arithmetic operations, and
that upper bounds on this computational effort can be estimated in advance, cf.\ Theorem~\ref{thm:alg_heuristic}.
However, a specific probability remains that Algorithm~\ref{alg:heuristic} does not construct any rank-1 lattice.
One can control this failure probability when choosing suitable parameters $K$ and $T$.
In fact, suitable parameters even allow to estimate the resulting lattice sizes in its order of magnitude and corresponding probabilities.
It turns out, that Algorithm~\ref{alg:heuristic} determines rank-1 lattices of sizes less than $4M_{\operatorname{lb}}$ with high probability when choosing suitable small parameters $K$ and $T$ without precomputing the number (or rather tight upper bounds on) $M_{\operatorname{lb}}$.

\begin{sloppypar}
At this point, it proves beneficial to simplify the notation by defining
\begin{equation}
h(I):=\max_{M \text{ prime}, M\le \tilde{M}_{\text{int}}}h(I,M),\label{def:hI}
\end{equation}
which is a universal upper bound on the computational complexity of a single call of \texttt{check\_exactness}($I_\ell$, $(z_1,\ldots,z_{\ell-1},y)^\top$,$M$) for each combination of $\ell\in\{1,\ldots,d\}$, $y\in\{0,\ldots,M-1\}$, and primes $M$ in the range $[2,\tilde{M}_{\text{int}}]$, i.e., the prime numbers that may be used by Algorithm~\ref{alg:heuristic} as lattice size candidates $\tilde{M}$, cf.\ Algorithm~\ref{alg:construct_rand_reco_r1l_I} for the definition of $h(I,M)$.
According to the discussions in Section~\ref{sec:check_exactness}, we can bound $h(I)$ by terms in $\OO{h(I,\tilde{M}_{\text{int}}}$
and the terms that contribute to $h(I)$ are all of (roughly) the same magnitude in the worst case.
\end{sloppypar}

\begin{Theorem}\label{thm:alg_heuristic}
Let $\delta\in(0,1)$, $d\in\N\setminus\{1\}$, $I\subset\Z^d$, $|I|<\infty$, $T\ge\log_2\frac{d-1}{\delta}$, and $K\ge 3$ be fixed.
Moreover, we fix $\tilde{M}_{\text{int}}$ as stated in \eqref{def:tildeM_int}.
Then Algorithm~\ref{alg:heuristic} determines an exactly integrating rank-1 lattice for $I$ of size
$M< 4M_{\operatorname{lb}}$ with probability at least $1-2\delta^K$.
The arithmetic complexity of Algorithm~\ref{alg:heuristic} is bounded by terms in
\begin{equation}
\OO{K\,T\,d\,\left(\log{T}+h(I)\right)\log\tilde{M}_{\text{int}} + \tilde{M}_{\text{int}}^{0.525+\varepsilon}}
\label{eq:complexity_heuristic_general}
\end{equation}
for each fixed and arbitrarily small $\varepsilon>0$.
\end{Theorem}
\begin{proof}
We need to discuss the success probability of Algorithm~\ref{alg:check_exactness_construct_tilde_I} and we number the lattice size candidates $\tilde{M}$ in the order in which they appear.
Since we choose the initial $\tilde{M}_1$ 
and $T$ in accordance with the requirements of Corollary~\ref{thm:alg1}, i.e.,
\begin{align*}
\tilde{M}_1&:=\operatorname{\texttt{nextprime}}(2\,\tilde{M}_{\text{int}})=\underbrace{\frac{\tilde{M}_1}{M_{\operatorname{lb}}}}_{=:c_1}M_{\operatorname{lb}}\ge 2\,M_{\operatorname{lb}}\\
T&\ge\log_2\frac{d-1}{\delta}\ge\log_{c_1}\frac{d-1}{\delta},
\end{align*}
we can determine an upper bound $\delta_{d,c_1,T}\in(0,\delta]$ on the failure probability of Algorithm~\ref{alg:check_exactness_construct_tilde_I} using 
$$
T=\log_{c_1}\frac{d-1}{\delta_{d,c_1,T}},
$$
which yields that
Algorithm~\ref{alg:check_exactness_construct_tilde_I} succeeds with probability at least $1-(d-1)c_1^{-T}$.
Due to the fact that we iteratively restart Algorithm~\ref{alg:check_exactness_construct_tilde_I} in cases of failures, we 
exponentially decrease the failure probability due to the independence of each restart.
We apply Algorithm~\ref{alg:check_exactness_construct_tilde_I} up to $K$ times for the initial $\tilde{M}_1$, i.e.,
we observe a failure probability of at most $((d-1)c_1^{-T})^K$ for the first considered lattice size $\tilde{M}_1$.

The argumentation for the initial $\tilde{M}_1$ applies to $\tilde{M}_j:=\operatorname{\texttt{nextprime}}(\tilde{M}_{j-1}/2)$, $j=2,\ldots$, as long as $c_j:=\tilde{M}_j/M_{\operatorname{lb}}\ge 2$ holds true, i.e., we observe 
$((d-1)c_j^{-T})^K$ as upper bound on the failure probability. We define $j_0:=\max\{j\colon c_j\ge 2\}$, which yields
$\tilde{M}_{j_0}<4 M_{\operatorname{lb}}$ since otherwise we obtain $c_{j_0+1}\ge 2$, which is in contradiction to the definition of $j_0$. It remains to estimate the probability that Algorithm~\ref{alg:heuristic} terminates after Algorithm~\ref{alg:check_exactness_construct_tilde_I} succeeds for $\tilde{M}_{j_0}$.
\begin{align*}
&\mathbb{P}\left( \bigcap_{j=1}^{j_0}\{\textnormal{at least one out of $K$ calls of Algorithm~\ref{alg:construct_rand_reco_r1l_I} succeeds for $\tilde{M}_j$}\}\right)\\
&\quad = 1- \mathbb{P}\left( \bigcup_{j=1}^{j_0}\{\textnormal{all $K$ calls of Algorithm~\ref{alg:construct_rand_reco_r1l_I} fails for $\tilde{M}_j$}\}\right)\\
&\quad\ge 1-\sum_{j=1}^{j_0}\mathbb{P}\left(\{\textnormal{all $K$ calls of Algorithm~\ref{alg:construct_rand_reco_r1l_I} fails for $\tilde{M}_j$}\}\right)\\
&\quad\ge 1-(d-1)^K\sum_{j=1}^{j_0}c_j^{-TK}
\end{align*}
Due to Lemma~\ref{lem:primes}, we observe $\tilde{M}_{j+1}\le\frac{5}{7}\tilde{M}_j$ which yields $c_j\ge \frac{7}{5}c_{j+1}$, and thus, we get 
$$c_j\ge \left(\frac{7}{5}\right)^{j_0-j}c_{j_0}\ge \left(\frac{7}{5}\right)^{j_0-j}2,$$
for $j=1,\ldots,j_0$. Accordingly, we observe
\begin{align*}
\sum_{j=1}^{j_0}c_j^{-TK}&\le 2^{-TK}\sum_{j=1}^{j_0}\left(\left(\frac{7}{5}\right)^{j_0-j}\right)^{-TK}
\le 2^{-TK}\frac{7^{TK}}{7^{TK}-5^{TK}},
\end{align*}
which yields
\begin{equation*}
\mathbb{P}\left( \bigcap_{j=1}^{j_0}\{\textnormal{$K$ calls of Algorithm~\ref{alg:construct_rand_reco_r1l_I} succeeds for $\tilde{M}_j$ at least once}\}\right)
\ge 1-2\delta^K,
\end{equation*}
for $K\ge 3$.

As already discussed above, the computational complexity of Algorithm~\ref{alg:heuristic} is bounded by $K$ times $\log \tilde{M}_{\text{int}}$ times the complexity $T\,d\,\left(\log{T}+h(I)\right)$ of Algorithm~\ref{alg:construct_rand_reco_r1l_I} plus the effort of all $\operatorname{\texttt{nextprime}}$ calls, which cause an arithmetic complexity in $\OO{\tilde{M}_{\text{int}}^{0.525+\varepsilon}}$
for  arbitrary fixed $\varepsilon>0$.
\end{proof}

Let us consider the computational complexity of Algorithm~\ref{alg:heuristic} in more detail.
We start with Remark~\ref{rem:heuristic_complex_int}, which discusses the complexity of the probabilistic construction of exactly integrating rank-1 lattices.
However, as mentioned in the introduction, the crucial advantages of the presented approach appears when determining exactly integrating rank-1 lattices for difference sets $I=D(I')$ of frequency sets $I'$. Remark~\ref{rem:heuristic_complex_reco} investigates this specific approach in detail.

\begin{Remark}\label{rem:heuristic_complex_int}
Due to the fact that the result on prime gaps in~\cite{BaHaPi01} seems to be non-optimal, cf.\ the conjectures on prime gaps in \cite{Mi14}, and that the gap between consecutive primes is less than $1550$ for large but reasonable numbers%
\footnote{For details on numerically determined maximal prime gaps we refer to the \emph{Prime gap list project}, \url{https://primegap-list-project.github.io}, and to \emph{OEIS Foundation Inc.\ (2020), The On-Line Encyclopedia of Integer Sequences}, \url{https://oeis.org/A002386} and \url{https://oeis.org/A005250}.}
up to at least $2^{64}\approx 1.84\cdot 10^{19}$,
we
leave the part out of consideration in the following discussion. In other words, we just assume the length of prime gaps $q-p$ bounded in $\OO{p^{0.5-\varepsilon}}$ for some $\varepsilon>0$,  where $p,q$, $p<q$, are consecutive prime numbers.
Then, the total complexity of all $\operatorname{\texttt{nextprime}}$ calls is bounded in $\OO{\sqrt{\tilde{M}_{\text{int}}}}=\OO{\sqrt{\max(|I|,\max(I))}}$, cf.~\eqref{def:tildeM_int}, which is dominated by the following estimates of the first summand in \eqref{eq:complexity_heuristic_general} anyway when assuming $|I|\gtrsim \max(I)$.

So, let us consider $\OO{K\,T\,d\,\left(\log{T}+h(I)\right)\log\tilde{M}_{\text{int}}}$ with $h(I)$ as stated in \eqref{def:hI}.
Certainly, for the construction of exactly integrating rank-1 lattices for a frequency set $I$, we apply Algorithm~\ref{alg:construct_rand_reco_r1l_I}
which itself uses Algorithm~\ref{alg:check_exactness1} for checking the exact integration property. Accordingly, we can estimate
$h(I)\in\OO{|I|}$ and thus we observe a complexity of Algorithm~\ref{alg:heuristic} in $\OO{K\,T\,d\,\left(\log{T}+|I|\right)\log(\max(|I|,\max(I)))}$.
Under the assumption of Theorem~\ref{thm:alg_heuristic} and additionally assuming $|I|\gtrsim \max(I)$, $\delta\le 1/\sqrt{2}$, $K=3$, and $T:=\ceil{\log_2(d-1)-\log_2{\delta}}$, we observe a computational complexity in
$\OO{\log(d/\delta)\,d\,\left(\log\log{d/\delta}+|I|\right)\log|I|}$ in the worst case.
In addition, Algorithm~\ref{alg:heuristic} succeeds  (i.e., we have that $M<4M_{\operatorname{lb}}$ holds) with probability at least $1-\delta$ due to Theorem~\ref{thm:alg_heuristic}.

We stress the fact that one can slightly improve the computational complexity when using a similar heuristic strategy as indicated in Algorithm~\ref{alg:heuristic} in combination with the elimination approach presented in \cite{KuoMiNoNu19}, cf.\ \cite[Remark 24]{KuoMiNoNu19} for details.
This improvement would affect the logarithmic terms depending on $d$ and $\delta$, i.e., the strategy would lead to a non-probabilistic approach with complexity in $\OO{d|I|\log|I|}$ and resulting exactly integrating rank-1 lattices of sizes $M< 2 M_{\operatorname{lb}}$.
\end{Remark}

\begin{Remark}\label{rem:heuristic_complex_reco}\label{rem:construct_rand_reco_r1l_I_2}
In order to efficiently determine exactly integrating rank-1 lattices for difference sets $I=D(I')$ we switch to an approach that directly determines reconstructing rank-1 lattices for the frequency set $I'$. Clearly, both tasks are equivalent but the latter may save significant computational effort.
To this end, three noteworthy changes need to be made in Algorithm~\ref{alg:heuristic}:
\begin{itemize}
\item Instead of $I=D(I')$, we use $I'$ as input.
\item In line~\ref{alg:heuristic_nextprime1}, we choose the initial $\tilde{M}$ according to~\eqref{def:tildeM_reco}, i.e., 
$$\tilde{M}:=\operatorname{\texttt{nextprime}}(2\tilde{M}_{\text{reco}}(I'))=\operatorname{\texttt{nextprime}}(\max(|I'|^2,2N_{I'})).$$
\item Algorithm~\ref{alg:construct_rand_reco_r1l_I} uses Algorithm~\ref{alg:check_exactness} in order to check the reconstruction property for $I'_\ell$, $\ell=2,\ldots,d$.
\end{itemize}
We stress on the fact that one does not need to construct the difference set $D(I')$, since all calculations can be done using the frequency set $I'$ itself.
Nevertheless, the modified Algorithm~\ref{alg:heuristic} will determine an exactly integrating rank-1 lattice for $D(I')$ or, equivalently, a reconstructing rank-1 lattice for $I'$.

Using the aforementioned modifications, taking the discussion on $\operatorname{\texttt{nextprime}}$ in Remark~\ref{rem:heuristic_complex_int} into account,
and assuming $|I'|\gtrsim N_{I'}$, we observe a computational complexity of Algorithm~\ref{alg:heuristic} in
\begin{equation}
\OO{K\,T\,d\,\left(\log{T}+|I'|\log|I'|\right)\log \tilde{M}_{\text{reco}}}.
\label{eq:compl_alg5_reco}
\end{equation}
Under the assumption of Theorem~\ref{thm:alg_heuristic} and additionally assuming $|I'|\gtrsim N_{I'}$, $\delta\le 1/\sqrt{2}$, $K=3$, and $T:=\ceil{\log_2(d-1)-\log_2{\delta}}$, we observe a computational complexity in
\begin{equation*}
\OO{\log(d/\delta)\,d\,\left(\log\log{d/\delta}+|I'|\log|I'|\right)\log|I'|} \subset \OO{d\log^{1+\varepsilon}(d/\delta)\,|I'|\,\log^2|I'|},
\end{equation*}
$\varepsilon>0$ arbitrarily small, in the worst case.
As already mentioned in Remark~\ref{rem:heuristic_complex_int}, the parameters are chosen such that Algorithm~\ref{alg:heuristic} succeeds (i.e., we have that $M<4M_{\operatorname{lb}}$ holds) with probability at least $1-\delta$ due to Theorem~\ref{thm:alg_heuristic}.

We compare our result to the recent results in~\cite{KuoMiNoNu19}, for details we refer to row (j) in \cite[Table~1]{KuoMiNoNu19} listing the computational costs of different CBC approaches. Therein, the computational costs for (non-probabilistic) state of the art CBC constructions for reconstructing rank-1 lattices are specified as bounded in $\OO{d|D(I')|}$.
Combining these CBC approaches with a similar strategy as indicated in Algorithm~\ref{alg:heuristic} in order to determine rank-1 lattices of sizes
near $M_{\operatorname{lb}}$ again leads to almost the same complexity $\OO{d|D(I')|\log|I'|}$, up to an additional factor that is logarithmic in $|I'|$.
Clearly, the computational complexity of Algorithm~\ref{alg:heuristic} is slightly more dependent on the spatial dimension $d$. However, we observe that the dependence on the frequency sets differ widely. For known non-probabilistic algorithms, the computational complexity depends linearly on the cardinality of the difference set $D(I')$, whereas Algorithm~\ref{alg:heuristic} depends linearly on the cardinality of the original frequency set $I'$
up to some logarithmic factors.
For determining reconstructing rank-1 lattices, Algorithm~\ref{alg:heuristic} is most advantageous when $|D(I')|$ and $|I'|$ differ widely. 
Note that the inequalities $\sqrt{|D(I')|}\le |I'|\le |D(I')|$ hold in general,
where the first inequality is sharp up to some universal constant for specific frequency sets $I'$, cf.\ Section~\ref{ssec:numerics_anova} for an example.
According to the last considerations, we proved a significantly lower computational complexity for Algorithm~\ref{alg:heuristic} compared to known non-probabilistic approaches in cases where $|I'|\ll |D(I')|$ holds.
\end{Remark}

\section{Numerics}\label{sec:numerics}

Since we expect significant numerical advantages of the presented probabilistic lattice search, cf. Algorithm~\ref{alg:construct_rand_reco_r1l_I_2} and Remark~\ref{rem:heuristic_complex_reco}, when searching for reconstructing rank-1 lattices, the following numerical tests focus on applications that allow for exploiting this fact.
We repeated each lattice search ten times and plotted the average values.
In addition, we illustrated the whole ranges of the different outcomes as bars.
In more detail, we used implementations of the presented algorithms in Julia 1.3.0.
The time measurements were performed on a computer with
2 x 6-core Intel Xeon CPU X5690 (3,47 GHz), 144 GB RAM, using a single thread.
Furthermore, we would like to point out the choice of fixed parameters $T=100$ and $K=5$ in all numerical tests, which
leads to failure probabilities less than $\delta'=(d-1)^K2^{1-TK}\le (d-1)^52^{-499}<10^{-133}$ for dimensions $d\le 2000$ according to Theorem~\ref{thm:alg_heuristic}, i.e.,
the resulting rank-1 lattice sizes are pretty sure less than $4M_{\operatorname{lb}}$.
Since $\delta'<10^{-120}$ even for all dimensions $d\le 10^6$, we do not see the necessity of varying
the parameters when the dimension $d$ is changing in our numerical tests.

\subsection{Integration lattices for functions of a specific ANOVA property}\label{ssec:numerics_anova}

Let us assume that the superposition dimension, cf.~\cite{CaMoOw97}, of the considered periodic function $f$ is only two, which means that all terms of effective dimension at most two within the
ANOVA decomposition covers almost entirely the function $f$, cf.~\cite{PoSchmi20} for details on ANOVA decompositions of Fourier series.
From a Fourier point of view, the function $f$ is, up to a small error, reconstructable using frequencies $\boldk\in\Z^d$ with at most two nonzero entries.
Roughly speaking, the large Fourier coefficients $c_\boldk(f)$ belong to the $\left(\begin{array}{c}d\\2\end{array}\right)$ two dimensional planes that are spanned by two axis. We join all these planes and define
\begin{align*}
\A_d&:=\bigcup_{\substack{1\le j\le d-1\\j<l\le d}}\A_{j,l}^{(d)},\\
\A_{j,l}^{(d)}&:=\underbrace{\{0\}\times\dots\{0\}}_{\subset\Z^{j-1}}\times\Z\times\underbrace{\{0\}\times\dots\times\{0\}}_{\subset\Z^{l-j-1}}\times\Z\times\underbrace{\{0\}\times\dots\times\{0\}}_{\subset\Z^{d-l}}.
\end{align*}
In addition, we assume that the function $f$ is isotropically smooth, which means, that the absolute values of the Fourier coefficients of $f$ are uniformly low outside of a cube $[-N,N]^d$, $N\in\N$.

\begin{figure}[tb]
\subfloat[$I_8^3$\label{fig:ANOVA_I}]{
\begin{tikzpicture}
\begin{axis}[axis background/.style={fill=white},
every axis/.append style={font=\footnotesize},
width=0.32\textwidth,
height=0.32\textwidth,
enlargelimits=false,
enlargelimits=false,
clip=false,
view={15}{15},
grid=major,
plot box ratio = 1 1 1,
clip mode=individual,
tickwidth=0pt,
z buffer=sort,
xmin=-16,xmax=16,
ymin=-16,ymax=16,
zmin=-16, zmax=16,
ytick={-10,10}
]
\addplot3+[only marks, mark size=0.5pt, mark=*, solid, ball color=black!75, mark options={black!75, draw=black}] coordinates{(0,0,0) (1,0,0) (-1,0,0) (2,0,0) (-2,0,0) (3,0,0) (-3,0,0) (4,0,0) (-4,0,0) (5,0,0) (-5,0,0) (6,0,0) (-6,0,0) (7,0,0) (-7,0,0) (8,0,0) (-8,0,0) (1,1,0) (-1,1,0) (1,-1,0) (-1,-1,0) (1,2,0) (-1,2,0) (1,-2,0) (-1,-2,0) (1,3,0) (-1,3,0) (1,-3,0) (-1,-3,0) (1,4,0) (-1,4,0) (1,-4,0) (-1,-4,0) (1,5,0) (-1,5,0) (1,-5,0) (-1,-5,0) (1,6,0) (-1,6,0) (1,-6,0) (-1,-6,0) (1,7,0) (-1,7,0) (1,-7,0) (-1,-7,0) (1,8,0) (-1,8,0) (1,-8,0) (-1,-8,0) (2,1,0) (-2,1,0) (2,-1,0) (-2,-1,0) (2,2,0) (-2,2,0) (2,-2,0) (-2,-2,0) (2,3,0) (-2,3,0) (2,-3,0) (-2,-3,0) (2,4,0) (-2,4,0) (2,-4,0) (-2,-4,0) (2,5,0) (-2,5,0) (2,-5,0) (-2,-5,0) (2,6,0) (-2,6,0) (2,-6,0) (-2,-6,0) (2,7,0) (-2,7,0) (2,-7,0) (-2,-7,0) (2,8,0) (-2,8,0) (2,-8,0) (-2,-8,0) (3,1,0) (-3,1,0) (3,-1,0) (-3,-1,0) (3,2,0) (-3,2,0) (3,-2,0) (-3,-2,0) (3,3,0) (-3,3,0) (3,-3,0) (-3,-3,0) (3,4,0) (-3,4,0) (3,-4,0) (-3,-4,0) (3,5,0) (-3,5,0) (3,-5,0) (-3,-5,0) (3,6,0) (-3,6,0) (3,-6,0) (-3,-6,0) (3,7,0) (-3,7,0) (3,-7,0) (-3,-7,0) (3,8,0) (-3,8,0) (3,-8,0) (-3,-8,0) (4,1,0) (-4,1,0) (4,-1,0) (-4,-1,0) (4,2,0) (-4,2,0) (4,-2,0) (-4,-2,0) (4,3,0) (-4,3,0) (4,-3,0) (-4,-3,0) (4,4,0) (-4,4,0) (4,-4,0) (-4,-4,0) (4,5,0) (-4,5,0) (4,-5,0) (-4,-5,0) (4,6,0) (-4,6,0) (4,-6,0) (-4,-6,0) (4,7,0) (-4,7,0) (4,-7,0) (-4,-7,0) (4,8,0) (-4,8,0) (4,-8,0) (-4,-8,0) (5,1,0) (-5,1,0) (5,-1,0) (-5,-1,0) (5,2,0) (-5,2,0) (5,-2,0) (-5,-2,0) (5,3,0) (-5,3,0) (5,-3,0) (-5,-3,0) (5,4,0) (-5,4,0) (5,-4,0) (-5,-4,0) (5,5,0) (-5,5,0) (5,-5,0) (-5,-5,0) (5,6,0) (-5,6,0) (5,-6,0) (-5,-6,0) (5,7,0) (-5,7,0) (5,-7,0) (-5,-7,0) (5,8,0) (-5,8,0) (5,-8,0) (-5,-8,0) (6,1,0) (-6,1,0) (6,-1,0) (-6,-1,0) (6,2,0) (-6,2,0) (6,-2,0) (-6,-2,0) (6,3,0) (-6,3,0) (6,-3,0) (-6,-3,0) (6,4,0) (-6,4,0) (6,-4,0) (-6,-4,0) (6,5,0) (-6,5,0) (6,-5,0) (-6,-5,0) (6,6,0) (-6,6,0) (6,-6,0) (-6,-6,0) (6,7,0) (-6,7,0) (6,-7,0) (-6,-7,0) (6,8,0) (-6,8,0) (6,-8,0) (-6,-8,0) (7,1,0) (-7,1,0) (7,-1,0) (-7,-1,0) (7,2,0) (-7,2,0) (7,-2,0) (-7,-2,0) (7,3,0) (-7,3,0) (7,-3,0) (-7,-3,0) (7,4,0) (-7,4,0) (7,-4,0) (-7,-4,0) (7,5,0) (-7,5,0) (7,-5,0) (-7,-5,0) (7,6,0) (-7,6,0) (7,-6,0) (-7,-6,0) (7,7,0) (-7,7,0) (7,-7,0) (-7,-7,0) (7,8,0) (-7,8,0) (7,-8,0) (-7,-8,0) (8,1,0) (-8,1,0) (8,-1,0) (-8,-1,0) (8,2,0) (-8,2,0) (8,-2,0) (-8,-2,0) (8,3,0) (-8,3,0) (8,-3,0) (-8,-3,0) (8,4,0) (-8,4,0) (8,-4,0) (-8,-4,0) (8,5,0) (-8,5,0) (8,-5,0) (-8,-5,0) (8,6,0) (-8,6,0) (8,-6,0) (-8,-6,0) (8,7,0) (-8,7,0) (8,-7,0) (-8,-7,0) (8,8,0) (-8,8,0) (8,-8,0) (-8,-8,0) (1,0,1) (-1,0,1) (1,0,-1) (-1,0,-1) (1,0,2) (-1,0,2) (1,0,-2) (-1,0,-2) (1,0,3) (-1,0,3) (1,0,-3) (-1,0,-3) (1,0,4) (-1,0,4) (1,0,-4) (-1,0,-4) (1,0,5) (-1,0,5) (1,0,-5) (-1,0,-5) (1,0,6) (-1,0,6) (1,0,-6) (-1,0,-6) (1,0,7) (-1,0,7) (1,0,-7) (-1,0,-7) (1,0,8) (-1,0,8) (1,0,-8) (-1,0,-8) (2,0,1) (-2,0,1) (2,0,-1) (-2,0,-1) (2,0,2) (-2,0,2) (2,0,-2) (-2,0,-2) (2,0,3) (-2,0,3) (2,0,-3) (-2,0,-3) (2,0,4) (-2,0,4) (2,0,-4) (-2,0,-4) (2,0,5) (-2,0,5) (2,0,-5) (-2,0,-5) (2,0,6) (-2,0,6) (2,0,-6) (-2,0,-6) (2,0,7) (-2,0,7) (2,0,-7) (-2,0,-7) (2,0,8) (-2,0,8) (2,0,-8) (-2,0,-8) (3,0,1) (-3,0,1) (3,0,-1) (-3,0,-1) (3,0,2) (-3,0,2) (3,0,-2) (-3,0,-2) (3,0,3) (-3,0,3) (3,0,-3) (-3,0,-3) (3,0,4) (-3,0,4) (3,0,-4) (-3,0,-4) (3,0,5) (-3,0,5) (3,0,-5) (-3,0,-5) (3,0,6) (-3,0,6) (3,0,-6) (-3,0,-6) (3,0,7) (-3,0,7) (3,0,-7) (-3,0,-7) (3,0,8) (-3,0,8) (3,0,-8) (-3,0,-8) (4,0,1) (-4,0,1) (4,0,-1) (-4,0,-1) (4,0,2) (-4,0,2) (4,0,-2) (-4,0,-2) (4,0,3) (-4,0,3) (4,0,-3) (-4,0,-3) (4,0,4) (-4,0,4) (4,0,-4) (-4,0,-4) (4,0,5) (-4,0,5) (4,0,-5) (-4,0,-5) (4,0,6) (-4,0,6) (4,0,-6) (-4,0,-6) (4,0,7) (-4,0,7) (4,0,-7) (-4,0,-7) (4,0,8) (-4,0,8) (4,0,-8) (-4,0,-8) (5,0,1) (-5,0,1) (5,0,-1) (-5,0,-1) (5,0,2) (-5,0,2) (5,0,-2) (-5,0,-2) (5,0,3) (-5,0,3) (5,0,-3) (-5,0,-3) (5,0,4) (-5,0,4) (5,0,-4) (-5,0,-4) (5,0,5) (-5,0,5) (5,0,-5) (-5,0,-5) (5,0,6) (-5,0,6) (5,0,-6) (-5,0,-6) (5,0,7) (-5,0,7) (5,0,-7) (-5,0,-7) (5,0,8) (-5,0,8) (5,0,-8) (-5,0,-8) (6,0,1) (-6,0,1) (6,0,-1) (-6,0,-1) (6,0,2) (-6,0,2) (6,0,-2) (-6,0,-2) (6,0,3) (-6,0,3) (6,0,-3) (-6,0,-3) (6,0,4) (-6,0,4) (6,0,-4) (-6,0,-4) (6,0,5) (-6,0,5) (6,0,-5) (-6,0,-5) (6,0,6) (-6,0,6) (6,0,-6) (-6,0,-6) (6,0,7) (-6,0,7) (6,0,-7) (-6,0,-7) (6,0,8) (-6,0,8) (6,0,-8) (-6,0,-8) (7,0,1) (-7,0,1) (7,0,-1) (-7,0,-1) (7,0,2) (-7,0,2) (7,0,-2) (-7,0,-2) (7,0,3) (-7,0,3) (7,0,-3) (-7,0,-3) (7,0,4) (-7,0,4) (7,0,-4) (-7,0,-4) (7,0,5) (-7,0,5) (7,0,-5) (-7,0,-5) (7,0,6) (-7,0,6) (7,0,-6) (-7,0,-6) (7,0,7) (-7,0,7) (7,0,-7) (-7,0,-7) (7,0,8) (-7,0,8) (7,0,-8) (-7,0,-8) (8,0,1) (-8,0,1) (8,0,-1) (-8,0,-1) (8,0,2) (-8,0,2) (8,0,-2) (-8,0,-2) (8,0,3) (-8,0,3) (8,0,-3) (-8,0,-3) (8,0,4) (-8,0,4) (8,0,-4) (-8,0,-4) (8,0,5) (-8,0,5) (8,0,-5) (-8,0,-5) (8,0,6) (-8,0,6) (8,0,-6) (-8,0,-6) (8,0,7) (-8,0,7) (8,0,-7) (-8,0,-7) (8,0,8) (-8,0,8) (8,0,-8) (-8,0,-8) (0,1,0) (0,-1,0) (0,2,0) (0,-2,0) (0,3,0) (0,-3,0) (0,4,0) (0,-4,0) (0,5,0) (0,-5,0) (0,6,0) (0,-6,0) (0,7,0) (0,-7,0) (0,8,0) (0,-8,0) (0,1,1) (0,-1,1) (0,1,-1) (0,-1,-1) (0,1,2) (0,-1,2) (0,1,-2) (0,-1,-2) (0,1,3) (0,-1,3) (0,1,-3) (0,-1,-3) (0,1,4) (0,-1,4) (0,1,-4) (0,-1,-4) (0,1,5) (0,-1,5) (0,1,-5) (0,-1,-5) (0,1,6) (0,-1,6) (0,1,-6) (0,-1,-6) (0,1,7) (0,-1,7) (0,1,-7) (0,-1,-7) (0,1,8) (0,-1,8) (0,1,-8) (0,-1,-8) (0,2,1) (0,-2,1) (0,2,-1) (0,-2,-1) (0,2,2) (0,-2,2) (0,2,-2) (0,-2,-2) (0,2,3) (0,-2,3) (0,2,-3) (0,-2,-3) (0,2,4) (0,-2,4) (0,2,-4) (0,-2,-4) (0,2,5) (0,-2,5) (0,2,-5) (0,-2,-5) (0,2,6) (0,-2,6) (0,2,-6) (0,-2,-6) (0,2,7) (0,-2,7) (0,2,-7) (0,-2,-7) (0,2,8) (0,-2,8) (0,2,-8) (0,-2,-8) (0,3,1) (0,-3,1) (0,3,-1) (0,-3,-1) (0,3,2) (0,-3,2) (0,3,-2) (0,-3,-2) (0,3,3) (0,-3,3) (0,3,-3) (0,-3,-3) (0,3,4) (0,-3,4) (0,3,-4) (0,-3,-4) (0,3,5) (0,-3,5) (0,3,-5) (0,-3,-5) (0,3,6) (0,-3,6) (0,3,-6) (0,-3,-6) (0,3,7) (0,-3,7) (0,3,-7) (0,-3,-7) (0,3,8) (0,-3,8) (0,3,-8) (0,-3,-8) (0,4,1) (0,-4,1) (0,4,-1) (0,-4,-1) (0,4,2) (0,-4,2) (0,4,-2) (0,-4,-2) (0,4,3) (0,-4,3) (0,4,-3) (0,-4,-3) (0,4,4) (0,-4,4) (0,4,-4) (0,-4,-4) (0,4,5) (0,-4,5) (0,4,-5) (0,-4,-5) (0,4,6) (0,-4,6) (0,4,-6) (0,-4,-6) (0,4,7) (0,-4,7) (0,4,-7) (0,-4,-7) (0,4,8) (0,-4,8) (0,4,-8) (0,-4,-8) (0,5,1) (0,-5,1) (0,5,-1) (0,-5,-1) (0,5,2) (0,-5,2) (0,5,-2) (0,-5,-2) (0,5,3) (0,-5,3) (0,5,-3) (0,-5,-3) (0,5,4) (0,-5,4) (0,5,-4) (0,-5,-4) (0,5,5) (0,-5,5) (0,5,-5) (0,-5,-5) (0,5,6) (0,-5,6) (0,5,-6) (0,-5,-6) (0,5,7) (0,-5,7) (0,5,-7) (0,-5,-7) (0,5,8) (0,-5,8) (0,5,-8) (0,-5,-8) (0,6,1) (0,-6,1) (0,6,-1) (0,-6,-1) (0,6,2) (0,-6,2) (0,6,-2) (0,-6,-2) (0,6,3) (0,-6,3) (0,6,-3) (0,-6,-3) (0,6,4) (0,-6,4) (0,6,-4) (0,-6,-4) (0,6,5) (0,-6,5) (0,6,-5) (0,-6,-5) (0,6,6) (0,-6,6) (0,6,-6) (0,-6,-6) (0,6,7) (0,-6,7) (0,6,-7) (0,-6,-7) (0,6,8) (0,-6,8) (0,6,-8) (0,-6,-8) (0,7,1) (0,-7,1) (0,7,-1) (0,-7,-1) (0,7,2) (0,-7,2) (0,7,-2) (0,-7,-2) (0,7,3) (0,-7,3) (0,7,-3) (0,-7,-3) (0,7,4) (0,-7,4) (0,7,-4) (0,-7,-4) (0,7,5) (0,-7,5) (0,7,-5) (0,-7,-5) (0,7,6) (0,-7,6) (0,7,-6) (0,-7,-6) (0,7,7) (0,-7,7) (0,7,-7) (0,-7,-7) (0,7,8) (0,-7,8) (0,7,-8) (0,-7,-8) (0,8,1) (0,-8,1) (0,8,-1) (0,-8,-1) (0,8,2) (0,-8,2) (0,8,-2) (0,-8,-2) (0,8,3) (0,-8,3) (0,8,-3) (0,-8,-3) (0,8,4) (0,-8,4) (0,8,-4) (0,-8,-4) (0,8,5) (0,-8,5) (0,8,-5) (0,-8,-5) (0,8,6) (0,-8,6) (0,8,-6) (0,-8,-6) (0,8,7) (0,-8,7) (0,8,-7) (0,-8,-7) (0,8,8) (0,-8,8) (0,8,-8) (0,-8,-8) (0,0,1) (0,0,-1) (0,0,2) (0,0,-2) (0,0,3) (0,0,-3) (0,0,4) (0,0,-4) (0,0,5) (0,0,-5) (0,0,6) (0,0,-6) (0,0,7) (0,0,-7) (0,0,8) (0,0,-8)};
\end{axis}
\end{tikzpicture}
}
\subfloat[$\tilde{I}_8^3$\label{fig:axis_cross}]{
\begin{tikzpicture}
\begin{axis}[axis background/.style={fill=white},
every axis/.append style={font=\footnotesize},
width=0.32\textwidth,
height=0.32\textwidth,
enlargelimits=false,
enlargelimits=false,
clip=false,
view={15}{15},
grid=major,
plot box ratio = 1 1 1,
clip mode=individual,
tickwidth=0pt,
z buffer=sort,
xmin=-16,xmax=16,
ymin=-16,ymax=16,
zmin=-16, zmax=16,
ytick={-10,10}
]
\addplot3+[only marks, mark size=0.5pt, mark=*, solid, ball color=black!75, mark options={black!75, draw=black}] coordinates{
(0,0,0) (1,0,0) (-1,0,0) (2,0,0) (-2,0,0) (3,0,0) (-3,0,0) (4,0,0) (-4,0,0) (5,0,0) (-5,0,0) (6,0,0) (-6,0,0) (7,0,0) (-7,0,0) (8,0,0) (-8,0,0) (0,1,0) (0,-1,0) (0,2,0) (0,-2,0) (0,3,0) (0,-3,0) (0,4,0) (0,-4,0) (0,5,0) (0,-5,0) (0,6,0) (0,-6,0) (0,7,0) (0,-7,0) (0,8,0) (0,-8,0) (0,0,1) (0,0,-1) (0,0,2) (0,0,-2) (0,0,3) (0,0,-3) (0,0,4) (0,0,-4) (0,0,5) (0,0,-5) (0,0,6) (0,0,-6) (0,0,7) (0,0,-7) (0,0,8) (0,0,-8)};
\end{axis}
\end{tikzpicture}
}
\subfloat[$D(\tilde{I}_8^3)$\label{fig:differenceset_axis_cross}]{
\begin{tikzpicture}
\begin{axis}[axis background/.style={fill=white},
every axis/.append style={font=\footnotesize},
width=0.32\textwidth,
height=0.32\textwidth,
enlargelimits=false,
enlargelimits=false,
clip=false,
view={15}{15},
grid=major,
plot box ratio = 1 1 1,
clip mode=individual,
tickwidth=0pt,
z buffer=sort,
xmin=-16,xmax=16,
ymin=-16,ymax=16,
zmin=-16, zmax=16,
ytick={-10,10}
]
\addplot3+[only marks, mark size=0.5pt, mark=*, solid, ball color=black!75, mark options={black!75, draw=black}] coordinates{
(0,0,0) (-1,0,0) (1,0,0) (-2,0,0) (2,0,0) (-3,0,0) (3,0,0) (-4,0,0) (4,0,0) (-5,0,0) (5,0,0) (-6,0,0) (6,0,0) (-7,0,0) (7,0,0) (-8,0,0) (8,0,0) (0,-1,0) (0,1,0) (0,-2,0) (0,2,0) (0,-3,0) (0,3,0) (0,-4,0) (0,4,0) (0,-5,0) (0,5,0) (0,-6,0) (0,6,0) (0,-7,0) (0,7,0) (0,-8,0) (0,8,0) (0,0,-1) (0,0,1) (0,0,-2) (0,0,2) (0,0,-3) (0,0,3) (0,0,-4) (0,0,4) (0,0,-5) (0,0,5) (0,0,-6) (0,0,6) (0,0,-7) (0,0,7) (0,0,-8) (0,0,8) (9,0,0) (1,-1,0) (1,1,0) (1,-2,0) (1,2,0) (1,-3,0) (1,3,0) (1,-4,0) (1,4,0) (1,-5,0) (1,5,0) (1,-6,0) (1,6,0) (1,-7,0) (1,7,0) (1,-8,0) (1,8,0) (1,0,-1) (1,0,1) (1,0,-2) (1,0,2) (1,0,-3) (1,0,3) (1,0,-4) (1,0,4) (1,0,-5) (1,0,5) (1,0,-6) (1,0,6) (1,0,-7) (1,0,7) (1,0,-8) (1,0,8) (-9,0,0) (-1,-1,0) (-1,1,0) (-1,-2,0) (-1,2,0) (-1,-3,0) (-1,3,0) (-1,-4,0) (-1,4,0) (-1,-5,0) (-1,5,0) (-1,-6,0) (-1,6,0) (-1,-7,0) (-1,7,0) (-1,-8,0) (-1,8,0) (-1,0,-1) (-1,0,1) (-1,0,-2) (-1,0,2) (-1,0,-3) (-1,0,3) (-1,0,-4) (-1,0,4) (-1,0,-5) (-1,0,5) (-1,0,-6) (-1,0,6) (-1,0,-7) (-1,0,7) (-1,0,-8) (-1,0,8) (10,0,0) (2,-1,0) (2,1,0) (2,-2,0) (2,2,0) (2,-3,0) (2,3,0) (2,-4,0) (2,4,0) (2,-5,0) (2,5,0) (2,-6,0) (2,6,0) (2,-7,0) (2,7,0) (2,-8,0) (2,8,0) (2,0,-1) (2,0,1) (2,0,-2) (2,0,2) (2,0,-3) (2,0,3) (2,0,-4) (2,0,4) (2,0,-5) (2,0,5) (2,0,-6) (2,0,6) (2,0,-7) (2,0,7) (2,0,-8) (2,0,8) (-10,0,0) (-2,-1,0) (-2,1,0) (-2,-2,0) (-2,2,0) (-2,-3,0) (-2,3,0) (-2,-4,0) (-2,4,0) (-2,-5,0) (-2,5,0) (-2,-6,0) (-2,6,0) (-2,-7,0) (-2,7,0) (-2,-8,0) (-2,8,0) (-2,0,-1) (-2,0,1) (-2,0,-2) (-2,0,2) (-2,0,-3) (-2,0,3) (-2,0,-4) (-2,0,4) (-2,0,-5) (-2,0,5) (-2,0,-6) (-2,0,6) (-2,0,-7) (-2,0,7) (-2,0,-8) (-2,0,8) (11,0,0) (3,-1,0) (3,1,0) (3,-2,0) (3,2,0) (3,-3,0) (3,3,0) (3,-4,0) (3,4,0) (3,-5,0) (3,5,0) (3,-6,0) (3,6,0) (3,-7,0) (3,7,0) (3,-8,0) (3,8,0) (3,0,-1) (3,0,1) (3,0,-2) (3,0,2) (3,0,-3) (3,0,3) (3,0,-4) (3,0,4) (3,0,-5) (3,0,5) (3,0,-6) (3,0,6) (3,0,-7) (3,0,7) (3,0,-8) (3,0,8) (-11,0,0) (-3,-1,0) (-3,1,0) (-3,-2,0) (-3,2,0) (-3,-3,0) (-3,3,0) (-3,-4,0) (-3,4,0) (-3,-5,0) (-3,5,0) (-3,-6,0) (-3,6,0) (-3,-7,0) (-3,7,0) (-3,-8,0) (-3,8,0) (-3,0,-1) (-3,0,1) (-3,0,-2) (-3,0,2) (-3,0,-3) (-3,0,3) (-3,0,-4) (-3,0,4) (-3,0,-5) (-3,0,5) (-3,0,-6) (-3,0,6) (-3,0,-7) (-3,0,7) (-3,0,-8) (-3,0,8) (12,0,0) (4,-1,0) (4,1,0) (4,-2,0) (4,2,0) (4,-3,0) (4,3,0) (4,-4,0) (4,4,0) (4,-5,0) (4,5,0) (4,-6,0) (4,6,0) (4,-7,0) (4,7,0) (4,-8,0) (4,8,0) (4,0,-1) (4,0,1) (4,0,-2) (4,0,2) (4,0,-3) (4,0,3) (4,0,-4) (4,0,4) (4,0,-5) (4,0,5) (4,0,-6) (4,0,6) (4,0,-7) (4,0,7) (4,0,-8) (4,0,8) (-12,0,0) (-4,-1,0) (-4,1,0) (-4,-2,0) (-4,2,0) (-4,-3,0) (-4,3,0) (-4,-4,0) (-4,4,0) (-4,-5,0) (-4,5,0) (-4,-6,0) (-4,6,0) (-4,-7,0) (-4,7,0) (-4,-8,0) (-4,8,0) (-4,0,-1) (-4,0,1) (-4,0,-2) (-4,0,2) (-4,0,-3) (-4,0,3) (-4,0,-4) (-4,0,4) (-4,0,-5) (-4,0,5) (-4,0,-6) (-4,0,6) (-4,0,-7) (-4,0,7) (-4,0,-8) (-4,0,8) (13,0,0) (5,-1,0) (5,1,0) (5,-2,0) (5,2,0) (5,-3,0) (5,3,0) (5,-4,0) (5,4,0) (5,-5,0) (5,5,0) (5,-6,0) (5,6,0) (5,-7,0) (5,7,0) (5,-8,0) (5,8,0) (5,0,-1) (5,0,1) (5,0,-2) (5,0,2) (5,0,-3) (5,0,3) (5,0,-4) (5,0,4) (5,0,-5) (5,0,5) (5,0,-6) (5,0,6) (5,0,-7) (5,0,7) (5,0,-8) (5,0,8) (-13,0,0) (-5,-1,0) (-5,1,0) (-5,-2,0) (-5,2,0) (-5,-3,0) (-5,3,0) (-5,-4,0) (-5,4,0) (-5,-5,0) (-5,5,0) (-5,-6,0) (-5,6,0) (-5,-7,0) (-5,7,0) (-5,-8,0) (-5,8,0) (-5,0,-1) (-5,0,1) (-5,0,-2) (-5,0,2) (-5,0,-3) (-5,0,3) (-5,0,-4) (-5,0,4) (-5,0,-5) (-5,0,5) (-5,0,-6) (-5,0,6) (-5,0,-7) (-5,0,7) (-5,0,-8) (-5,0,8) (14,0,0) (6,-1,0) (6,1,0) (6,-2,0) (6,2,0) (6,-3,0) (6,3,0) (6,-4,0) (6,4,0) (6,-5,0) (6,5,0) (6,-6,0) (6,6,0) (6,-7,0) (6,7,0) (6,-8,0) (6,8,0) (6,0,-1) (6,0,1) (6,0,-2) (6,0,2) (6,0,-3) (6,0,3) (6,0,-4) (6,0,4) (6,0,-5) (6,0,5) (6,0,-6) (6,0,6) (6,0,-7) (6,0,7) (6,0,-8) (6,0,8) (-14,0,0) (-6,-1,0) (-6,1,0) (-6,-2,0) (-6,2,0) (-6,-3,0) (-6,3,0) (-6,-4,0) (-6,4,0) (-6,-5,0) (-6,5,0) (-6,-6,0) (-6,6,0) (-6,-7,0) (-6,7,0) (-6,-8,0) (-6,8,0) (-6,0,-1) (-6,0,1) (-6,0,-2) (-6,0,2) (-6,0,-3) (-6,0,3) (-6,0,-4) (-6,0,4) (-6,0,-5) (-6,0,5) (-6,0,-6) (-6,0,6) (-6,0,-7) (-6,0,7) (-6,0,-8) (-6,0,8) (15,0,0) (7,-1,0) (7,1,0) (7,-2,0) (7,2,0) (7,-3,0) (7,3,0) (7,-4,0) (7,4,0) (7,-5,0) (7,5,0) (7,-6,0) (7,6,0) (7,-7,0) (7,7,0) (7,-8,0) (7,8,0) (7,0,-1) (7,0,1) (7,0,-2) (7,0,2) (7,0,-3) (7,0,3) (7,0,-4) (7,0,4) (7,0,-5) (7,0,5) (7,0,-6) (7,0,6) (7,0,-7) (7,0,7) (7,0,-8) (7,0,8) (-15,0,0) (-7,-1,0) (-7,1,0) (-7,-2,0) (-7,2,0) (-7,-3,0) (-7,3,0) (-7,-4,0) (-7,4,0) (-7,-5,0) (-7,5,0) (-7,-6,0) (-7,6,0) (-7,-7,0) (-7,7,0) (-7,-8,0) (-7,8,0) (-7,0,-1) (-7,0,1) (-7,0,-2) (-7,0,2) (-7,0,-3) (-7,0,3) (-7,0,-4) (-7,0,4) (-7,0,-5) (-7,0,5) (-7,0,-6) (-7,0,6) (-7,0,-7) (-7,0,7) (-7,0,-8) (-7,0,8) (16,0,0) (8,-1,0) (8,1,0) (8,-2,0) (8,2,0) (8,-3,0) (8,3,0) (8,-4,0) (8,4,0) (8,-5,0) (8,5,0) (8,-6,0) (8,6,0) (8,-7,0) (8,7,0) (8,-8,0) (8,8,0) (8,0,-1) (8,0,1) (8,0,-2) (8,0,2) (8,0,-3) (8,0,3) (8,0,-4) (8,0,4) (8,0,-5) (8,0,5) (8,0,-6) (8,0,6) (8,0,-7) (8,0,7) (8,0,-8) (8,0,8) (-16,0,0) (-8,-1,0) (-8,1,0) (-8,-2,0) (-8,2,0) (-8,-3,0) (-8,3,0) (-8,-4,0) (-8,4,0) (-8,-5,0) (-8,5,0) (-8,-6,0) (-8,6,0) (-8,-7,0) (-8,7,0) (-8,-8,0) (-8,8,0) (-8,0,-1) (-8,0,1) (-8,0,-2) (-8,0,2) (-8,0,-3) (-8,0,3) (-8,0,-4) (-8,0,4) (-8,0,-5) (-8,0,5) (-8,0,-6) (-8,0,6) (-8,0,-7) (-8,0,7) (-8,0,-8) (-8,0,8) (0,9,0) (0,1,-1) (0,1,1) (0,1,-2) (0,1,2) (0,1,-3) (0,1,3) (0,1,-4) (0,1,4) (0,1,-5) (0,1,5) (0,1,-6) (0,1,6) (0,1,-7) (0,1,7) (0,1,-8) (0,1,8) (0,-9,0) (0,-1,-1) (0,-1,1) (0,-1,-2) (0,-1,2) (0,-1,-3) (0,-1,3) (0,-1,-4) (0,-1,4) (0,-1,-5) (0,-1,5) (0,-1,-6) (0,-1,6) (0,-1,-7) (0,-1,7) (0,-1,-8) (0,-1,8) (0,10,0) (0,2,-1) (0,2,1) (0,2,-2) (0,2,2) (0,2,-3) (0,2,3) (0,2,-4) (0,2,4) (0,2,-5) (0,2,5) (0,2,-6) (0,2,6) (0,2,-7) (0,2,7) (0,2,-8) (0,2,8) (0,-10,0) (0,-2,-1) (0,-2,1) (0,-2,-2) (0,-2,2) (0,-2,-3) (0,-2,3) (0,-2,-4) (0,-2,4) (0,-2,-5) (0,-2,5) (0,-2,-6) (0,-2,6) (0,-2,-7) (0,-2,7) (0,-2,-8) (0,-2,8) (0,11,0) (0,3,-1) (0,3,1) (0,3,-2) (0,3,2) (0,3,-3) (0,3,3) (0,3,-4) (0,3,4) (0,3,-5) (0,3,5) (0,3,-6) (0,3,6) (0,3,-7) (0,3,7) (0,3,-8) (0,3,8) (0,-11,0) (0,-3,-1) (0,-3,1) (0,-3,-2) (0,-3,2) (0,-3,-3) (0,-3,3) (0,-3,-4) (0,-3,4) (0,-3,-5) (0,-3,5) (0,-3,-6) (0,-3,6) (0,-3,-7) (0,-3,7) (0,-3,-8) (0,-3,8) (0,12,0) (0,4,-1) (0,4,1) (0,4,-2) (0,4,2) (0,4,-3) (0,4,3) (0,4,-4) (0,4,4) (0,4,-5) (0,4,5) (0,4,-6) (0,4,6) (0,4,-7) (0,4,7) (0,4,-8) (0,4,8) (0,-12,0) (0,-4,-1) (0,-4,1) (0,-4,-2) (0,-4,2) (0,-4,-3) (0,-4,3) (0,-4,-4) (0,-4,4) (0,-4,-5) (0,-4,5) (0,-4,-6) (0,-4,6) (0,-4,-7) (0,-4,7) (0,-4,-8) (0,-4,8) (0,13,0) (0,5,-1) (0,5,1) (0,5,-2) (0,5,2) (0,5,-3) (0,5,3) (0,5,-4) (0,5,4) (0,5,-5) (0,5,5) (0,5,-6) (0,5,6) (0,5,-7) (0,5,7) (0,5,-8) (0,5,8) (0,-13,0) (0,-5,-1) (0,-5,1) (0,-5,-2) (0,-5,2) (0,-5,-3) (0,-5,3) (0,-5,-4) (0,-5,4) (0,-5,-5) (0,-5,5) (0,-5,-6) (0,-5,6) (0,-5,-7) (0,-5,7) (0,-5,-8) (0,-5,8) (0,14,0) (0,6,-1) (0,6,1) (0,6,-2) (0,6,2) (0,6,-3) (0,6,3) (0,6,-4) (0,6,4) (0,6,-5) (0,6,5) (0,6,-6) (0,6,6) (0,6,-7) (0,6,7) (0,6,-8) (0,6,8) (0,-14,0) (0,-6,-1) (0,-6,1) (0,-6,-2) (0,-6,2) (0,-6,-3) (0,-6,3) (0,-6,-4) (0,-6,4) (0,-6,-5) (0,-6,5) (0,-6,-6) (0,-6,6) (0,-6,-7) (0,-6,7) (0,-6,-8) (0,-6,8) (0,15,0) (0,7,-1) (0,7,1) (0,7,-2) (0,7,2) (0,7,-3) (0,7,3) (0,7,-4) (0,7,4) (0,7,-5) (0,7,5) (0,7,-6) (0,7,6) (0,7,-7) (0,7,7) (0,7,-8) (0,7,8) (0,-15,0) (0,-7,-1) (0,-7,1) (0,-7,-2) (0,-7,2) (0,-7,-3) (0,-7,3) (0,-7,-4) (0,-7,4) (0,-7,-5) (0,-7,5) (0,-7,-6) (0,-7,6) (0,-7,-7) (0,-7,7) (0,-7,-8) (0,-7,8) (0,16,0) (0,8,-1) (0,8,1) (0,8,-2) (0,8,2) (0,8,-3) (0,8,3) (0,8,-4) (0,8,4) (0,8,-5) (0,8,5) (0,8,-6) (0,8,6) (0,8,-7) (0,8,7) (0,8,-8) (0,8,8) (0,-16,0) (0,-8,-1) (0,-8,1) (0,-8,-2) (0,-8,2) (0,-8,-3) (0,-8,3) (0,-8,-4) (0,-8,4) (0,-8,-5) (0,-8,5) (0,-8,-6) (0,-8,6) (0,-8,-7) (0,-8,7) (0,-8,-8) (0,-8,8) (0,0,9) (0,0,-9) (0,0,10) (0,0,-10) (0,0,11) (0,0,-11) (0,0,12) (0,0,-12) (0,0,13) (0,0,-13) (0,0,14) (0,0,-14) (0,0,15) (0,0,-15) (0,0,16) (0,0,-16)};
\end{axis}
\end{tikzpicture}
}
\vspace*{-1em}
\caption{Three-dimensional frequency set $I_8^3$ of superposition dimension two and the corresponding axis cross $\tilde{I}_8^3$, whose difference set $D(\tilde{I}_8^3)$ is a superset of the original set $I_8^3$, in comparison.}\label{fig:2danova}
\end{figure}
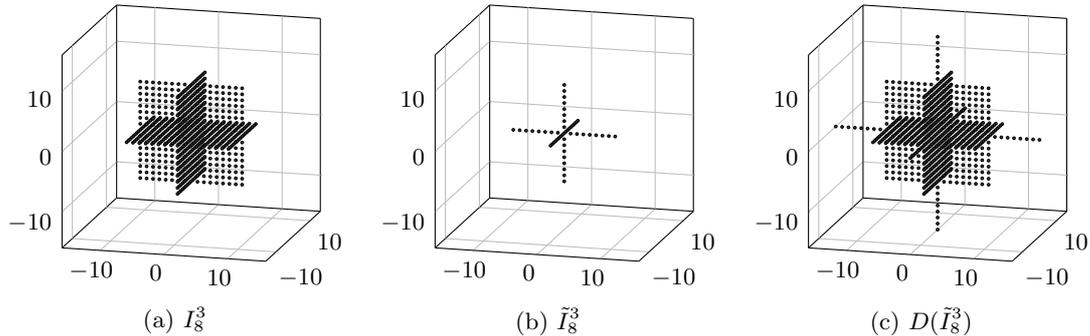

According to our assumptions, we expect small integration errors of $f$ for rank-1 lattice rules that integrates exactly trigonometric polynomials with frequency support at $\A_d$ within the cube $[-N,N]^d$.
Thus, we are interested in exactly integrating rank-1 lattices for $I_N^d:=\A_d\cap[-N,N]^d$, cf. Figure~\ref{fig:ANOVA_I} for illustration.
We determine the number of distinct elements within $I_N^d$
$$
|I_N^d|:=4\nchoosek{d}{2} N^2+2\nchoosek{d}{1} N + 1=2Nd(1+(d-1)N)+1.
$$
On the one hand, the presented approach allows for the direct computation of exactly integrating rank-1 lattices for the frequency set $I_N^d$.
On the other hand, we can take Remark~\ref{rem:diff_insetad_of_direct} into account and observe for the $d$-dimensional axis cross of expansion $2N$
$$
\tilde{I}_N^d:=\{\boldk\in\Z^d\colon \|\boldk\|_1=\|\boldk\|_\infty\le N\}
$$
that the corresponding difference set $D(\tilde{I}_N^d)$ is a superset of $I_N^d$, cf.
Figures~\ref{fig:axis_cross} and~\ref{fig:differenceset_axis_cross} for illustration.

Moreover, we observe
$$
|\tilde{I}_N^d|:=2dN+1\qquad\text{and}\qquad |D(\tilde{I}_N^d)|=2Nd(2+(d-1)N)+1.
$$
Accordingly, an exactly integrating rank-1 lattice for $D(\tilde{I}_N^d)$ needs to fulfill
the same requirements as an
exactly integrating rank-1 lattice for $I_N^d$ and, in addition, a few requirements more.
From a theoretical point of view, searching for an exactly integrating rank-1 lattice $D(\tilde{I}_N^d)$ should
result in at most slightly larger lattice sizes than the direct search for $I_N^d$.
From a practical point of view, the search for an exactly integrating rank-1 lattice for
$D(\tilde{I}_N^d)$, i.e., searching for a reconstructing rank-1 lattice for $\tilde{I}_N^d$, will lead to 
substantially reduced computational costs as mentioned in Remark~\ref{rem:diff_insetad_of_direct}.

\begin{figure}[tb]
\centering
\begin{tikzpicture}[baseline=(current axis.south)]%
\begin{semilogyaxis}[%
    scale only axis,
    xmin=1,xmax=51,
    ymin=10^2, ymax=6*10^5,
    title={\rule{6cm}{0pt}\bfseries Lattice sizes},
    x tick label style={align=center,font=\scriptsize},
    y tick label style={font=\scriptsize},
    width=0.85\textwidth,
    height=0.35\textwidth,
    xlabel={dimension $d$},
    legend style={font=\scriptsize, at={(0.03,1.0)},anchor=west, align=left},
    ]
\addplot[blue,mark=diamond,mark size=2pt,mark options={solid},error bars/.cd, y dir=both, y explicit, error bar style={solid}, error mark options={rotate=90,mark size=2,thick}] coordinates {
(2,8.3290e+03) -= (0,0.0000e+00) += (0,0.0000e+00) %
(3,6.1970e+03) -= (0,0.0000e+00) += (0,0.0000e+00) %
(4,6.8146e+03) -= (0,6.1760e+02) += (0,5.5584e+03) %
(5,1.0301e+04) -= (0,0.0000e+00) += (0,0.0000e+00) %
(6,1.1572e+04) -= (0,3.8550e+03) += (0,3.8550e+03) %
(7,1.0789e+04) -= (0,0.0000e+00) += (0,0.0000e+00) %
(8,1.4387e+04) -= (0,0.0000e+00) += (0,0.0000e+00) %
(9,1.8481e+04) -= (0,0.0000e+00) += (0,0.0000e+00) %
(10,2.1933e+04) -= (0,1.0384e+04) += (0,1.1538e+03) %
(11,1.5551e+04) -= (0,1.4076e+03) += (0,1.2668e+04) %
(12,1.6927e+04) -= (0,0.0000e+00) += (0,0.0000e+00) %
(13,2.0011e+04) -= (0,0.0000e+00) += (0,0.0000e+00) %
(14,2.3333e+04) -= (0,0.0000e+00) += (0,0.0000e+00) %
(15,2.6927e+04) -= (0,0.0000e+00) += (0,0.0000e+00) %
(16,3.0763e+04) -= (0,0.0000e+00) += (0,0.0000e+00) %
(17,3.3128e+04) -= (0,1.5685e+04) += (0,1.7428e+03) %
(18,3.1374e+04) -= (0,1.1765e+04) += (0,7.8432e+03) %
(19,3.5083e+04) -= (0,1.3154e+04) += (0,8.7696e+03) %
(20,2.6807e+04) -= (0,2.4360e+03) += (0,2.1924e+04) %
(21,2.6921e+04) -= (0,0.0000e+00) += (0,0.0000e+00) %
(22,2.9599e+04) -= (0,0.0000e+00) += (0,0.0000e+00) %
(23,3.2429e+04) -= (0,0.0000e+00) += (0,0.0000e+00) %
(24,3.5381e+04) -= (0,0.0000e+00) += (0,0.0000e+00) %
(25,3.8447e+04) -= (0,0.0000e+00) += (0,0.0000e+00) %
(26,4.1641e+04) -= (0,0.0000e+00) += (0,0.0000e+00) %
(27,4.4959e+04) -= (0,0.0000e+00) += (0,0.0000e+00) %
(28,4.8437e+04) -= (0,0.0000e+00) += (0,0.0000e+00) %
(29,5.2009e+04) -= (0,0.0000e+00) += (0,0.0000e+00) %
(30,5.5717e+04) -= (0,0.0000e+00) += (0,0.0000e+00) %
(31,5.9557e+04) -= (0,0.0000e+00) += (0,0.0000e+00) %
(32,6.3541e+04) -= (0,0.0000e+00) += (0,0.0000e+00) %
(33,6.4251e+04) -= (0,3.0424e+04) += (0,3.3804e+03) %
(34,7.1861e+04) -= (0,0.0000e+00) += (0,0.0000e+00) %
(35,7.2398e+04) -= (0,3.4285e+04) += (0,3.8094e+03) %
(36,7.2631e+04) -= (0,3.2280e+04) += (0,8.0700e+03) %
(37,7.2512e+04) -= (0,2.9845e+04) += (0,1.2791e+04) %
(38,8.5553e+04) -= (0,4.0500e+04) += (0,4.5000e+03) %
(39,8.0670e+04) -= (0,3.3211e+04) += (0,1.4233e+04) %
(40,7.9923e+04) -= (0,2.9966e+04) += (0,1.9978e+04) %
(41,7.8765e+04) -= (0,2.6254e+04) += (0,2.6254e+04) %
(42,6.0658e+04) -= (0,5.5114e+03) += (0,4.9603e+04) %
(43,5.7829e+04) -= (0,0.0000e+00) += (0,0.0000e+00) %
(44,6.0589e+04) -= (0,0.0000e+00) += (0,0.0000e+00) %
(45,6.3391e+04) -= (0,0.0000e+00) += (0,0.0000e+00) %
(46,7.2898e+04) -= (0,6.6270e+03) += (0,5.9643e+04) %
(47,6.9221e+04) -= (0,0.0000e+00) += (0,0.0000e+00) %
(48,7.2221e+04) -= (0,0.0000e+00) += (0,0.0000e+00) %
(49,7.5307e+04) -= (0,0.0000e+00) += (0,0.0000e+00) %
(50,7.8467e+04) -= (0,0.0000e+00) += (0,0.0000e+00) %
};
\addlegendentry{found lattice sizes for $I_{64}^d$}
\addplot[red,mark=diamond,mark size=2pt,mark options={solid},error bars/.cd, y dir=both, y explicit, error bar style={solid}, error mark options={rotate=90,mark size=2,thick}] coordinates {
(2,8.2630e+03) -= (0,0.0000e+00) += (0,0.0000e+00) %
(3,6.9580e+03) -= (0,2.3190e+03) += (0,2.3190e+03) %
(4,8.2430e+03) -= (0,0.0000e+00) += (0,0.0000e+00) %
(5,1.2853e+04) -= (0,0.0000e+00) += (0,0.0000e+00) %
(6,9.2570e+03) -= (0,0.0000e+00) += (0,0.0000e+00) %
(7,1.2583e+04) -= (0,0.0000e+00) += (0,0.0000e+00) %
(8,1.6447e+04) -= (0,0.0000e+00) += (0,0.0000e+00) %
(9,1.8731e+04) -= (0,8.3040e+03) += (0,2.0760e+03) %
(10,1.2829e+04) -= (0,0.0000e+00) += (0,0.0000e+00) %
(11,1.5527e+04) -= (0,0.0000e+00) += (0,0.0000e+00) %
(12,1.8481e+04) -= (0,0.0000e+00) += (0,0.0000e+00) %
(13,2.1701e+04) -= (0,0.0000e+00) += (0,0.0000e+00) %
(14,2.5147e+04) -= (0,0.0000e+00) += (0,0.0000e+00) %
(15,2.8843e+04) -= (0,0.0000e+00) += (0,0.0000e+00) %
(16,2.9548e+04) -= (0,1.3131e+04) += (0,3.2828e+03) %
(17,3.5187e+04) -= (0,1.6666e+04) += (0,1.8518e+03) %
(18,3.7362e+04) -= (0,1.6603e+04) += (0,4.1508e+03) %
(19,2.5444e+04) -= (0,2.3130e+03) += (0,2.0817e+04) %
(20,3.5883e+04) -= (0,1.0250e+04) += (0,1.5374e+04) %
(21,3.1099e+04) -= (0,2.8224e+03) += (0,2.5402e+04) %
(22,3.1013e+04) -= (0,0.0000e+00) += (0,0.0000e+00) %
(23,3.3889e+04) -= (0,0.0000e+00) += (0,0.0000e+00) %
(24,3.6913e+04) -= (0,0.0000e+00) += (0,0.0000e+00) %
(25,4.0037e+04) -= (0,0.0000e+00) += (0,0.0000e+00) %
(26,4.3313e+04) -= (0,0.0000e+00) += (0,0.0000e+00) %
(27,4.6703e+04) -= (0,0.0000e+00) += (0,0.0000e+00) %
(28,5.0227e+04) -= (0,0.0000e+00) += (0,0.0000e+00) %
(29,5.3881e+04) -= (0,0.0000e+00) += (0,0.0000e+00) %
(30,5.7641e+04) -= (0,0.0000e+00) += (0,0.0000e+00) %
(31,6.1543e+04) -= (0,0.0000e+00) += (0,0.0000e+00) %
(32,6.5587e+04) -= (0,0.0000e+00) += (0,0.0000e+00) %
(33,6.9761e+04) -= (0,0.0000e+00) += (0,0.0000e+00) %
(34,7.4047e+04) -= (0,0.0000e+00) += (0,0.0000e+00) %
(35,7.8467e+04) -= (0,0.0000e+00) += (0,0.0000e+00) %
(36,7.4699e+04) -= (0,3.3192e+04) += (0,8.2980e+03) %
(37,7.8907e+04) -= (0,3.5054e+04) += (0,8.7636e+03) %
(38,6.9370e+04) -= (0,2.3109e+04) += (0,2.3109e+04) %
(39,6.3339e+04) -= (0,1.4608e+04) += (0,3.4084e+04) %
(40,6.6606e+04) -= (0,1.5367e+04) += (0,3.5855e+04) %
(41,5.9212e+04) -= (0,5.3810e+03) += (0,4.8429e+04) %
(42,5.6477e+04) -= (0,0.0000e+00) += (0,0.0000e+00) %
(43,6.5126e+04) -= (0,5.9192e+03) += (0,5.3273e+04) %
(44,6.1987e+04) -= (0,0.0000e+00) += (0,0.0000e+00) %
(45,6.4849e+04) -= (0,0.0000e+00) += (0,0.0000e+00) %
(46,6.7741e+04) -= (0,0.0000e+00) += (0,0.0000e+00) %
(47,7.0729e+04) -= (0,0.0000e+00) += (0,0.0000e+00) %
(48,7.3771e+04) -= (0,0.0000e+00) += (0,0.0000e+00) %
(49,7.6871e+04) -= (0,0.0000e+00) += (0,0.0000e+00) %
(50,8.0039e+04) -= (0,0.0000e+00) += (0,0.0000e+00) %
};
\addlegendentry{found lattice sizes for $D(\tilde{I}_{64}^d)$}
\addplot[black,thick,domain=2:50, dashed, no markers,samples=200] {8192*(x-1)+2}; 
\addlegendentry{upper bound $8192\,(d-1)+2$ on $M_{\operatorname{lb}}$}
\addplot[black,thick,domain=2:50, dotted, no markers,samples=200] {64*x+1}; 
\addlegendentry{rough absolute lower bound $|\tilde{I}_{\floor{64/2}}^d|= 64\,d+1$}
\end{semilogyaxis}
\end{tikzpicture}
\caption{Comparison of results of Algorithm~\ref{alg:heuristic}. As expected, the lattice sizes of exactly integrating rank-1 lattices constructed for $I_{64}^d$ differ only slightly from the lattice sizes of reconstructing rank-1 lattices constructed for $\tilde{I}_{64}^d$.}
\label{fig:anova_M}
\end{figure}
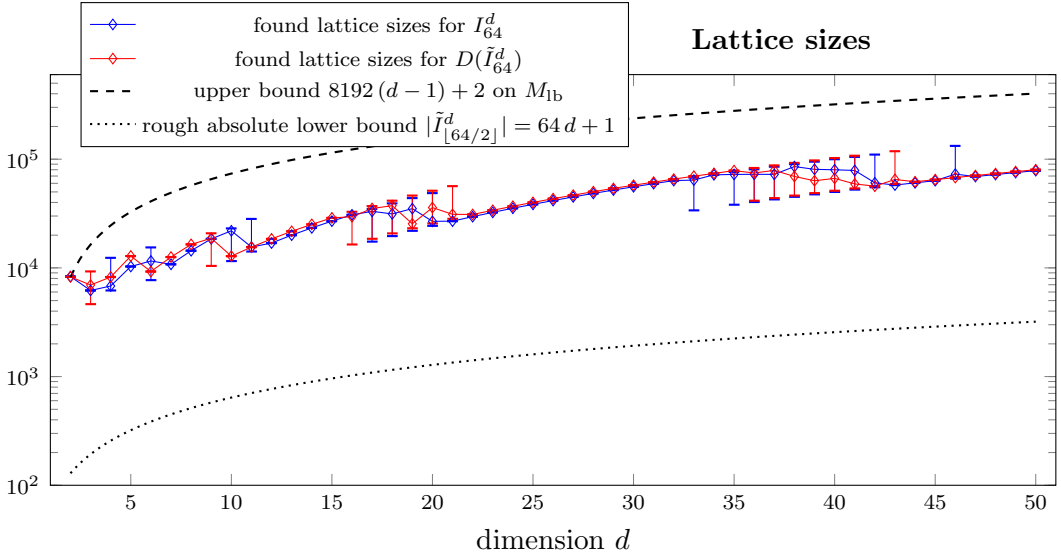

\begin{figure}[tb]
\centering
\begin{tikzpicture}[baseline=(current axis.south)]%
\begin{loglogaxis}[%
    scale only axis,
    xmin=1,xmax=2200,
    ymin=10^2, ymax=2*10^7,
    title={\rule{6cm}{0pt}\bfseries Lattice sizes},
    x tick label style={align=center,font=\scriptsize},
    y tick label style={font=\scriptsize},
    width=0.85\textwidth,
    height=0.35\textwidth,
    xlabel={dimension $d$},
    legend style={font=\scriptsize, at={(0.03,0.96)},anchor=west, align=left},
    ]
\addplot[blue,mark=diamond,mark size=2pt,mark options={solid},error bars/.cd, y dir=both, y explicit, error bar style={solid}, error mark options={rotate=90,mark size=2,thick}] coordinates {
(2,8.3290e+03) -= (0,0.0000e+00) += (0,0.0000e+00) %
(3,6.1970e+03) -= (0,0.0000e+00) += (0,0.0000e+00) %
(4,6.8146e+03) -= (0,6.1760e+02) += (0,5.5584e+03) %
(5,1.0301e+04) -= (0,0.0000e+00) += (0,0.0000e+00) %
(6,1.1572e+04) -= (0,3.8550e+03) += (0,3.8550e+03) %
(7,1.0789e+04) -= (0,0.0000e+00) += (0,0.0000e+00) %
(8,1.4387e+04) -= (0,0.0000e+00) += (0,0.0000e+00) %
(9,1.8481e+04) -= (0,0.0000e+00) += (0,0.0000e+00) %
(10,2.1933e+04) -= (0,1.0384e+04) += (0,1.1538e+03) %
(11,1.5551e+04) -= (0,1.4076e+03) += (0,1.2668e+04) %
(12,1.6927e+04) -= (0,0.0000e+00) += (0,0.0000e+00) %
(13,2.0011e+04) -= (0,0.0000e+00) += (0,0.0000e+00) %
(14,2.3333e+04) -= (0,0.0000e+00) += (0,0.0000e+00) %
(15,2.6927e+04) -= (0,0.0000e+00) += (0,0.0000e+00) %
(16,3.0763e+04) -= (0,0.0000e+00) += (0,0.0000e+00) %
(17,3.3128e+04) -= (0,1.5685e+04) += (0,1.7428e+03) %
(18,3.1374e+04) -= (0,1.1765e+04) += (0,7.8432e+03) %
(19,3.5083e+04) -= (0,1.3154e+04) += (0,8.7696e+03) %
(20,2.6807e+04) -= (0,2.4360e+03) += (0,2.1924e+04) %
(21,2.6921e+04) -= (0,0.0000e+00) += (0,0.0000e+00) %
(22,2.9599e+04) -= (0,0.0000e+00) += (0,0.0000e+00) %
(23,3.2429e+04) -= (0,0.0000e+00) += (0,0.0000e+00) %
(24,3.5381e+04) -= (0,0.0000e+00) += (0,0.0000e+00) %
(25,3.8447e+04) -= (0,0.0000e+00) += (0,0.0000e+00) %
(26,4.1641e+04) -= (0,0.0000e+00) += (0,0.0000e+00) %
(27,4.4959e+04) -= (0,0.0000e+00) += (0,0.0000e+00) %
(28,4.8437e+04) -= (0,0.0000e+00) += (0,0.0000e+00) %
(29,5.2009e+04) -= (0,0.0000e+00) += (0,0.0000e+00) %
(30,5.5717e+04) -= (0,0.0000e+00) += (0,0.0000e+00) %
(31,5.9557e+04) -= (0,0.0000e+00) += (0,0.0000e+00) %
(32,6.3541e+04) -= (0,0.0000e+00) += (0,0.0000e+00) %
(33,6.4251e+04) -= (0,3.0424e+04) += (0,3.3804e+03) %
(34,7.1861e+04) -= (0,0.0000e+00) += (0,0.0000e+00) %
(35,7.2398e+04) -= (0,3.4285e+04) += (0,3.8094e+03) %
(36,7.2631e+04) -= (0,3.2280e+04) += (0,8.0700e+03) %
(37,7.2512e+04) -= (0,2.9845e+04) += (0,1.2791e+04) %
(38,8.5553e+04) -= (0,4.0500e+04) += (0,4.5000e+03) %
(39,8.0670e+04) -= (0,3.3211e+04) += (0,1.4233e+04) %
(40,7.9923e+04) -= (0,2.9966e+04) += (0,1.9978e+04) %
(41,7.8765e+04) -= (0,2.6254e+04) += (0,2.6254e+04) %
(42,6.0658e+04) -= (0,5.5114e+03) += (0,4.9603e+04) %
(43,5.7829e+04) -= (0,0.0000e+00) += (0,0.0000e+00) %
(44,6.0589e+04) -= (0,0.0000e+00) += (0,0.0000e+00) %
(45,6.3391e+04) -= (0,0.0000e+00) += (0,0.0000e+00) %
(46,7.2898e+04) -= (0,6.6270e+03) += (0,5.9643e+04) %
(47,6.9221e+04) -= (0,0.0000e+00) += (0,0.0000e+00) %
(48,7.2221e+04) -= (0,0.0000e+00) += (0,0.0000e+00) %
(49,7.5307e+04) -= (0,0.0000e+00) += (0,0.0000e+00) %
(50,7.8467e+04) -= (0,0.0000e+00) += (0,0.0000e+00) %
(60,1.1333e+05) -= (0,0.0000e+00) += (0,0.0000e+00) %
(70,1.5461e+05) -= (0,0.0000e+00) += (0,0.0000e+00) %
(80,2.0231e+05) -= (0,0.0000e+00) += (0,0.0000e+00) %
(90,1.5384e+05) -= (0,2.5638e+04) += (0,1.0255e+05) %
(100,1.5844e+05) -= (0,0.0000e+00) += (0,0.0000e+00) %
};
\addlegendentry{found lattice sizes for $I_{64}^d$}
\addplot[red,mark=diamond,mark size=2pt,mark options={solid},error bars/.cd, y dir=both, y explicit, error bar style={solid}, error mark options={rotate=90,mark size=2,thick}] coordinates {
(2,8.2630e+03) -= (0,0.0000e+00) += (0,0.0000e+00) %
(3,6.9580e+03) -= (0,2.3190e+03) += (0,2.3190e+03) %
(4,8.2430e+03) -= (0,0.0000e+00) += (0,0.0000e+00) %
(5,1.2853e+04) -= (0,0.0000e+00) += (0,0.0000e+00) %
(6,9.2570e+03) -= (0,0.0000e+00) += (0,0.0000e+00) %
(7,1.2583e+04) -= (0,0.0000e+00) += (0,0.0000e+00) %
(8,1.6447e+04) -= (0,0.0000e+00) += (0,0.0000e+00) %
(9,1.8731e+04) -= (0,8.3040e+03) += (0,2.0760e+03) %
(10,1.2829e+04) -= (0,0.0000e+00) += (0,0.0000e+00) %
(11,1.5527e+04) -= (0,0.0000e+00) += (0,0.0000e+00) %
(12,1.8481e+04) -= (0,0.0000e+00) += (0,0.0000e+00) %
(13,2.1701e+04) -= (0,0.0000e+00) += (0,0.0000e+00) %
(14,2.5147e+04) -= (0,0.0000e+00) += (0,0.0000e+00) %
(15,2.8843e+04) -= (0,0.0000e+00) += (0,0.0000e+00) %
(16,2.9548e+04) -= (0,1.3131e+04) += (0,3.2828e+03) %
(17,3.5187e+04) -= (0,1.6666e+04) += (0,1.8518e+03) %
(18,3.7362e+04) -= (0,1.6603e+04) += (0,4.1508e+03) %
(19,2.5444e+04) -= (0,2.3130e+03) += (0,2.0817e+04) %
(20,3.5883e+04) -= (0,1.0250e+04) += (0,1.5374e+04) %
(21,3.1099e+04) -= (0,2.8224e+03) += (0,2.5402e+04) %
(22,3.1013e+04) -= (0,0.0000e+00) += (0,0.0000e+00) %
(23,3.3889e+04) -= (0,0.0000e+00) += (0,0.0000e+00) %
(24,3.6913e+04) -= (0,0.0000e+00) += (0,0.0000e+00) %
(25,4.0037e+04) -= (0,0.0000e+00) += (0,0.0000e+00) %
(26,4.3313e+04) -= (0,0.0000e+00) += (0,0.0000e+00) %
(27,4.6703e+04) -= (0,0.0000e+00) += (0,0.0000e+00) %
(28,5.0227e+04) -= (0,0.0000e+00) += (0,0.0000e+00) %
(29,5.3881e+04) -= (0,0.0000e+00) += (0,0.0000e+00) %
(30,5.7641e+04) -= (0,0.0000e+00) += (0,0.0000e+00) %
(31,6.1543e+04) -= (0,0.0000e+00) += (0,0.0000e+00) %
(32,6.5587e+04) -= (0,0.0000e+00) += (0,0.0000e+00) %
(33,6.9761e+04) -= (0,0.0000e+00) += (0,0.0000e+00) %
(34,7.4047e+04) -= (0,0.0000e+00) += (0,0.0000e+00) %
(35,7.8467e+04) -= (0,0.0000e+00) += (0,0.0000e+00) %
(36,7.4699e+04) -= (0,3.3192e+04) += (0,8.2980e+03) %
(37,7.8907e+04) -= (0,3.5054e+04) += (0,8.7636e+03) %
(38,6.9370e+04) -= (0,2.3109e+04) += (0,2.3109e+04) %
(39,6.3339e+04) -= (0,1.4608e+04) += (0,3.4084e+04) %
(40,6.6606e+04) -= (0,1.5367e+04) += (0,3.5855e+04) %
(41,5.9212e+04) -= (0,5.3810e+03) += (0,4.8429e+04) %
(42,5.6477e+04) -= (0,0.0000e+00) += (0,0.0000e+00) %
(43,6.5126e+04) -= (0,5.9192e+03) += (0,5.3273e+04) %
(44,6.1987e+04) -= (0,0.0000e+00) += (0,0.0000e+00) %
(45,6.4849e+04) -= (0,0.0000e+00) += (0,0.0000e+00) %
(46,6.7741e+04) -= (0,0.0000e+00) += (0,0.0000e+00) %
(47,7.0729e+04) -= (0,0.0000e+00) += (0,0.0000e+00) %
(48,7.3771e+04) -= (0,0.0000e+00) += (0,0.0000e+00) %
(49,7.6871e+04) -= (0,0.0000e+00) += (0,0.0000e+00) %
(50,8.0039e+04) -= (0,0.0000e+00) += (0,0.0000e+00) %
(60,1.1524e+05) -= (0,0.0000e+00) += (0,0.0000e+00) %
(70,1.5689e+05) -= (0,0.0000e+00) += (0,0.0000e+00) %
(80,1.8437e+05) -= (0,8.1939e+04) += (0,2.0485e+04) %
(90,1.2964e+05) -= (0,0.0000e+00) += (0,0.0000e+00) %
(100,1.6003e+05) -= (0,0.0000e+00) += (0,0.0000e+00) %
(150,3.6007e+05) -= (0,0.0000e+00) += (0,0.0000e+00) %
(200,3.2004e+05) -= (0,0.0000e+00) += (0,0.0000e+00) %
(250,5.0004e+05) -= (0,0.0000e+00) += (0,0.0000e+00) %
(300,7.2006e+05) -= (0,0.0000e+00) += (0,0.0000e+00) %
(350,9.8007e+05) -= (0,0.0000e+00) += (0,0.0000e+00) %
(400,8.9608e+05) -= (0,2.5602e+05) += (0,3.8402e+05) %
(450,8.1005e+05) -= (0,0.0000e+00) += (0,0.0000e+00) %
(500,1.0001e+06) -= (0,0.0000e+00) += (0,0.0000e+00) %
(600,1.4401e+06) -= (0,0.0000e+00) += (0,0.0000e+00) %
(700,1.9601e+06) -= (0,0.0000e+00) += (0,0.0000e+00) %
(800,2.1761e+06) -= (0,8.9599e+05) += (0,3.8400e+05) %
(900,1.6200e+06) -= (0,0.0000e+00) += (0,0.0000e+00) %
(1000,2.0001e+06) -= (0,0.0000e+00) += (0,0.0000e+00) %
(1500,4.5001e+06) -= (0,0.0000e+00) += (0,0.0000e+00) %
(2000,4.0001e+06) -= (0,0.0000e+00) += (0,0.0000e+00) %
};
\addlegendentry{found lattice sizes for $D(\tilde{I}_{64}^d)$}
\addplot[black,thick,domain=2:2000, dashed, no markers,samples=200] {8192*(x-1)+2}; 
\addlegendentry{upper bound $8192\,(d-1)+2$ on $M_{\operatorname{lb}}$}
\addplot[black,thick,domain=2:2000, dotted, no markers,samples=200] {64*x+1}; 
\addlegendentry{rough absolute lower bound $|\tilde{I}_{\floor{64/2}}^d|= 64\,d+1$}
\addplot[red,mark=diamond*,mark size=2pt,mark options={solid},error bars/.cd, y dir=both, y explicit, error bar style={solid}, error mark options={rotate=90,mark size=2,thick}] coordinates{
(2,4.1315e+03) -= (0,0.0000e+00) += (0,0.0000e+00) %
(3,2.3193e+03) -= (0,7.7300e+02) += (0,7.7300e+02) %
(4,2.0608e+03) -= (0,0.0000e+00) += (0,0.0000e+00) %
(5,2.5706e+03) -= (0,0.0000e+00) += (0,0.0000e+00) %
(6,1.5428e+03) -= (0,0.0000e+00) += (0,0.0000e+00) %
(7,1.7976e+03) -= (0,0.0000e+00) += (0,0.0000e+00) %
(8,2.0559e+03) -= (0,0.0000e+00) += (0,0.0000e+00) %
(9,2.0812e+03) -= (0,9.2267e+02) += (0,2.3067e+02) %
(10,1.2829e+03) -= (0,0.0000e+00) += (0,0.0000e+00) %
(11,1.4115e+03) -= (0,0.0000e+00) += (0,0.0000e+00) %
(12,1.5401e+03) -= (0,0.0000e+00) += (0,0.0000e+00) %
(13,1.6693e+03) -= (0,0.0000e+00) += (0,0.0000e+00) %
(14,1.7962e+03) -= (0,0.0000e+00) += (0,0.0000e+00) %
(15,1.9229e+03) -= (0,0.0000e+00) += (0,0.0000e+00) %
(16,1.8468e+03) -= (0,8.2070e+02) += (0,2.0517e+02) %
(17,2.0698e+03) -= (0,9.8036e+02) += (0,1.0893e+02) %
(18,2.0757e+03) -= (0,9.2240e+02) += (0,2.3060e+02) %
(19,1.3392e+03) -= (0,1.2174e+02) += (0,1.0956e+03) %
(20,1.7941e+03) -= (0,5.1248e+02) += (0,7.6872e+02) %
(21,1.4809e+03) -= (0,1.3440e+02) += (0,1.2096e+03) %
(22,1.4097e+03) -= (0,0.0000e+00) += (0,0.0000e+00) %
(23,1.4734e+03) -= (0,0.0000e+00) += (0,0.0000e+00) %
(24,1.5380e+03) -= (0,0.0000e+00) += (0,0.0000e+00) %
(25,1.6015e+03) -= (0,0.0000e+00) += (0,0.0000e+00) %
(26,1.6659e+03) -= (0,0.0000e+00) += (0,0.0000e+00) %
(27,1.7297e+03) -= (0,0.0000e+00) += (0,0.0000e+00) %
(28,1.7938e+03) -= (0,0.0000e+00) += (0,0.0000e+00) %
(29,1.8580e+03) -= (0,0.0000e+00) += (0,0.0000e+00) %
(30,1.9214e+03) -= (0,0.0000e+00) += (0,0.0000e+00) %
(31,1.9853e+03) -= (0,0.0000e+00) += (0,0.0000e+00) %
(32,2.0496e+03) -= (0,0.0000e+00) += (0,0.0000e+00) %
(33,2.1140e+03) -= (0,0.0000e+00) += (0,0.0000e+00) %
(34,2.1779e+03) -= (0,0.0000e+00) += (0,0.0000e+00) %
(35,2.2419e+03) -= (0,0.0000e+00) += (0,0.0000e+00) %
(36,2.0750e+03) -= (0,9.2200e+02) += (0,2.3050e+02) %
(37,2.1326e+03) -= (0,9.4742e+02) += (0,2.3685e+02) %
(38,1.8255e+03) -= (0,6.0813e+02) += (0,6.0813e+02) %
(39,1.6241e+03) -= (0,3.7455e+02) += (0,8.7396e+02) %
(40,1.6651e+03) -= (0,3.8417e+02) += (0,8.9638e+02) %
(41,1.4442e+03) -= (0,1.3124e+02) += (0,1.1812e+03) %
(42,1.3447e+03) -= (0,0.0000e+00) += (0,0.0000e+00) %
(43,1.5146e+03) -= (0,1.3766e+02) += (0,1.2389e+03) %
(44,1.4088e+03) -= (0,0.0000e+00) += (0,0.0000e+00) %
(45,1.4411e+03) -= (0,0.0000e+00) += (0,0.0000e+00) %
(46,1.4726e+03) -= (0,0.0000e+00) += (0,0.0000e+00) %
(47,1.5049e+03) -= (0,0.0000e+00) += (0,0.0000e+00) %
(48,1.5369e+03) -= (0,0.0000e+00) += (0,0.0000e+00) %
(49,1.5688e+03) -= (0,0.0000e+00) += (0,0.0000e+00) %
(50,1.6008e+03) -= (0,0.0000e+00) += (0,0.0000e+00) %
(60,1.9206e+03) -= (0,0.0000e+00) += (0,0.0000e+00) %
(70,2.2412e+03) -= (0,0.0000e+00) += (0,0.0000e+00) %
(80,2.3047e+03) -= (0,1.0242e+03) += (0,2.5606e+02) %
(90,1.4405e+03) -= (0,0.0000e+00) += (0,0.0000e+00) %
(100,1.6003e+03) -= (0,0.0000e+00) += (0,0.0000e+00) %
(150,2.4005e+03) -= (0,0.0000e+00) += (0,0.0000e+00) %
(200,1.6002e+03) -= (0,0.0000e+00) += (0,0.0000e+00) %
(250,2.0002e+03) -= (0,0.0000e+00) += (0,0.0000e+00) %
(300,2.4002e+03) -= (0,0.0000e+00) += (0,0.0000e+00) %
(350,2.8002e+03) -= (0,0.0000e+00) += (0,0.0000e+00) %
(400,2.2402e+03) -= (0,6.4004e+02) += (0,9.6006e+02) %
(450,1.8001e+03) -= (0,0.0000e+00) += (0,0.0000e+00) %
(500,2.0002e+03) -= (0,0.0000e+00) += (0,0.0000e+00) %
(600,2.4001e+03) -= (0,0.0000e+00) += (0,0.0000e+00) %
(700,2.8001e+03) -= (0,0.0000e+00) += (0,0.0000e+00) %
(800,2.7201e+03) -= (0,1.1200e+03) += (0,4.7999e+02) %
(900,1.8000e+03) -= (0,0.0000e+00) += (0,0.0000e+00) %
(1000,2.0001e+03) -= (0,0.0000e+00) += (0,0.0000e+00) %
(1500,3.0000e+03) -= (0,0.0000e+00) += (0,0.0000e+00) %
(2000,2.0000e+03) -= (0,0.0000e+00) += (0,0.0000e+00) %
};
\addlegendentry{$d^{-1}$ $\times$ found lattice sizes for $D(\tilde{I}_{64}^d)$}
\end{loglogaxis}
\end{tikzpicture}
\caption{Double-log depiction and continuation of Figure~\ref{fig:anova_M}. The lattice sizes are determined using Algorithm~\ref{alg:heuristic}.}
\label{fig:anova_M_cont}
\end{figure}
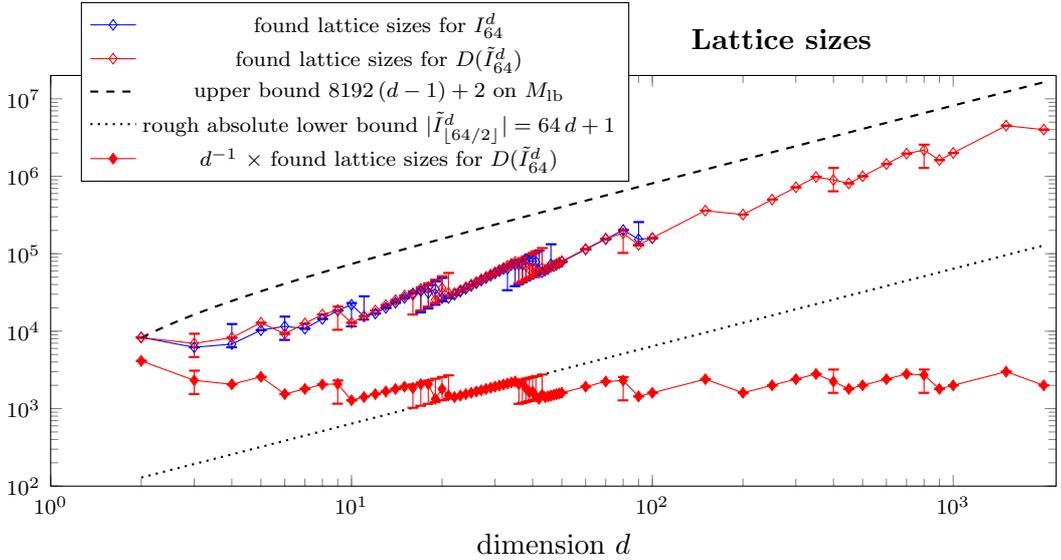

A further detailed analysis of the frequency sets $I_N^d$ and $\tilde{I}_N^d$ leads to the following two observations.
First, we have
$D(\tilde{I}_{\floor{N/2}}^d)\subset I_N^d$, which yields that an exactly integrating rank-1 lattice for $I_N^d$
requires at least $M \ge |\tilde{I}_{\floor{N/2}}^d|=2d\floor{N/2}+1$ due to
the fact that this rank-1 lattice is also a reconstructing rank-1 lattice for $\tilde{I}_{\floor{N/2}}^d$, cf.~\ref{eq:reco_property}. 
Second, taking \cite[Theorem 1]{Kae2013} into account, we observe $M_{\operatorname{lb}}\le 2(d-1)N^2+2$.

The numerically detected lattice sizes determined by Algorithm~\ref{alg:heuristic} are plotted in Figures~\ref{fig:anova_M}
and~\ref{fig:anova_M_cont} for $N=64$ and dimensions up to 50 and 2000, respectively.
The lattice sizes of the reconstructing rank-1 lattices for $\tilde{I}_{64}^d$ and the exactly integrating
rank-1 lattices for $I_{64}^d$ differ only slightly. Moreover, we observe that
the determined rank-1 lattice sizes are less than the known upper bounds on
$M_{\operatorname{lb}}$, cf.~\cite{Kae2013, KuoMiNoNu19}.
In Figure~\ref{fig:anova_M_cont}, we additionally plotted the lattice sizes scaled by $d^{-1}$, which
leads to a bounded graph, i.e., the lattice sizes seem to depend linearly on the dimension $d$
just as the above determined bounds on $M_{\operatorname{lb}}$ indicate.

The crucial difference between the two approaches of determining exactly integrating rank\mbox{-}1 lattices for $I_N^d$
and reconstructing rank-1 lattices for $\tilde{I}_N^d$ are the numbers of elements that one has to process.
Since $|\tilde{I}_N^d|\sim |I_N^d|^2$ holds, we observe substantially reduced computational costs for the second approach.
\begin{figure}[tb]
\centering
\begin{tikzpicture}[baseline=(current axis.south)]%
\begin{loglogaxis}[%
    scale only axis,
    xmin=1,xmax=2200,
    ymin=10^-3, ymax=6*10^3,
    title={\bfseries Runtimes},
    x tick label style={align=center,font=\scriptsize},
    y tick label style={font=\scriptsize},
    xlabel={dimension $d$},
    ylabel={runtime in seconds},
    width=0.85\textwidth,
    height=0.35\textwidth,
    legend pos=south east,
    legend style={font=\scriptsize},
    ]
\addplot[blue,mark=diamond,mark size=2pt,mark options={solid},error bars/.cd, y dir=both, y explicit, error bar style={solid}, error mark options={rotate=90,mark size=2,thick}] coordinates {
(2,1.5207e-01) -= (0,3.5804e-03) += (0,1.6021e-02) %
(3,3.4669e-01) -= (0,2.4583e-02) += (0,3.4476e-02) %
(4,7.3252e-01) -= (0,1.5944e-01) += (0,2.0654e-01) %
(5,1.0862e+00) -= (0,1.5259e-01) += (0,2.5499e-01) %
(6,1.7667e+00) -= (0,4.2625e-01) += (0,4.8334e-01) %
(7,2.1956e+00) -= (0,2.4804e-01) += (0,1.6368e-01) %
(8,3.0291e+00) -= (0,1.7940e-01) += (0,5.1329e-01) %
(9,4.4467e+00) -= (0,2.9436e-01) += (0,5.7906e-01) %
(10,5.5487e+00) -= (0,6.0922e-01) += (0,4.0769e-01) %
(11,6.8257e+00) -= (0,1.0162e+00) += (0,2.2466e+00) %
(12,7.9995e+00) -= (0,8.2896e-01) += (0,1.7508e+00) %
(13,9.8104e+00) -= (0,9.8010e-01) += (0,7.4846e-01) %
(14,1.1927e+01) -= (0,5.9108e-01) += (0,9.9198e-01) %
(15,1.4299e+01) -= (0,7.5891e-01) += (0,1.0383e+00) %
(16,1.7168e+01) -= (0,8.6944e-01) += (0,9.0778e-01) %
(17,2.1203e+01) -= (0,1.7329e+00) += (0,1.6841e+00) %
(18,2.4718e+01) -= (0,2.9985e+00) += (0,5.5711e+00) %
(19,2.7392e+01) -= (0,2.4770e+00) += (0,3.2926e+00) %
(20,3.4227e+01) -= (0,6.0844e+00) += (0,6.5624e+00) %
(21,3.7314e+01) -= (0,5.8941e+00) += (0,7.2638e+00) %
(22,4.0467e+01) -= (0,4.5308e+00) += (0,7.4821e+00) %
(23,4.5740e+01) -= (0,4.5720e+00) += (0,5.5935e+00) %
(24,5.0101e+01) -= (0,4.5377e+00) += (0,5.2706e+00) %
(25,5.4418e+01) -= (0,3.7996e+00) += (0,6.4059e+00) %
(26,5.9125e+01) -= (0,1.6101e+00) += (0,1.1848e+00) %
(27,6.5594e+01) -= (0,2.1910e+00) += (0,2.4932e+00) %
(28,7.4502e+01) -= (0,1.9288e+00) += (0,2.4286e+00) %
(29,8.6063e+01) -= (0,4.2513e+00) += (0,3.7696e+00) %
(30,9.3229e+01) -= (0,3.0454e+00) += (0,3.4229e+00) %
(31,1.0303e+02) -= (0,3.9460e+00) += (0,6.0717e+00) %
(32,1.1093e+02) -= (0,3.8086e+00) += (0,3.0533e+00) %
(33,1.2314e+02) -= (0,8.5367e+00) += (0,8.9493e+00) %
(34,1.3631e+02) -= (0,3.9772e+00) += (0,5.8676e+00) %
(35,1.4598e+02) -= (0,9.1629e+00) += (0,5.1885e+00) %
(36,1.6388e+02) -= (0,7.2036e+00) += (0,1.4188e+01) %
(37,1.7876e+02) -= (0,1.4005e+01) += (0,3.8142e+01) %
(38,1.8563e+02) -= (0,6.9138e+00) += (0,1.5061e+01) %
(39,2.0929e+02) -= (0,1.0617e+01) += (0,3.2768e+01) %
(40,2.2210e+02) -= (0,2.6285e+01) += (0,3.4812e+01) %
(41,2.3683e+02) -= (0,2.7900e+01) += (0,2.8290e+01) %
(42,2.4878e+02) -= (0,2.0857e+01) += (0,2.4881e+01) %
(43,3.0662e+02) -= (0,4.9245e+01) += (0,3.8946e+01) %
(44,2.9109e+02) -= (0,3.6831e+01) += (0,6.0529e+01) %
(45,3.0805e+02) -= (0,3.4764e+01) += (0,4.8703e+01) %
(46,3.1676e+02) -= (0,1.7037e+01) += (0,3.8355e+01) %
(47,3.2559e+02) -= (0,2.6297e+01) += (0,2.9357e+01) %
(48,3.4422e+02) -= (0,1.4564e+01) += (0,9.0924e+01) %
(49,3.6771e+02) -= (0,1.2283e+01) += (0,4.6532e+01) %
(50,3.9171e+02) -= (0,1.0493e+01) += (0,2.2231e+01) %
(60,6.8069e+02) -= (0,1.7945e+01) += (0,1.3888e+01) %
(70,1.0898e+03) -= (0,3.6105e+01) += (0,4.9150e+01) %
(80,1.6371e+03) -= (0,5.9412e+01) += (0,5.0607e+01) %
(90,2.3798e+03) -= (0,2.4067e+02) += (0,4.4164e+02) %
(100,3.0406e+03) -= (0,1.4379e+02) += (0,2.9787e+02) %
};
\addlegendentry{av.\ runtime $I_{64}^d$}
\addplot[red,mark=diamond,mark size=2pt,mark options={solid},error bars/.cd, y dir=both, y explicit, error bar style={solid}, error mark options={rotate=90,mark size=2,thick}] coordinates {
(2,6.5085e-03) -= (0,6.9969e-04) += (0,5.1055e-03) %
(3,9.9089e-03) -= (0,2.6229e-03) += (0,7.1290e-03) %
(4,9.1002e-03) -= (0,7.1266e-04) += (0,2.1145e-03) %
(5,1.5722e-02) -= (0,2.1947e-03) += (0,3.3301e-03) %
(6,1.7104e-02) -= (0,1.9430e-03) += (0,3.0889e-03) %
(7,2.1670e-02) -= (0,2.3056e-03) += (0,3.1305e-03) %
(8,2.7739e-02) -= (0,4.1194e-03) += (0,2.2224e-03) %
(9,3.4518e-02) -= (0,8.7526e-03) += (0,8.6528e-03) %
(10,3.8288e-02) -= (0,5.2778e-03) += (0,1.0680e-02) %
(11,4.2562e-02) -= (0,6.4631e-03) += (0,5.4914e-03) %
(12,4.8079e-02) -= (0,2.8627e-03) += (0,7.3064e-03) %
(13,5.4084e-02) -= (0,5.9393e-03) += (0,3.9053e-03) %
(14,6.3871e-02) -= (0,5.5002e-03) += (0,3.3959e-03) %
(15,7.0339e-02) -= (0,5.6764e-03) += (0,5.0343e-03) %
(16,8.2777e-02) -= (0,8.8067e-03) += (0,7.7551e-03) %
(17,9.3473e-02) -= (0,1.3907e-02) += (0,8.9992e-03) %
(18,1.0402e-01) -= (0,1.0669e-02) += (0,1.8949e-02) %
(19,1.1756e-01) -= (0,8.9298e-03) += (0,8.1844e-03) %
(20,1.2126e-01) -= (0,1.8912e-02) += (0,2.7412e-02) %
(21,1.3109e-01) -= (0,1.7674e-02) += (0,2.8827e-02) %
(22,1.3807e-01) -= (0,1.2798e-02) += (0,2.0462e-02) %
(23,1.4748e-01) -= (0,1.3165e-02) += (0,2.5099e-02) %
(24,1.6594e-01) -= (0,1.8727e-02) += (0,2.3751e-02) %
(25,1.7447e-01) -= (0,6.4345e-03) += (0,1.3011e-02) %
(26,1.8429e-01) -= (0,7.7604e-03) += (0,7.7570e-03) %
(27,2.0429e-01) -= (0,1.6606e-02) += (0,1.6650e-02) %
(28,2.1469e-01) -= (0,8.9783e-03) += (0,1.2906e-02) %
(29,2.3174e-01) -= (0,1.5910e-02) += (0,1.5239e-02) %
(30,2.4617e-01) -= (0,1.7797e-02) += (0,1.0562e-02) %
(31,2.7252e-01) -= (0,1.3287e-02) += (0,2.0507e-02) %
(32,2.8165e-01) -= (0,2.0855e-02) += (0,2.2477e-02) %
(33,3.0880e-01) -= (0,4.0279e-02) += (0,2.0129e-02) %
(34,3.2419e-01) -= (0,2.9500e-02) += (0,2.4079e-02) %
(35,3.6550e-01) -= (0,3.5475e-02) += (0,9.5855e-02) %
(36,3.6785e-01) -= (0,3.3813e-02) += (0,2.4244e-02) %
(37,3.7955e-01) -= (0,6.6259e-02) += (0,3.4908e-02) %
(38,4.2167e-01) -= (0,6.1569e-02) += (0,6.3846e-02) %
(39,4.1175e-01) -= (0,2.9545e-02) += (0,5.5424e-02) %
(40,4.5177e-01) -= (0,6.8723e-02) += (0,9.0281e-02) %
(41,4.4352e-01) -= (0,5.6011e-02) += (0,7.5418e-02) %
(42,5.1366e-01) -= (0,1.0305e-01) += (0,1.3244e-01) %
(43,5.0555e-01) -= (0,7.2132e-02) += (0,1.1197e-01) %
(44,5.1140e-01) -= (0,6.7306e-02) += (0,1.0730e-01) %
(45,4.9820e-01) -= (0,3.2539e-02) += (0,3.5209e-02) %
(46,5.4616e-01) -= (0,5.4015e-02) += (0,9.5414e-02) %
(47,5.5572e-01) -= (0,4.8839e-02) += (0,6.2446e-02) %
(48,5.9855e-01) -= (0,6.2558e-02) += (0,1.1199e-01) %
(49,6.2760e-01) -= (0,5.4398e-02) += (0,9.9274e-02) %
(50,6.0956e-01) -= (0,3.0531e-02) += (0,7.7435e-02) %
(60,9.0772e-01) -= (0,3.4631e-02) += (0,3.3507e-02) %
(70,1.2857e+00) -= (0,5.9264e-02) += (0,1.7421e-01) %
(80,1.7585e+00) -= (0,2.6680e-01) += (0,4.7071e-01) %
(90,2.1569e+00) -= (0,2.6222e-01) += (0,4.7722e-01) %
(100,2.5107e+00) -= (0,1.7960e-01) += (0,3.6488e-01) %
(150,6.1151e+00) -= (0,2.3935e-01) += (0,3.2548e-01) %
(200,1.1071e+01) -= (0,1.1113e+00) += (0,2.1875e+00) %
(250,1.6357e+01) -= (0,3.1745e-01) += (0,3.1405e-01) %
(300,2.4455e+01) -= (0,7.0802e-01) += (0,5.9579e-01) %
(350,3.5283e+01) -= (0,9.4622e-01) += (0,9.5695e-01) %
(400,4.9079e+01) -= (0,5.3234e+00) += (0,5.3541e+00) %
(450,5.5538e+01) -= (0,1.4119e+00) += (0,3.9098e+00) %
(500,7.0397e+01) -= (0,9.1898e-01) += (0,9.1743e-01) %
(600,1.0699e+02) -= (0,1.7709e+00) += (0,1.4543e+00) %
(700,1.5235e+02) -= (0,4.1290e+00) += (0,4.1462e+00) %
(800,2.1434e+02) -= (0,2.2052e+01) += (0,2.8223e+01) %
(900,2.6267e+02) -= (0,2.2435e+01) += (0,3.9410e+01) %
(1000,3.3779e+02) -= (0,1.9302e+01) += (0,1.2751e+01) %
(1500,8.0329e+02) -= (0,1.3755e+01) += (0,1.0167e+01) %
(2000,1.4421e+03) -= (0,1.7945e+01) += (0,5.7737e+01) %
};
\addlegendentry{av.\ runtime $D(\tilde{I}_{64}^d)$}
\addplot[black, very thick,domain=2:2000, dotted, no markers,samples=200] {x^3/500}; 
\addlegendentry{$d^3/500$}
\addplot[black,very thick,domain=2:2000, dashed, no markers,samples=200] {x^2/2000}; 
\addlegendentry{$d^2/2000$}
\end{loglogaxis}
\end{tikzpicture}
\caption{Runtimes in comparison. Blue: Algorithm~\ref{alg:construct_rand_reco_r1l_I_2} applied to $I_{64}^d$, $K=5$, $T=100$. Red: Algorithm~\ref{alg:construct_rand_reco_r1l_I_2} with modifications explained in Remark~\ref{rem:heuristic_complex_reco} applied to find reconstructing rank-1 lattices for $\tilde{I}_{64}^d$, $K=5$, $T=100$.}\label{fig:anova_comptime}
\end{figure}
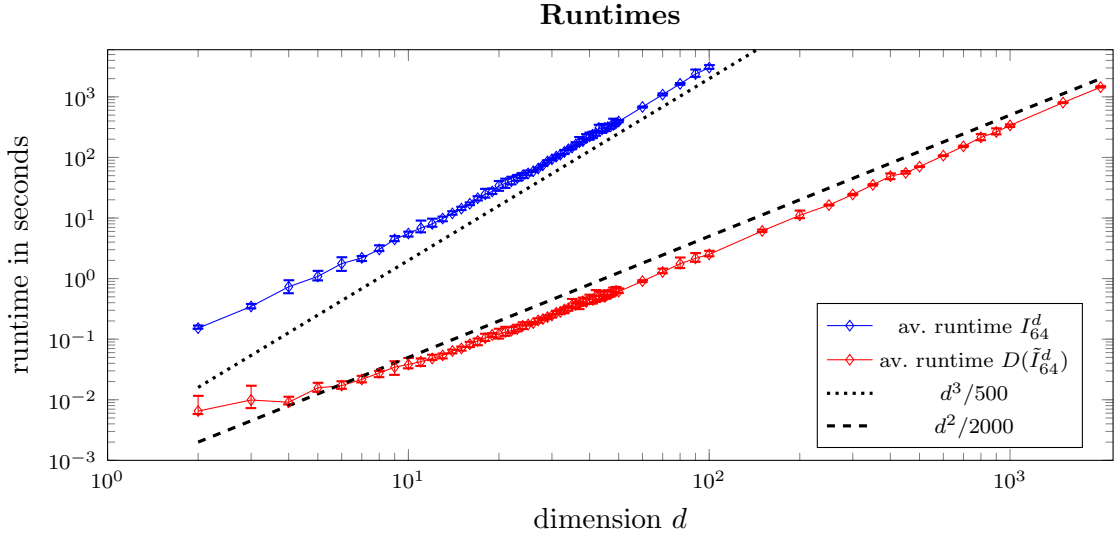
Figure~\ref{fig:anova_comptime} illustrates the different runtimes of Algorithm~\ref{alg:heuristic} applied to $I_{64}^d$ using
Algorithm~\ref{alg:check_exactness1} and Algorithm~\ref{alg:heuristic} applied to $\tilde{I}_{64}^d$ using Algorithm~\ref{alg:check_exactness}.
Obviously, both approaches scales as expected in the dimension $d$. Due to the lower complexity of the second approach,
one can determine exactly integrating rank-1 lattices for $I_{64}^d$ for significantly higher dimensions $d$ when detouring via reconstructing rank-1 lattices for $\tilde{I}_{64}^d$. E.g., in the case $d=350$, it took less than 40 seconds to detect a rank-1 lattice of size $M=980069$ that exactly integrates more than one billion monomials.

At this point, we would like to point out that we observe similar behavior of the runtimes
for fixed dimension $d$ and growing $N$ due to the different dependencies of $|I_{N}^d|$ and $|\tilde{I}_{N}^d|$ on $N$.

\subsection{Reconstructing rank-1 lattices for weighted hyperbolic crosses}

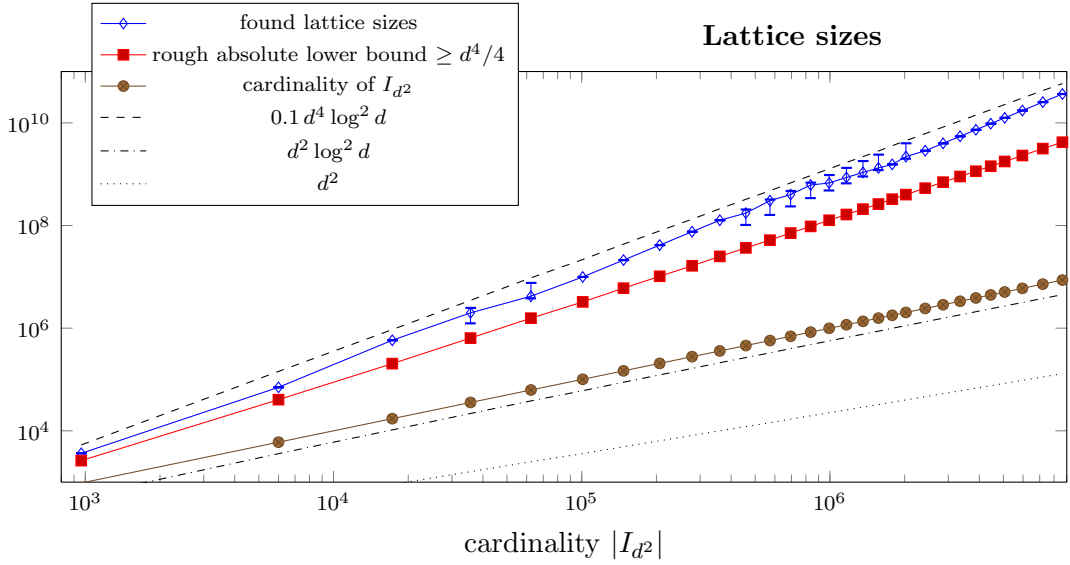
\begin{figure}[tb]
\centering
\begin{tikzpicture}[baseline=(current axis.south)]%
\begin{loglogaxis}[%
    scale only axis,
    xmin=800,xmax=9000000,
    ymin=1000, ymax=10^11,
    title={\rule{6cm}{0pt}\bfseries Lattice sizes},
    x tick label style={align=center,font=\scriptsize},
    y tick label style={font=\scriptsize},
    width=0.85\textwidth,
    height=0.35\textwidth,
    xlabel={cardinality $|I_{d^2}|$},
    legend style={font=\scriptsize, at={(0.03,0.92)},anchor=west, align=left}
    ]
\addplot[blue,mark=diamond,mark size=2pt,mark options={solid},error bars/.cd, y dir=both, y explicit, error bar style={solid}, error mark options={rotate=90,mark size=2,thick}] coordinates {
(963,3.6430e+03) -= (0,0.0000e+00) += (0,0.0000e+00) %
(6003,7.0393e+04) -= (0,0.0000e+00) += (0,0.0000e+00) %
(17251,5.8126e+05) -= (0,0.0000e+00) += (0,0.0000e+00) %
(35655,1.9864e+06) -= (0,7.4489e+05) += (0,4.9660e+05) %
(62449,4.1893e+06) -= (0,3.8085e+05) += (0,3.4276e+06) %
(100951,9.9523e+06) -= (0,0.0000e+00) += (0,0.0000e+00) %
(147365,2.1207e+07) -= (0,0.0000e+00) += (0,0.0000e+00) %
(206139,4.1497e+07) -= (0,0.0000e+00) += (0,0.0000e+00) %
(278439,7.5711e+07) -= (0,0.0000e+00) += (0,0.0000e+00) %
(360067,1.2661e+08) -= (0,0.0000e+00) += (0,0.0000e+00) %
(458327,1.7437e+08) -= (0,7.1799e+07) += (0,3.0771e+07) %
(573453,3.0508e+08) -= (0,1.4451e+08) += (0,1.6057e+07) %
(694751,4.0066e+08) -= (0,1.6498e+08) += (0,7.0705e+07) %
(836491,6.1499e+08) -= (0,2.7333e+08) += (0,6.8332e+07) %
(993025,6.7409e+08) -= (0,1.9260e+08) += (0,2.8890e+08) %
(1165203,8.6182e+08) -= (0,1.9888e+08) += (0,4.6406e+08) %
(1359001,1.0822e+09) -= (0,1.8036e+08) += (0,7.2144e+08) %
(1568741,1.3218e+09) -= (0,1.2016e+08) += (0,1.0815e+09) %
(1782309,1.5511e+09) -= (0,0.0000e+00) += (0,0.0000e+00) %
(2023921,2.2001e+09) -= (0,2.0001e+08) += (0,1.8001e+09) %
(2420663,2.8611e+09) -= (0,0.0000e+00) += (0,0.0000e+00) %
(2856721,3.9848e+09) -= (0,0.0000e+00) += (0,0.0000e+00) %
(3344661,5.4623e+09) -= (0,0.0000e+00) += (0,0.0000e+00) %
(3873041,7.3244e+09) -= (0,0.0000e+00) += (0,0.0000e+00) %
(4445339,9.6489e+09) -= (0,0.0000e+00) += (0,0.0000e+00) %
(5061373,1.2509e+10) -= (0,0.0000e+00) += (0,0.0000e+00) %
(5962465,1.7359e+10) -= (0,0.0000e+00) += (0,0.0000e+00) %
(7205821,2.5353e+10) -= (0,0.0000e+00) += (0,0.0000e+00) %
(8644653,3.6489e+10) -= (0,0.0000e+00) += (0,0.0000e+00) %
};
\addlegendentry{found lattice sizes}
\addplot coordinates{
(963,2.6260e+03) %
(6003,4.0501e+04) %
(17251,2.0363e+05) %
(35655,6.4200e+05) %
(62449,1.5656e+06) %
(100951,3.2445e+06) %
(147365,6.0086e+06) %
(206139,1.0248e+07) %
(278439,1.6413e+07) %
(360067,2.5013e+07) %
(458327,3.6618e+07) %
(573453,5.1858e+07) %
(694751,7.1424e+07) %
(836491,9.6065e+07) %
(993025,1.2659e+08) %
(1165203,1.6387e+08) %
(1359001,2.0884e+08) %
(1568741,2.6248e+08) %
(1782309,3.2585e+08) %
(2023921,4.0005e+08) %
(2420663,5.3423e+08) %
(2856721,6.9967e+08) %
(3344661,9.0081e+08) %
(3873041,1.1425e+09) %
(4445339,1.4299e+09) %
(5061373,1.7683e+09) %
(5962465,2.3089e+09) %
(7205821,3.1487e+09) %
(8644653,4.1992e+09) %
};
\addlegendentry{rough absolute lower bound $\ge d^4/4$}
\addplot coordinates{
(963,963) %
(6003,6003) %
(17251,17251) %
(35655,35655) %
(62449,62449) %
(100951,100951) %
(147365,147365) %
(206139,206139) %
(278439,278439) %
(360067,360067) %
(458327,458327) %
(573453,573453) %
(694751,694751) %
(836491,836491) %
(993025,993025) %
(1165203,1165203) %
(1359001,1359001) %
(1568741,1568741) %
(1782309,1782309) %
(2023921,2023921) %
(2420663,2420663) %
(2856721,2856721) %
(3344661,3344661) %
(3873041,3873041) %
(4445339,4445339) %
(5061373,5061373) %
(5962465,5962465) %
(7205821,7205821) %
(8644653,8644653) %
};
\addlegendentry{cardinality of $I_{d^2}$}
\addplot[black, dashed, no markers] coordinates{
(963,5.3019e+03) %
(6003,1.4359e+05) %
(17251,9.3702e+05) %
(35655,3.4836e+06) %
(62449,9.5650e+06) %
(100951,2.1726e+07) %
(147365,4.3337e+07) %
(206139,7.8652e+07) %
(278439,1.3285e+08) %
(360067,2.1208e+08) %
(458327,3.2349e+08) %
(573453,4.7527e+08) %
(694751,6.7669e+08) %
(836491,9.3811e+08) %
(993025,1.2710e+09) %
(1165203,1.6880e+09) %
(1359001,2.2030e+09) %
(1568741,2.8309e+09) %
(1782309,3.5879e+09) %
(2023921,4.4915e+09) %
(2420663,6.1632e+09) %
(2856721,8.2757e+09) %
(3344661,1.0904e+10) %
(3873041,1.4130e+10) %
(4445339,1.8043e+10) %
(5061373,2.2737e+10) %
(5962465,3.0391e+10) %
(7205821,4.2574e+10) %
(8644653,5.8192e+10) %
};
\addlegendentry{$0.1\,d^4\log^2{d}$}
\addplot[black, dashdotted, no markers] coordinates{
(963,5.3019e+02) %
(6003,3.5898e+03) %
(17251,1.0411e+04) %
(35655,2.1773e+04) %
(62449,3.8260e+04) %
(100951,6.0349e+04) %
(147365,8.8444e+04) %
(206139,1.2289e+05) %
(278439,1.6401e+05) %
(360067,2.1208e+05) %
(458327,2.6734e+05) %
(573453,3.3005e+05) %
(694751,4.0041e+05) %
(836491,4.7863e+05) %
(993025,5.6490e+05) %
(1165203,6.5939e+05) %
(1359001,7.6228e+05) %
(1568741,8.7372e+05) %
(1782309,9.9388e+05) %
(2023921,1.1229e+06) %
(2420663,1.3333e+06) %
(2856721,1.5644e+06) %
(3344661,1.8166e+06) %
(3873041,2.0903e+06) %
(4445339,2.3858e+06) %
(5061373,2.7036e+06) %
(5962465,3.1625e+06) %
(7205821,3.7937e+06) %
(8644653,4.4902e+06) %
};
\addlegendentry{$d^2\log^2{d}$}
\addplot[black, dotted, no markers] coordinates{
(963,1.0000e+02) %
(6003,4.0000e+02) %
(17251,9.0000e+02) %
(35655,1.6000e+03) %
(62449,2.5000e+03) %
(100951,3.6000e+03) %
(147365,4.9000e+03) %
(206139,6.4000e+03) %
(278439,8.1000e+03) %
(360067,1.0000e+04) %
(458327,1.2100e+04) %
(573453,1.4400e+04) %
(694751,1.6900e+04) %
(836491,1.9600e+04) %
(993025,2.2500e+04) %
(1165203,2.5600e+04) %
(1359001,2.8900e+04) %
(1568741,3.2400e+04) %
(1782309,3.6100e+04) %
(2023921,4.0000e+04) %
(2420663,4.6225e+04) %
(2856721,5.2900e+04) %
(3344661,6.0025e+04) %
(3873041,6.7600e+04) %
(4445339,7.5625e+04) %
(5061373,8.4100e+04) %
(5962465,9.6100e+04) %
(7205821,1.1223e+05) %
(8644653,1.2960e+05) %
};
\addlegendentry{$d^2$}
\end{loglogaxis}
\end{tikzpicture}
\caption{Cardinalities of weighted hyperbolic cross frequency sets $I_{d^2}$ and corresponding reconstructing rank-1 lattice sizes determined using Algorithm~\ref{alg:construct_rand_reco_r1l_I_2} taking Remark~\ref{rem:heuristic_complex_reco} into account.}
\label{fig:j2_M}
\end{figure}

We consider weighted hyperbolic cross frequency sets, which are widely used examples in the literature on rank-1 lattices. In particular, the product weights $\omega(\boldk):=\prod_{j=1}^\infty\max(1,\gamma_j^{-1}|k_j|)$ using the weight vector
$\gamma_j=j^{-2}$  were already used for numerical tests in, e.g., \cite{Kuo03}. 
In short words, the weight function $\omega$ specifies the importance of the frequencies within the whole $\Z^\infty$, cf.\ e.g.~\cite{KaVo19, CoKuNuSl19} for more details on different types of weights and corresponding function spaces. The larger the weight $\omega(\boldk)$ the lower is the importance of the frequency $\boldk$. Accordingly, we collect the more important frequencies in a frequency set 
$$
I_N:=\{\boldk\in\Z^\infty\colon\omega(\boldk)\le N\}.
$$
Considering the frequencies $\boldk\in I_N$ in detail, one observes that they have nonzero components only in the first few dimensions which is caused by the decay of the weight vector $\boldgamma$ that implies that $I_N\subset\Z^{\floor{\sqrt{N}}}\times\{0\}^\infty$ and, thus, $I_N$ can be considered as $\floor{\sqrt{N}}$-dimensional frequency set.
For that reason, we choose $N=d^2$ and treat $I_{d^2}$ as a $d$-dimensional frequency set in the following. Due to the decay of the weight vector $\boldgamma$, we know that for each $q>2$ there exists $c_q<\infty$ such that $|I_{d^2}|\le c_{q} d^q$. Moreover, using a similar argumentation as in \cite[Theorem 3.5]{KaKuPo10} yields that a reconstructing rank-1 lattice for $I_{d^2}$ necessarily requires a lattice size $M$ of at least
$(d^2+1)(\floor{d^2/4}+1)>d^4/4$. We just applied Algorithm~\ref{alg:construct_rand_reco_r1l_I_2}
ten times and obtained reconstructing rank-1 lattices with average lattice sizes depicted as diamonds in Figure~\ref{fig:j2_M}.
In addition, we illustrate the whole range of the computed lattice sizes as bars. In more detail, the ten different applications of Algorithm~\ref{alg:construct_rand_reco_r1l_I_2} to a fixed frequency set $I_{d^2}$ determined reconstructing rank-1 lattices of sizes $M$ that differ only by a factor of at most approximately two.
In fact, this observation indicates
that the outer while loop in Algorithm~\ref{alg:construct_rand_reco_r1l_I_2} terminated at only two different -- and successive -- iteration levels for fixed frequency set $I_{d^2}$. Accordingly, the output of Algorithm~\ref{alg:construct_rand_reco_r1l_I_2} seems to be highly reliable with respect to the order of magnitude of the determined rank-1 lattice size.

\begin{figure}[tb]
\centering
\begin{tikzpicture}[baseline=(current axis.south)]%
\begin{semilogyaxis}[%
    scale only axis,
    xmin=10,xmax=360,
    ymin=1, ymax=5000,
    title={\bfseries Oversampling factors},
    xtick={10,50,100,150,200,250,300,350},
    x tick label style={align=center,font=\scriptsize},
    y tick label style={font=\scriptsize},
    xlabel={dimension $d$},
    width=0.85\textwidth,
    height=0.35\textwidth,
    legend pos=south east,
    legend style={font=\scriptsize},
    ]
\addplot[blue,mark=diamond,mark size=2pt,mark options={solid},error bars/.cd, y dir=both, y explicit, error bar style={solid}, error mark options={rotate=90,mark size=2,thick}] coordinates {
(10,3.7830e+00) -= (0,0.0000e+00) += (0,0.0000e+00) %
(20,1.1726e+01) -= (0,0.0000e+00) += (0,0.0000e+00) %
(30,3.3694e+01) -= (0,0.0000e+00) += (0,0.0000e+00) %
(40,5.5711e+01) -= (0,2.0892e+01) += (0,1.3928e+01) %
(50,6.7084e+01) -= (0,6.0985e+00) += (0,5.4887e+01) %
(60,9.8585e+01) -= (0,0.0000e+00) += (0,0.0000e+00) %
(70,1.4391e+02) -= (0,0.0000e+00) += (0,0.0000e+00) %
(80,2.0131e+02) -= (0,0.0000e+00) += (0,0.0000e+00) %
(90,2.7191e+02) -= (0,0.0000e+00) += (0,0.0000e+00) %
(100,3.5163e+02) -= (0,0.0000e+00) += (0,0.0000e+00) %
(110,3.8045e+02) -= (0,1.5665e+02) += (0,6.7138e+01) %
(120,5.3201e+02) -= (0,2.5201e+02) += (0,2.8001e+01) %
(130,5.7670e+02) -= (0,2.3746e+02) += (0,1.0177e+02) %
(140,7.3520e+02) -= (0,3.2675e+02) += (0,8.1689e+01) %
(150,6.7883e+02) -= (0,1.9395e+02) += (0,2.9093e+02) %
(160,7.3963e+02) -= (0,1.7068e+02) += (0,3.9826e+02) %
(170,7.9629e+02) -= (0,1.3271e+02) += (0,5.3086e+02) %
(180,8.4259e+02) -= (0,7.6599e+01) += (0,6.8939e+02) %
(190,8.7027e+02) -= (0,0.0000e+00) += (0,0.0000e+00) %
(200,1.0871e+03) -= (0,9.8824e+01) += (0,8.8942e+02) %
(215,1.1820e+03) -= (0,0.0000e+00) += (0,0.0000e+00) %
(230,1.3949e+03) -= (0,0.0000e+00) += (0,0.0000e+00) %
(245,1.6331e+03) -= (0,0.0000e+00) += (0,0.0000e+00) %
(260,1.8911e+03) -= (0,0.0000e+00) += (0,0.0000e+00) %
(275,2.1706e+03) -= (0,0.0000e+00) += (0,0.0000e+00) %
(290,2.4714e+03) -= (0,0.0000e+00) += (0,0.0000e+00) %
(310,2.9114e+03) -= (0,0.0000e+00) += (0,0.0000e+00) %
(335,3.5185e+03) -= (0,0.0000e+00) += (0,0.0000e+00) %
(360,4.2210e+03) -= (0,0.0000e+00) += (0,0.0000e+00) %
};
\addlegendentry{found lattices}
\addplot coordinates{
(10,2.7269e+00) %
(20,6.7468e+00) %
(30,1.1804e+01) %
(40,1.8006e+01) %
(50,2.5070e+01) %
(60,3.2139e+01) %
(70,4.0774e+01) %
(80,4.9714e+01) %
(90,5.8945e+01) %
(100,6.9466e+01) %
(110,7.9894e+01) %
(120,9.0431e+01) %
(130,1.0280e+02) %
(140,1.1484e+02) %
(150,1.2748e+02) %
(160,1.4064e+02) %
(170,1.5367e+02) %
(180,1.6732e+02) %
(190,1.8282e+02) %
(200,1.9766e+02) %
(215,2.2070e+02) %
(230,2.4492e+02) %
(245,2.6933e+02) %
(260,2.9499e+02) %
(275,3.2165e+02) %
(290,3.4937e+02) %
(310,3.8724e+02) %
(335,4.3697e+02) %
(360,4.8576e+02) %
};
\addlegendentry{rough absolute lower bound}
\end{semilogyaxis}
\end{tikzpicture}
\caption{Oversampling factors $M/|I_{d^2}|$ of reconstructing rank-1 lattices $\Lambda(\boldz,M)$ constructed by Algorithm~\ref{alg:construct_rand_reco_r1l_I_2} taking Remark~\ref{rem:heuristic_complex_reco} into account and
simple theoretical lower bound on the oversampling factor.}
\label{fig:j2_M2}
\end{figure}
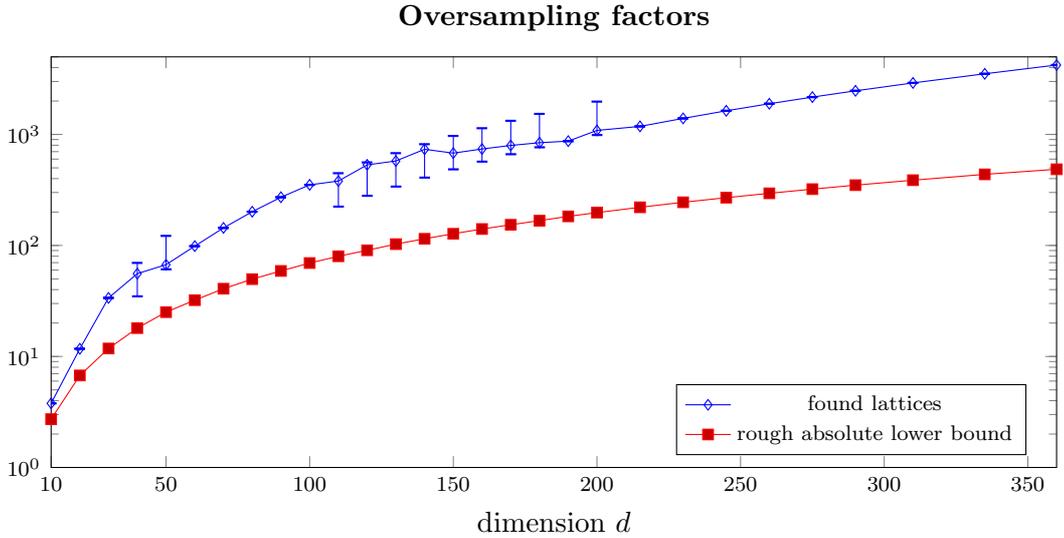

A further observation of this numerical test is that the determined rank-1 lattices are of sizes that are relatively close to the rough lower bound $(d^2+1)(\floor{d^2/4}+1)>d^4/4$. In particular in our numerical tests, the lattice sizes are less than twenty times this lower bound, cf. Figure~\ref{fig:j2_M2}, and much less than the starting point, cf. Remark~\ref{rem:heuristic_complex_reco}, which is approximately $|I_{d^2}|^2$.

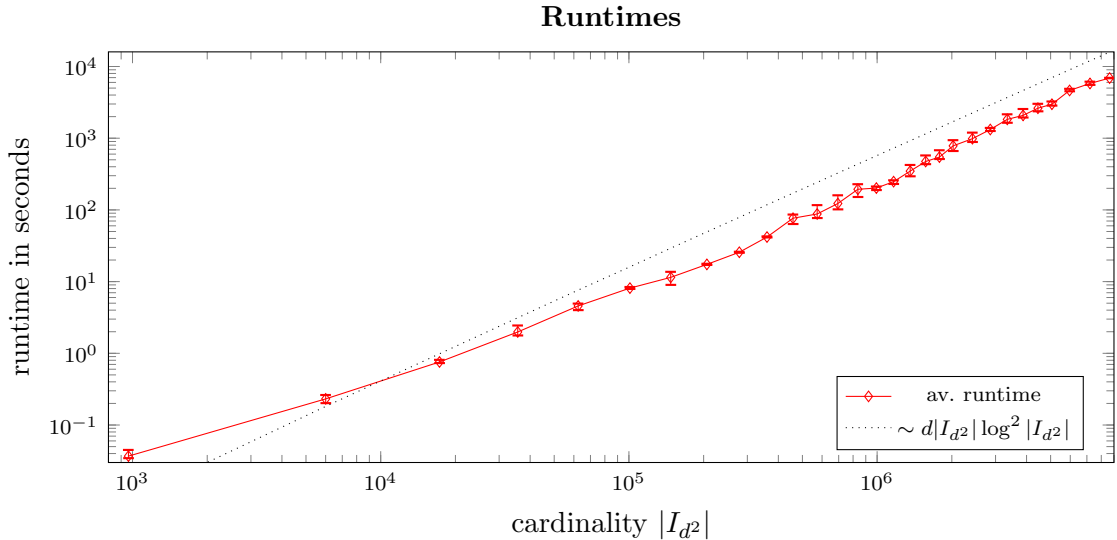
\begin{figure}[tb]
\centering
\begin{tikzpicture}[baseline=(current axis.south)]%
\begin{loglogaxis}[%
    scale only axis,
    xmin=800,xmax=9000000,
    ymin=0.03, ymax=16000,
    title={\bfseries Runtimes},
    xlabel={cardinality $|I_{d^2}|$},
    ylabel={runtime in seconds},
    x tick label style={align=center,font=\scriptsize},
    y tick label style={font=\scriptsize},
    width=0.85\textwidth,
    height=0.35\textwidth,
    legend pos=south east,
    legend style={font=\scriptsize}
    ]
  \addplot[red,mark=diamond,mark size=2pt,mark options={solid},error bars/.cd, y dir=both, y explicit, error bar style={solid}, error mark options={rotate=90,mark size=2,thick}] coordinates {
(963,3.6820e-02) -= (0,2.3277e-03) += (0,8.1122e-03) %
(6003,2.3145e-01) -= (0,2.9539e-02) += (0,3.0772e-02) %
(17251,7.5704e-01) -= (0,2.6318e-02) += (0,5.3660e-02) %
(35655,1.9845e+00) -= (0,2.1037e-01) += (0,4.5934e-01) %
(62449,4.5777e+00) -= (0,5.7365e-01) += (0,3.5989e-01) %
(100951,8.1144e+00) -= (0,2.1238e-01) += (0,3.0070e-01) %
(147365,1.1454e+01) -= (0,2.4492e+00) += (0,2.2399e+00) %
(206139,1.7385e+01) -= (0,3.0148e-01) += (0,3.5919e-01) %
(278439,2.5677e+01) -= (0,4.3844e-01) += (0,4.1796e-01) %
(360067,4.1829e+01) -= (0,5.3677e-01) += (0,9.1596e-01) %
(458327,7.6608e+01) -= (0,1.3067e+01) += (0,9.5418e+00) %
(573453,8.7744e+01) -= (0,1.0731e+01) += (0,2.8845e+01) %
(694751,1.2281e+02) -= (0,2.1065e+01) += (0,3.7071e+01) %
(836491,1.9405e+02) -= (0,4.2656e+01) += (0,3.4274e+01) %
(993025,2.0163e+02) -= (0,1.2508e+01) += (0,1.0235e+01) %
(1165203,2.4773e+02) -= (0,1.7695e+01) += (0,1.1097e+01) %
(1359001,3.4742e+02) -= (0,5.3087e+01) += (0,7.6150e+01) %
(1568741,4.7348e+02) -= (0,3.8889e+01) += (0,1.0144e+02) %
(1782309,5.4524e+02) -= (0,3.2921e+01) += (0,1.3398e+02) %
(2023921,7.8227e+02) -= (0,1.1958e+02) += (0,1.5679e+02) %
(2420663,9.8484e+02) -= (0,9.9054e+01) += (0,2.1129e+02) %
(2856721,1.3221e+03) -= (0,6.2031e+01) += (0,6.9702e+01) %
(3344661,1.8310e+03) -= (0,1.8931e+02) += (0,3.2398e+02) %
(3873041,2.0810e+03) -= (0,1.3775e+02) += (0,4.7392e+02) %
(4445339,2.6207e+03) -= (0,2.4914e+02) += (0,3.9741e+02) %
(5061373,2.9792e+03) -= (0,1.3218e+02) += (0,2.7282e+02) %
(5962465,4.6355e+03) -= (0,1.2394e+02) += (0,2.4983e+02) %
(7205821,5.8003e+03) -= (0,2.5980e+02) += (0,3.4463e+02) %
(8644653,6.8990e+03) -= (0,4.4115e+01) += (0,5.8528e+01) %
};
\addlegendentry{av.\ runtime}
\addplot[black, dotted, no markers] coordinates{
(963,9.0903e-04) %
(963,9.0903e-03) %
(6003,1.8175e-01) %
(17251,9.8509e-01) %
(35655,3.1338e+00) %
(62449,7.6143e+00) %
(100951,1.6083e+01) %
(147365,2.9219e+01) %
(206139,4.9383e+01) %
(278439,7.8775e+01) %
(360067,1.1788e+02) %
(458327,1.7133e+02) %
(573453,2.4197e+02) %
(694751,3.2684e+02) %
(836491,4.3557e+02) %
(993025,5.6804e+02) %
(1165203,7.2752e+02) %
(1359001,9.2153e+02) %
(1568741,1.1493e+03) %
(1782309,1.4031e+03) %
(2023921,1.7069e+03) %
(2420663,2.2491e+03) %
(2856721,2.9038e+03) %
(3344661,3.6987e+03) %
(3873041,4.6345e+03) %
(4445339,5.7289e+03) %
(5061373,6.9957e+03) %
(5962465,8.9975e+03) %
(7205821,1.2038e+04) %
(8644653,1.5879e+04) %
};
\addlegendentry{$\sim d|I_{d^2}|\log^2|I_{d^2}|$}
\end{loglogaxis}
\end{tikzpicture}
\caption{Runtimes of Algorithm~\ref{alg:construct_rand_reco_r1l_I_2} determining reconstructing rank-1 lattices for $I_{d^2}$ of reasonable lattice sizes, $K=5$, $T=100$.}\label{fig:hc_comptime}
\end{figure}

Figure~\ref{fig:hc_comptime} depicts the runtimes needed by Algorithm~\ref{alg:construct_rand_reco_r1l_I_2} taking Remark~\ref{rem:heuristic_complex_reco} into account for determining reconstructing rank-1 lattices
for $I_{d^2}$. In addition, we added a dotted line representing the asymptotic behaviour of the computational complexity of Algorithm~\ref{alg:construct_rand_reco_r1l_I_2} for large enough and fixed parameters $K$ and $T$. The slopes of both lines are very close, which indicates that the complexity of our implementation of Algorithm~\ref{alg:construct_rand_reco_r1l_I_2} is indeed in $\OO{d|I|\log^2|I|}$ for fixed $T$ and $K$, cf.~\eqref{eq:compl_alg5_reco}.

\subsection*{Acknowledgment}
The author was supported by the Deutsche Forschungsgemeinschaft (DFG – German Research Foundation, project number 380648269).

\scriptsize
\setlength{\bibsep}{2.5pt}

\end{document}